%% file: main_arxiv.tex
\documentclass{article}
\usepackage{times}
\usepackage{soul}
\usepackage{url}
\usepackage[utf8]{inputenc}
\usepackage[small]{caption}
\usepackage{graphicx}
\usepackage{amsmath}
\usepackage{amsthm}
\usepackage{booktabs}
\usepackage{algorithm}
\usepackage{algorithmic}
\usepackage[algo2e]{algorithm2e}
\usepackage{verbatim}
\usepackage{graphicx}
\usepackage{amsmath}
\usepackage{bm}
\usepackage{booktabs}
\usepackage{subcaption}
\usepackage{textcomp}
\usepackage{multirow}
\usepackage{color}

\usepackage{wrapfig}
\usepackage{color,url}
\usepackage{multicol}
\usepackage{multirow}
\usepackage{microtype}
\usepackage{adjustbox}
\usepackage{graphics}
\usepackage{amsmath}
\usepackage{amssymb}
\usepackage{mathtools}
\usepackage{amsthm}
\usepackage{subcaption}
\usepackage{booktabs} 

\usepackage{pifont}

\usepackage{hyperref}
\usepackage[capitalize,noabbrev]{cleveref}
\theoremstyle{plain}
\newtheorem{theorem}{Theorem}

\newtheorem{lemma}{Lemma}

\theoremstyle{definition}
\newtheorem{definition}{Definition}
\newtheorem{assumption}{Assumption}
\theoremstyle{remark}

\usepackage[textsize=tiny]{todonotes}
\usepackage{yhmath}
\usepackage{color}
\newcommand{\E}{\mathbb{E}}
\usepackage{amsmath,epsfig,bm, amssymb,graphicx,algorithm,algorithmic,color}

\input{mathdef}

\advance\oddsidemargin-1in
\author{
Dinesh Singh$^1$, Hardik Tankariya$^2$, Makoto~Yamada$^{1,2}$\\
$^1$RIKEN AIP,
$^2$Kyoto University,
}

\title{Nys-Newton: Nystr\"om-Approximated Curvature\\ for Stochastic Optimization}

\begin{document}
	
	\maketitle
	
	\begin{abstract}
		Second-order optimization methods are among the most widely used optimization approaches for convex optimization problems, and have recently been used to optimize non-convex optimization problems such as deep learning models. The widely used second-order optimization methods such as quasi-Newton methods generally provide curvature information by approximating the Hessian using the secant equation. However, the secant equation becomes insipid in approximating the Newton step owing to its use of the first-order derivatives.
		In this study, we propose an approximate Newton sketch-based stochastic optimization algorithm for large-scale empirical risk minimization. Specifically, we compute a partial column Hessian of size ($d\times m$) with $m\ll d$ randomly selected variables, then use the \emph{Nystr\"om method} to better approximate the full Hessian matrix. To further reduce the computational complexity per iteration, we directly compute the update step ($\Delta\boldw$) without computing and storing the full Hessian or its inverse. We then  integrate our approximated Hessian with stochastic gradient descent and stochastic variance-reduced gradient methods. 
		The results of numerical experiments on both convex and non-convex functions show that the proposed approach was able to obtain a better approximation of Newton\textquotesingle s method, exhibiting performance competitive with that of state-of-the-art first-order and stochastic quasi-Newton methods. Furthermore, we provide a theoretical convergence analysis for convex functions.
	\end{abstract}
	
	\section{Introduction}
	The  problem  of  the optimization of various function is among  the  most  critical  and  popular topics in machine learning and mathematical optimization. Let $\{(\boldx_i,y_i)\}_{i=1}^{n}$ be given $n$ training samples, where $\boldx_i \in \mathbb{R}^d $ and $y_i \in \{-1,1\}$, and $f$ be an objective function defined as follows.
	\begin{equation}\label{objective_fun}
		\min_{\boldw \in \mathbb{R}^d} f(\boldw) = \frac{1}{n} \sum_{i=1}^{n} f(\boldw;\boldx_i,y_i) =\frac{1}{n}\sum_{i=1}^{n} f_i(\boldw),   
	\end{equation}
	where $f_i(\boldw) = f(\boldw;\boldx_i,y_i):\mathbb{R}^d \to \mathbb{R}$ is the loss function.
	Well-known convex loss functions include logistic loss  $f_i(\boldw) = \log(1+\exp{(-y_i\boldx_i^\top\boldw)})$.
	
	To optimize~\eqref{objective_fun}, first-order optimization methods such as stochastic gradient descent (SGD)~\cite{Robbins1951sgd}, AdaGrad~\cite{Duchi2011adaptive}, stochastic variance-reduced gradient (SVRG)~\cite{Johnson2013accelerating}, SAGA~\cite{Defazio2014saga}, Adam~\cite{Kingma2015adam}, and the stochastic recursive gradient algorithm (SARAH), ~\cite{Nguyen2017sarah}, possibly augmented with momentum, are preferred for large-scale optimization problems owing to their more affordable computational costs, which are linear in dimensions per epoch $O(nd)$. However, the convergence of the first-order methods is notably slow, and they are sensitive to hyperparameter choices and ineffective for ill-conditioned problems.
	
	In contrast, Newton\textquotesingle s method does not depend on the parameters of specific problems and requires only minimal hyperparameter tuning for self-concordant functions, such as $\ell_2$-regularized logistic regression. However, Newton\textquotesingle s method involves a computational complexity of $\Omega(nd^2+d^{2.37})$~\cite{AgarwalBH17jmlr} per iteration and thus is not suitable for large-scale settings. To reduce this computational complexity, the sub-sampled Newton’s method and random projection (or sketching) are commonly used to reduce the dimensionality of the problem and solve it in a lower-dimensional subspace. The sub-sampled Newton method performs well for large-scale but relatively low-dimensional problems by computing the Hessian matrix on a relatively small sample. However, the approximated Hessian can deviate significantly from the full Hessian in high-dimensional problems. Randomized algorithms~\cite{Lacotte2021icml, pilanci2017newton} estimate the Hessian in Newton’s method using a random embedding matrix $\boldS\in\mathbb{R}^{m\times n}$ $\boldH_\boldS(\boldw):=(\nabla^2f(\boldw)^{\frac{1}{2}})^\top\boldS^\top\boldS\nabla^2f(\boldw)^{\frac{1}{2}}$. Specifically, their approximation used the square root of the generalized Gauss-Newton (GGN) matrix as a low-rank approximation instead of deriving it from actual curvature information, whereas $\boldS$ is a random projection matrix of size ($m\times n$). Thus, the approximation may deviate significantly from the actual Hessian. Other sketch matrix-based approximations to reduce the computational complexity of the Newton step calculation have been proposed with deterministic updates which are therefore unsuitable for large-scale stochastic settings where $n$ is significantly larger.
	
	The limited-memory Broyden-Fletcher-Goldfarb-Shanno (LBFGS) algorithm~\cite{Liu1989limited} is a widely used stochastic quasi-Newton method. More specifically, it estimates the Hessian inverse using the past difference of gradients and updates. The online BFGS (oBFGS)~\cite{Schraudolph2007stochastic} method is a stochastic version of regularized BFGS and L-BFGS with gradient descent. \cite{Kolte2015accelerating} proposed two variants of a stochastic quasi-Newton method incorporating a variance-reduced gradient. The first variant used a sub-sampled Hessian with singular value thresholding, which is numerically weaker. The second variant used the LBFGS method to approximate the Hessian inverse. Because these quasi-Newton methods approximate the Hessian inverse using first-order gradients, they are often unable to provide significant curvature information. The stochastic quasi-Newton method (SQN)~\cite{Byrd2016stochastic} used the Hessian vector product computed on a subset of each mini-batch instead of approximating the Hessian inverse from the difference between the current and previous gradients, as in LBFGS. SVRG-SQN~\cite{Moritz2016linearly} also incorporated variance-reduced gradients. Hence, both SQN and SVRG-SQN perform well with low-dimensional datasets, whereas their computational time cost is drastically increased for high-dimensional datasets. 
	
	In this study, we overcome the limitations of the random projection/sketch-based and quasi-Newton methods to further enhance the curvature information used to approximate the Newton step through the $k$-rank approximation of the Hessian matrix. In contrast to quasi-Newton methods that estimate the curvature information using first-order gradients only, we use Nystr\"om approximation on the partial Hessian matrix ($d\times m$) constructed for the $m\ll d$ randomly selected columns only to estimate the Hessian matrix. Thus, our approximation is also more suitable for high-dimensional cases. Note that the proposed method can be represented in terms of the sketch method; more details can be found in Section ~\ref{sec:Newton-sketch}. We used our approximation with stochastic gradient descent (SGD) and stochastic variance-reduced gradient methods (SVRG), and theoretically prove their convergence for convex functions. In addition, we evaluated the performance of the proposed approach on several large-scale datasets with a wide range of dimensions for both convex and nonconvex functions. The experimental results show that the proposed methods consistently performed better than or competitively with existing methods, whereas the behavior of existing methods depended significantly on the problem.
	
	\noindent
	{\bf Contribution:} The contribution of this study are summarized as follows.
	\begin{itemize} 
		\item We propose Nystr\"{o}m-approximated Newton sketch-based methods to perform stochastic optimization.
		\item We show theoretically that the proposed optimization technique can achieve linear and linear-quadratic convergence in a stochastic convex optimization setup for SVRG and SGD, respectively.
		\item We empirically show that the performance of the proposed method compared favorably with that of existing methods in training deep learning models. 
	\end{itemize}
	
	\section{Proposed method}
	In this section, we propose Nystr\"om stochastic gradient descent (Nystr\"om-SGD) and stochastic variance-reduced gradient (Nystr\"om-SVRG) methods. The Nystr\"om method is a widely used kernel technique.  However, surprisingly, Nystr\"om method for second-order stochastic optimization has not been well studied.
	
	\subsection{Nystr\"om SGD and SVRG methods} In the proposed method, to incorporate a better approximation of the curvature information, we use a low-rank approximation of the Hessian matrix. For example, the Hessian matrix of a binary logistic regression is positive semi-definite.
	However, if we use the low-rank approximation of the Hessian matrix, the approximated Hessian cannot be invertible. Moreover, the computation of the approximated Hessian is computationally expensive if computed in each step. Thus, we propose the following regularized variant as a Newton-like update.
	\begin{align}\label{eq:update_rule} 
		\boldw_t=\boldw_{t-1}-\eta(\boldN_{\tau}+\rho \boldI)^{-1}\boldv_{t-1},
	\end{align}
	where $\boldv_t$ denotes the appropriate stochastic gradient, $\boldN_\tau$ is the rank $k$-approximated Hessian matrix at epoch $\tau$, $\eta \geq 0$ is a hyperparameter, and $\rho > 0$ is a regularization parameter. In this study, we denote $\boldB_{\tau} = (\boldN_{\tau}+\rho \boldI)^{-1}$. As $(\boldN_{\tau}+\rho \boldI)$ is positive definite, its matrix inverse can be computed. 
	
	The key challenge is to estimate a low-rank Hessian matrix for each step. To this end, we employed the Nystr\"om method and computed the approximated Hessian at each \emph{epoch} $\tau$.
	\begin{definition}[Nystr\"om approximation]
		Let $\boldH\in\mathbb{R}^{d\times d}$ be a symmetric positive semi-definite matrix. Then, the 
		$k$-rank $(k\ll d)$ approximation $\boldN_k$ of the matrix $\boldH$ is given by 
		\begin{align}
			\boldN_k & = \boldC\boldM^\dagger_k\boldC^\top= \boldZ\boldZ^\top,\label{def:Nystrom}
		\end{align}
		where $\boldZ = \boldC \boldU_k\boldSigma_k^{-1/2} \in \mathbbR^{d \times k}$, and $\boldC \in \mathbb{R}^{d\times m}$ is a matrix consisting of $m$ columns $(m\ll d)$ of $\boldH$, 
		and $\boldM_k$ is the best $k$-rank approximation of $\boldM$, which is formed by the 
		intersection between those $m$ columns of $\boldH$ and the corresponding $m$ rows of $\boldH$,
		and $\boldM^\dagger_k$ is the pseudo-inverse of $\boldM_k$, and the rank of $\boldM $ is $k\leq m$. Note that the number of columns $m$ is a hyperparameter. To obtain the best $k$ rank approximation, it can be computed using the singular value decomposition (SVD) of $\boldM_k$ as $\boldM_k = \boldU_k\boldSigma_k \boldU_k^\top$, where  $\boldU_k \in \mathbbR^{m \times k}$ are singular vectors and  $\boldSigma_k \in \mathbbR^{k \times k}$ are singular values. The pseudo-inverse can be computed as $\boldM^\dagger_k = \boldU_k \boldSigma_k^{-1}\boldU_k^\top$.
	\end{definition}
	At each epoch, we uniformly sample $m$ columns as $\boldOmega\subset\{1,\cdots,d\}$, where $m\ll d$, unlike all columns, as in~\cite{Byrd2016stochastic,yang2020structured}. Then, we compute a partial Hessian matrix $\boldC \in \mathbbR^{d \times m}$
	as the Jacobian of $\nabla f(\boldw)$ with respect to $\boldw_{\boldOmega}$. Then, by using Nystr\"om method in~\eqref{def:Nystrom}, the $k$-rank approximate Hessian can be computed as
	\begin{align*}
		\boldN_\tau = \boldZ\boldZ^\top~\text{and}~ \boldZ = \boldC \boldU_k \boldSigma_k^{-1/2}. 
	\end{align*} 
	To efficiently compute the inverse in~\eqref{eq:update_rule}, we compute
	$\boldB_{\tau}\boldv_{t-1}$ in~\eqref{eq:update_rule} using the inversion lemma as
	\begin{align}\label{eq:vect_mult}
		\hspace{-.05in}\boldB_\tau\boldv_{t-1}\!=\!(\boldN_\tau\!+\!\rho\boldI)^{-1}\boldv_{t-1}\!=\! {\frac{1}{\rho}}\boldv_{t-1} \!-\! {\boldQ_{\tau}\boldZ_{\tau}^\top }\boldv_{t-1},
	\end{align}
	where $\boldQ_{\tau} = \frac{1}{\rho^2}\boldZ_{\tau}(\boldI_k +\frac{1}{\rho} \boldZ_{\tau}^\top\boldZ_{\tau})^{-1}$. Here, $(\boldI + \frac{1}{\rho}\boldZ_{\tau}\boldZ_{\tau}^\top) \in \mathbb{R}^{k\times k}$, and its inverse can be computed much more quickly than the inverse of $(\boldN_{\tau}+\rho \boldI)$ directly. 
	Note that estimating the Hessian matrix via a low-rank approximation is a common technique, \emph{e.g.} \cite{daxberger2021laplace}. The key contribution here is the use of the Nystr\"om approximation for the Hessian approximation.   
	
	Algorithm \ref{alg:Nys-SGD} and \ref{algNys-SVRG} are the Nystr\"om SGD and the Naystr\"om SVRG methods, respectively. 
	
	\begin{algorithm}[!t]
		\caption{Nystr\"om-SGD Algorithm}\label{alg:Nys-SGD}
		\begin{algorithmic}[1]
			\STATE \textbf{Initialize} ${\boldw}_0, \tau=1$, $\eta_0=\beta$, and update frequency $\ell$.
			\FOR{$t = 1, 2, \ldots$}
			\STATE randomly pick batch $\calB \sim \{1,\ldots,n\}$
			\STATE $\boldv_{t-1} = \nabla f_{\calB}({\boldw_{t-1}}) $
			\IF{$(t-1)\textbf{ mod }\ell = 0$}
			\vspace{0.1cm}
			\vspace{0.1cm}
			\STATE  $\boldC = ${\large $\frac{\partial \nabla f({\boldw})}{\partial{\boldw}_{\boldOmega}}$}, where $ \boldOmega
			\subset \{1,\ldots,d\}, |\boldOmega| = m $
			\STATE compute $\boldZ_{\tau}$ using ~\eqref{def:Nystrom} 
			\STATE $\boldQ_{\tau} = \frac{1}{\rho^2}\boldZ_{\tau}({\boldI_k} +\frac{1}{\rho} \boldZ_{\tau}^\top\boldZ_{\tau})^{-1}$
			\STATE $\tau = \tau + 1$
			\ENDIF
			\STATE Compute $\boldB_{\tau-1}\boldv_{t-1}$ using ~\eqref{eq:vect_mult}
			\STATE \textbf{$\boldw_{t} = \boldw_{t-1} - \eta \boldB_{\tau-1}\boldv_{t-1}  $}
			\ENDFOR
		\end{algorithmic}
	\end{algorithm}
	
	\begin{algorithm}[!t]
		\caption{Nystr\"om-SVRG Algorithm}\label{algNys-SVRG}
		\begin{algorithmic}[1]
			\STATE \textbf{Initialize} $\Tilde{\boldw}_0$, $\eta_0$, and update frequency $\ell$.
			\FOR{$\tau = 1, 2, \ldots$}
			\STATE $\Tilde{\boldw} = \Tilde{\boldw}_{\tau-1}$
			\STATE $\Tilde{\boldg} = \frac{1}{n} \sum_{i=1}^{n} \nabla f_i({\Tilde{\boldw}})$
			\vspace{0.1cm}
			\vspace{0.1cm}
			\STATE  $\boldC = ${\large$\frac{\partial\nabla f(\Tilde{\boldw})}{\partial\Tilde{\boldw}_{\boldOmega}}$}, where $\boldOmega\subset\{1,\ldots,d\},|\boldOmega|=m$
			\vspace{0.1cm}
			\STATE compute $\boldZ_{\tau}$ using~\eqref{def:Nystrom}
			\STATE $\boldQ_{\tau} = \frac{1}{\rho^2}\boldZ_{\tau}(\boldI_k+\frac{1}{\rho} \boldZ_{\tau}^\top\boldZ_{\tau})^{-1}$
			\vspace{0.1cm}
			\STATE $\boldw_0 = \Tilde{\boldw}$
			\FOR{$t = 1,\ldots,\ell$}
			\STATE randomly pick batch $\calB \sim \{1,\ldots,n\}$
			\STATE $\boldv_{t-1} = \nabla f_\calB({\boldw_{t-1}}) - \nabla f_\calB({\Tilde{\boldw})} + \Tilde{\boldg}$
			\STATE Compute $\boldB_{\tau}\boldv_{t-1}$ using~\eqref{eq:vect_mult}
			\STATE \textbf{$\boldw_{t} = \boldw_{t-1} - \eta \boldB_{\tau}\boldv_{t-1}  $}
			\ENDFOR
			\STATE $\Tilde{\boldw}_{\tau} = \boldw_{t}$ for randomly chosen $t \in \{1, \ldots,\ell\}$
			\ENDFOR
		\end{algorithmic}
	\end{algorithm}
	
	\noindent {\bf Computational complexity:}
	Here, we analyze the per-epoch computational complexity of the proposed method. First, it is important to note that $\boldZ_\tau,\boldQ_\tau \in \mathbb{R}^{d \times k} $ is computed at each epoch, whereas ~\eqref{eq:vect_mult} is computed at each iteration. The cost of matrix-vector multiplication $\boldB_\tau \boldv_t, i.e., $ \eqref{eq:vect_mult} is $O(dk)$ at each iteration, and therefore $O(\ell dk) $ per epoch, where $\ell$ denotes the number of iterations per epoch. The cost of computing $\boldQ_\tau$ is $O(dk^2)$ at each epoch. The cost of computing $\boldZ_\tau$ is $O(dmk)$. The computational cost of constructing the matrix $\boldC$ is $O(ndm)$. Thus, over the course of all epochs, the  construction of the matrix $\boldC$ is associated with the highest computational cost; therefore, the overall time and space complexity are $O(ndm)$ and $O(dm)$, respectively. 
	
	\subsection{Nystr\"om approximation as Newton-sketch} \label{sec:Newton-sketch}
	In this section, we show that Nystr\"om approximation can be interpreted as a Newton sketch-based method~\cite{pilanci2017newton, Lacotte2021icml}. This sketching representation of the Nystr\"om approximation is a key to analyze the convergence of the Nystr\"om-SGD method. See Section \ref{sec:theory} for a theoretical analysis of Nystr\"om-SGD.

	Let $\boldH = \nabla^2 f(\boldw)$ be a Hessian of $f(\boldw)$ of the form $\boldH = \boldX^\top\boldX, $ where $\boldX$ is an $n \times d$ matrix. It is always possible to assume that $\boldH = \boldX^\top\boldX$ because $\boldH$ is a symmetric positive semi-definite (SPSD). The embedding $\boldW \in \mathbb{R}^{d \times m}$ can be constructed as follows.
	\begin{equation}\label{eq:SamplingConstruction}
		\boldW(i,j) = \begin{cases}
			1 & \mbox{if the $i$-th column is chosen in}\\
			& \mbox{the $j$-th random trail},\\
			0 & \text{otherwise}.
		\end{cases}
	\end{equation}
	Moreover, $\boldC := \boldH\boldW$, which is the sampled column matrix of the true Hessian, and $\boldM := \boldW^\top\boldH\boldW $, which is the intersection matrix in~\eqref{def:Nystrom}. We define  $\boldX\boldW = \widehat{\boldU} \widehat{\boldSigma} \widehat{\boldV}^\top$, and $\boldM  = \widehat{\boldV} {\widehat{\boldSigma}}^2 \widehat{\boldV}^\top$.
	Then, we obtain,
	\begin{align}\label{eq:inter}
		\nonumber
		\boldC (\boldM_k)^\dagger \boldC^\top & = (\boldH\boldW) (\boldW^\top\boldH\boldW)_k^\dagger    (\boldH\boldW)^\top \\
		& = \boldX^\top \widehat{\boldU}_k \widehat{\boldU}_k^\top \boldX.
	\end{align}
	The right-hand side of \eqref{eq:inter} is similar to the Newton sketch with two differences, in that 1) $\boldX$ is replaced with the square root of the GGN matrix, and 2) the  natural orthogonal matrix $\widehat{\boldU}_k$  in proposed method  is replaced by a randomized embedding matrix $\boldS^\top$, which is expected to be orthogonal in principle, whereas in practice, it may deviate significantly.
	
	If we let $\boldX = \nabla^2 f(\boldw)^{1/2}$ then, our approximation is of the form of
	\begin{equation}\label{eq:Newton_sketched_Hessian}
		\boldH_S =  (\nabla^2 f(\boldw) ^{1/2})^\top {\boldS^\top} \boldS (\nabla^2 f(\boldw) ^{1/2}) + \rho \boldI.
	\end{equation}
	More generally, the approximation above can be written in the form of an embedding matrix as follows. Let $\boldG = \rho \boldI_d,$ and let $\boldG^{1/2} = \sqrt{\rho}\cdot\boldI_d  $ be an $d \times d$ matrix. Then, by defining the embedding matrix
	$\Bar{\boldS} = \begin{bmatrix}
		\boldS_{m \times n} & \boldzero_{m \times d} \\
		\boldzero_{d \times n} & \boldI_{d} 
	\end{bmatrix}$  and partial Hessian
	$\Bar{\boldH} = \begin{bmatrix}
		\nabla^2 f(\boldw)^{1/2} \\
		\boldG^{1/2} 
	\end{bmatrix}$, we get 
	\begin{equation}\label{eq:regularized_Newton_sketch_app}
		\boldH_{S} = \Bar{\boldH}^\top\Bar{\boldS}^\top\Bar{\boldS}\Bar{\boldH}
	\end{equation}
	which is identical to the \eqref{eq:Newton_sketched_Hessian} and hence $\boldB_{\tau} = \boldH_{\tau,S}^{-1}$ is non-singular, where $\boldH_{\tau,S}$ is the Nystr\"om approximation $\boldH_S$ at iteration $\tau$.
	
	\subsection{ Relation to $\ell_2$ regularization:} 
	The Hessian of $\ell_2$-regularized function, \emph{i.e.}, $f(\boldw) + \frac{\lambda}{2} \|\boldw\|^2$ can be given as $\boldH = \boldH_f + \lambda \boldI$, where $\lambda \geq 0$.
	In the Nystr\"om approximation for the regularized objective function, we approximate the $\boldH$ by the regularized Nystr\"om approximation by adding an independent regularizer $\rho$:
	\begin{align*}
		\boldH + \rho \boldI \approx \boldC \boldM_k^{\dagger}\boldC^\top + (\lambda+\rho) \boldI = \boldZ \boldZ^\top  \ +\ (\lambda+\rho) \boldI .
	\end{align*}
	That is, one can treat $\boldZ\boldZ^T + (\lambda +\rho)\boldI$ as the approximated Hessian of $f(\boldw) + \frac{(\lambda +\rho)}{2}\|\boldw\|_2^2$, which makes approximation invertible  even in the case where $\lambda =0$. Usually, $\rho$ is obtained by trust region techniques which also depends on user defined constant. Therefore, we treat $\rho$ as a hyperparameter and we can obtain the best $\rho$ using grid search.

	\section{Theoretical Analysis\label{sec:theory}}
	We show the theoretical analysis of Nystr\"om-SGD and Nystr\"om-SVRG methods for convex functions. Note that we provide a theoretical analysis for the objective function $f$ given in \eqref{objective_fun}, which does not include the $\ell_2$ regularization.

	\subsection{Assumptions and Lemmas} 
	In this section, we present a convergence analysis, which involves the following assumptions.
	
	\begin{assumption}\label{asm:convex}
		We assume that each $f_i$ is convex and twice continuously differentiable.
	\end{assumption}
	
	\begin{assumption}\label{asm:Hessian_bound}
		There exists two positive constants
		$\mu$ and $\Lambda$ such that
		\begin{equation*}
			\mu \boldI \preceq \nabla^2 f(\boldw) \preceq \Lambda \boldI,\quad \boldw\in \mathbb{R}^d,
		\end{equation*}
		The lower bound always holds in the regularized case.
	\end{assumption}
	
	\begin{assumption}\label{asm:Lipschitz}
		We assume that the gradient of each $f_i$ is $\Lambda$-Lipschitz continuous, $i.e.,$
		\begin{equation*}
			\|\nabla f_i(\boldw_a) - \nabla f_i(\boldw_b)\| \leq \Lambda ~\|\boldw_a-\boldw_b\| ~\forall \boldw_a,\boldw_b \in \mathbb{R}^d.
		\end{equation*}
		Under this assumption, it is clear that $\nabla f$ is also $\Lambda$-Lipschitz continuous.
		\begin{equation*}
			\|\nabla f(\boldw_a) - \nabla f(\boldw_b)\| \leq \Lambda~\|\boldw_a-\boldw_b\| ~\forall \boldw_a,\boldw_b \in \mathbb{R}^d.
		\end{equation*}
	\end{assumption}
	Note that these are the standards assumptions which have been used in~\cite{Moritz2016linearly,Byrd2016stochastic}. Moreover, SQN~\cite{Moritz2016linearly} and SVRG-SQN~\cite{Moritz2016linearly} both converge linearly.
	
	Theorem~\ref{thm:nystromerrorbound} provides the error bound for the Nystr\"om approximation with respect to the best $k$-rank approximation.
	
	\begin{theorem}\cite{drineas2005nystrom}\label{thm:nystromerrorbound}
		Let $\boldH$ be a $d\times d$ matrix and let $\boldH_k$ be the best $k$-rank approximation
		of the $\boldH$. Then, for $O(k/\epsilon^4)$ columns
		\begin{align}\label{eq:Nystromerrorbound}
			\|\boldH - \boldC\boldM^\dagger_k\boldC^\top\|_{\nu} \leq \|\boldH - \boldH_k\|_{\nu} + \epsilon \sum_{i=1}^{d} \boldH^2_{ii},
		\end{align}
		where $\epsilon>0$ and $\nu = 2$ (spectral) or $\nu = F$ (Frobenius). 
	\end{theorem}
	Next, we connect the Nystr\"om error bound~(Theorem~\ref{thm:nystromerrorbound}) with the upper bound of the Nystr\"om approximation.
	\begin{lemma}\label{lem:lowerboundNystrom}
		Let $\boldN = \boldZ\boldZ^\top$ satisfy Theorem~\ref{thm:nystromerrorbound}: Then, Nystr\"om approximation matrix $\boldN$ is bounded, as follows:
		\begin{equation}
			\label{asm:approx_hess_bound_main}
			\boldzero \preceq \boldN \preceq \Gamma \boldI
		\end{equation}
	\end{lemma}

	Then, we provide the bound of the regularized Nystr\"om approximation with the help of the Assumption~\ref{asm:convex} and \ref{asm:Hessian_bound}.
	
	\begin{lemma} \label{lem:Nystrom_bound}
		Let $\boldN$ satisfy ~\eqref{asm:approx_hess_bound_main}. Then the following bounds exist.
		\begin{equation*}
			\rho \boldI \preceq (\boldN + \rho \boldI)
			\preceq (\Gamma + \rho \boldI), \text{~ for } \rho > 0.
		\end{equation*}
	\end{lemma}
	
	The following lemma provides the bounds to the inverse approximation $\boldB = (\boldN + \rho\boldI)^{-1}$ using the above lemma.
	\begin{lemma}\label{lem:InverseNB}
		Let Assumption~\ref{asm:convex}, \ref{asm:Hessian_bound} and Lemma~\ref{lem:lowerboundNystrom}
		Hold. Subsequently, the regularized Nystr\"om approximation 
		$(\boldN + \rho \boldI)^{-1}$ is bounded.
	\end{lemma}
	
	The following lemma provides the upper bound of the variance of the variance-reduced gradient:
	
	\begin{lemma}\label{lem:expectationbound}
		Let $\boldw_*$ be a unique minimizer of $f$ and let $\boldv_t = \nabla f_\calB(\boldw_t) - 
		\nabla f_\calB(\Tilde{\boldw}) + \nabla f(\Tilde{\boldw})$ be a variance-reduced stochastic 
		gradient with mini batch $\calB \subseteq \{1,\ldots,n\}$. Then, the expectation with 
		respect to $\calB$ is bounded;
		\begin{equation*}
			\mathbb{E}\|\boldv_t\|^2 \leq 4\Lambda (f(\boldw_t) - f(\boldw_*) + f(\boldw_{\tau}) - f(\boldw_*)).
		\end{equation*}
	\end{lemma}
	This bound follows from the~\cite[Lemma 6]{Moritz2016linearly} which closely follows from
	~\cite[Theorem 1]{Johnson2013accelerating}.

	At this point, we apply Lemma~\ref{lem:gradient_bound}, which states a result of a strongly convex function.
	\begin{lemma} \label{lem:gradient_bound}
		Suppose that $f$ is continuously differentiable and strongly convex with parameter $\mu$.
		Let $\boldw_*$ be a unique minimizer of $f$. Subsequently, for any $\boldw \in \mathbb{R}$, we have 
		\begin{equation*}
			\|\nabla f(\boldw) \|^2 \geq 2\mu (f(\boldw) - f(\boldw_*)).
		\end{equation*}
	\end{lemma}
	
	\subsection{Convergence analysis of Nystr\"om SGD} Here, we show the convergence of Algorithm \ref{alg:Nys-SGD}.
	\begin{assumption}\label{asm:Lipshictz_Hessian}
		We assume that the Hessian of the objective function is Lipschitz continuous with Lipschitz constant $L_H$, such that
		$$ \|\nabla^2 f(\boldw_1) - \nabla^2 f(\boldw_2) \| \leq L_H \|\boldw_1 - \boldw_2\|.$$
	\end{assumption}
	Because the approximation \eqref{eq:Newton_sketched_Hessian} is similar to that of the Newton sketch, we follow a similar analysis for the Newton-like iteration for the Nystr\"om SGD algorithm. We consider the Newton-like iteration: 
	\begin{equation}\label{eq:Newton_like_iteration}
		\boldw_{t} = \boldw_{t-1} - \eta \boldB_{t-1}\nabla f(\boldw_{t-1}),
	\end{equation}
	as Nystr\"om-SGD is a specific case of Newton-like iteration. Note that Algorithm~\ref{alg:Nys-SGD} is a special case of \eqref{eq:Newton_like_iteration} where $\boldB_{\tau}$ is constant for $\ell$ iterations.
	
	\begin{definition}[Gaussian width]
		For a $d$-dimensional compact set $\mathcal{C} $, Gaussian width is given by
		\begin{equation}
			\mathcal{W}(\mathcal{C}) := \mathbb{E}_\boldy\left[ \max_{z\in \mathcal{C}} |\langle \boldy,\boldz \rangle| \right],
		\end{equation}
		where $\boldy \in \mathbb{R}^d$ is an i.i.d. sequence of $N(0,1)$ variables. 
	\end{definition}
	
	Let $\calL$ is constraint set for the given problem $f(\boldw_*)$ and $\calK$ be a tangent cone at $\boldw_*$, such that 
	\begin{equation}
		\calK = \{\boldv \in \mathbb{R}^d | \boldw_* + \eta \boldv \in \calL, \text{\ for\ } \eta > 0 \}
	\end{equation}
	Note that in the case of unconstrained optimization problem, $\calL = \mathbb{R}^d$ and therefore $\calK = \mathbb{R}^d$.
	
	Sketch dimension $m$ depends on the Gaussian width. Since, our approximation is similar to Newton-sketch, we take the same lower bound for $m$,
	\begin{equation}\label{eq:Gaussian_width_lower}
		m \geq \frac{c}{\epsilon^2} \max_{w\in \mathbb{R}^d} \calW^2(\nabla^2 f(\boldw)^{1/2}\calK)
	\end{equation}
	where $\hat{\epsilon} $ is user defined tolerance and $c$ is constant.
	
	The next theorem shows linear-quadratic convergence for the Nystr\"om SGD algorithms.
	\begin{theorem}
		\label{thm:linear-quadratic}
		Let Assumption \ref{asm:Lipshictz_Hessian} holds. $\gamma = \lambda_{\min} (\nabla^2f(\boldw_*)), \beta = \lambda_{\max} (\nabla^2 f(\boldw_*))$, and $\boldH_s$ be an approximation given in \eqref{eq:regularized_Newton_sketch_app} 
		and $\hat{\epsilon} \in \left(0,\frac{2\gamma}{9\beta} \right)$. Consider the 
		Nystr\"om sketch update ~\eqref{eq:Newton_like_iteration} initialized with $\boldw_0$ such that $\|\boldw_0 - \boldw_*\| \leq \frac{\gamma}{8L_H}$, and a sketch dimension $m$ satisfy the lower bound 
		\begin{equation*}
			m \geq \frac{c}{\hat{\epsilon}^2} \max_{w\in \mathbb{R}^d} \calW^2(\nabla^2 f(\boldw)^{1/2}\calK).
		\end{equation*}
		Then with probability at least $1-c_1Te^{-c_2m},$ the Eucliden error satisfies the bound
		\begin{equation}
			\|\boldw_{t+1} - \boldw_*\| \leq \hat{\epsilon} \frac{\beta}{\gamma} \|\boldw_t - \boldw_*\| + \frac{4L_H}{\gamma}\|\boldw_t - \boldw_*\|^2
		\end{equation}
		where $T$ is the total iteration, $\mathcal{W}$ is the Gaussian width and $\mathcal{K}$ is the tangent cone at $\boldw_*$.
	\end{theorem} 
	Since the Nystr\"om approximation can be interpreted as a Newton sketch \eqref{eq:regularized_Newton_sketch_app}, we can obtain Theorem~\ref{thm:linear-quadratic} by using \cite[Theorem 3.1]{pilanci2017newton}. Note that we can use the same lower bound on $m$ as our embedding matrix $\boldS$ has similar properties as of \cite{pilanci2017newton}.
	
	\subsection{Convergence analysis of Nystr\"om-SVRG}   The next theorem states the main result for the convergence of Algorithm~\ref{algNys-SVRG}. We follow an analysis similar to that of \cite{Moritz2016linearly}.
	\begin{theorem}\label{thmNys-SVRG}
		Suppose Assumptions~\ref{asm:convex} - \ref{asm:Lipschitz}, and Lemma~\ref{lem:lowerboundNystrom} hold, and let $\boldw_*$ be a unique minimizer of the objective function~\eqref{objective_fun}. Subsequently, for all $\tau \geq 0$,
		\begin{equation}
			\mathbb{E}[f(\boldw_{\tau}) - f(\boldw_*)] \leq \alpha^{\tau} ~\mathbb{E}[f(\boldw_0) - f(\boldw_*)],
		\end{equation}
		where
		$\alpha =  \frac{1 + 2 \ell \eta^2 \Lambda^2 \delta^2}{2\ell \eta (\mu \Delta - \eta \Lambda^2 \delta^2)} < 1$, $\Delta = {1}/{(\Gamma + \rho )}$ and $\delta = {1}/{\rho}$, assuming $\eta < \mu \Delta / 2 \Lambda^2 \delta^2$ and choosing a sufficiently large $\ell$ to satisfy
		\begin{equation} \label{eq:l_bound}
			\frac{1}{2\ell \eta} + 2 \eta \Lambda^2 \delta^2 < \mu \Delta.
		\end{equation}
	\end{theorem}
	
	\begin{theorem}\label{Th:New_linearconvergence_svrg}
		Suppose Assumptions ~\ref{asm:convex} - \ref{asm:Lipschitz}, and Lemma~\ref{lem:lowerboundNystrom} hold, and let $\boldw_*$ be an optimal point of the objective function~\eqref{objective_fun}. Then, the distance between two consecutive iterations to the optimal is bounded in expectation as follows.
		\begin{equation}
			\E \|\boldw_{\tau} - \boldw_*\|^2  <  \zeta ~\|\boldw_{\tau -1} - \boldw_*\|^2,
		\end{equation}
		where 
		\begin{equation}
			\zeta = \left[ \left( 1 - 2 \eta (\mu \Delta - \eta \Lambda^2 \delta)\right)^{\ell}
			+ \frac{\eta \Lambda^2 \delta}{ (\mu \Delta - 2\eta \Lambda^2 \delta)} \right].
		\end{equation}
	\end{theorem}
	
	\subsection{Closeness to Newton's method}
	Let $\boldH$ be the Hessian of the $\ell_2$-regularized objective function where $\lambda$ is $\ell_2$-regularizer. Then, the inverse of Hessian of is given by $ \boldH^{-1}_\boldw = (\nabla^2 f(\boldw) +\lambda \boldI)^{-1}$. Let the approximated Nystr\"om at $\boldw$ be given by $ (\boldZ_\boldw\boldZ_\boldw^\top + \lambda \boldI)^{-1}$. The distance of the regularized inverse matrix is then given as
	\begin{align}\label{eq:Newtoncloseness}
		\|(\boldZ_\boldw\boldZ_\boldw^\top + \lambda \boldI)^{-1} - \boldH_\boldw^{-1}\| \leq \frac{\|\boldJ_\boldw\|}{\lambda (\|\boldJ_\boldw\| + \lambda)} ,
	\end{align}
	where $0 < \| \boldJ_\boldw\| = \| \nabla^2 f(\boldw) - \boldZ_\boldw \boldZ_\boldw^\top\| \leq \sigma_{k+1} + \epsilon \sum_{i=1}^{d}\boldH^2_{ii} $;
	which follows from ~\eqref{eq:Nystromerrorbound}, whereas ~\eqref{eq:Newtoncloseness} follows from \cite[Proposition 3.1]{frangella2021randomized}.
	
	Moreover, if $\nabla^2 f(\boldw)$ has rank $k$, then $\|\boldJ_{\boldw}\|\leq \epsilon \sum_{i=1}^{d}\boldH^2_{ii}$. Hence, the upper bound of \eqref{eq:Newtoncloseness} depends on error $\epsilon > 0$. One can establish the superlinear convergence by showing $\|\boldv_{ny} - \boldv_{ne}\| = O(\|\boldv_{ny}\|)$~\cite{Nocedal2006numerical} using~\eqref{eq:Newtoncloseness}, where $\boldv_{ny}$ and $\boldv_{ne}$ are the search direction of the Nystr\"om-SGD and Newton's method, respectively. This implies that superlinear convergence is quite difficult in the case of Nystr\"om approximation. This problem remains as an interesting future work.

	\section{Related work}
	
	\noindent {\bf Newton methods:} 
	The evolution of the approximation of the Hessian or its inverse began with DFP~\cite{osti_4252678, fletcher1963rapidly}, Broyden~\cite{broyden1965class}, SR1~\cite{byrd1996analysis}, and well-known BFGS methods~\cite{broyden1969new, fletcher1970new, goldfarb1970family, shanno1970conditioning}. The BFGS method uses secant equations to approximate the Hessian inverse, and requires $O(d^2)$ memory to store the approximate Hessian.  \cite{Liu1989limited} proposed a limited-memory BFGS (L-BFGS) to solve large-scale unconstrained optimization problems. They used L-BFGS to compute an approximation of the Hessian inverse by using $O(md)$ memory instead of $O(d^2)$, where $m \ll d$. 
	
	To address the increased demand for optimization in machine and deep learning, various first-order stochastic optimization methods~\cite{Johnson2013accelerating, Duchi2011adaptive, Defazio2014saga, Kingma2015adam, Nguyen2017sarah} have been proposed, but their sensitivity to the choice of hyperparameters remains a concern. Owing to the low sensitivity of the choice of hyperparameters, stochastic quasi-Newton methods have gained popularity in recent years. A stochastic quasi-Newton method by~\cite{Byrd2016stochastic} proposed an approximated sub-sampled Hessian using the Hessian vector product instead of using the difference of gradients in the secant equation. \cite{Moritz2016linearly} proposed a stochastic L-BFGS using a similar technique along with mini-batch variance-reduced gradients. \cite{Kolte2015accelerating} proposed a strategy to compute the sub-sample Hessian followed by the singular value thresholding. As a generalization of quasi-Newton methods,~\cite{scieur2021generalization} proposed the use of robust symmetric multi-secant updates to approximate the Hessian.
	
	Recently, \cite{yang2020structured} approximated the Hessian inverse $\boldB$ by decomposing it into two parts such that $\boldB = \boldB_1 + \boldB_2$, where $\boldB_1$ is either a sub-sampled Hessian, generalized Gauss-Newton matrix, Fisher information matrix (FIM), or any other approximation and $\boldB_2$ is the L-BFGS approximation. Observe that using the sub-sampled Hessian can efficiently compute $\boldB_1$ for low-dimensional problems. However,  for high-dimensional data where $d$ is significantly large, computing the subsample Hessian or FIM $\boldB_1$ can be computationally expensive because it requires $O(nd^2)$ time and $O(d^2)$ memory. Additionally, computing the sub-sample Hessian or FIM at each iteration makes it more computationally intensive. AdaHessian \cite{yao2021adahessian} is an alternative Hessian approximation method, which estimates the diagonal of Hessian matrix
	with Hutchinson’s method. In contrast, our proposed methods approximates the full Hessian matrix by using the Nystr\"om method.
	
	\cite{talwalkar2010matrix} proposed the Nystr\"om logistic regression algorithm, where the Nystr\"om method is used to approximate the Hessian of the regularized logistic regression. Thus, it can be regarded as a variant of Nystr\"om-SGD. However, \cite{talwalkar2010matrix} only considered the regularized logistic regression, in which the Hessian can be explicitly obtained, with deterministic optimization. In contrast, we propose the Nystr\"om method for a general Hessian matrix for stochastic optimization and show its efficacy for deep learning models. Moreover, we elucidate the theoretical properties of the Nystr\"om-SGD algorithm for convex functions (See Appendix \ref{sec:nys-lr} for details).

	\noindent {\bf Natural gradient methods:} The Newton method is closely related the natural gradient descent \cite{amari1998natural} methods, which can be used for probabilistic models. More specifically, in probabilistic models, because the negative Hessian is equivalent to FIM, we can update the Newton-like algorithm using FIM. However, FIM computations are expensive. To handle this computational problem, several structural approximations of FIM have been developed \cite{GrosseM16ICML,mishkin2018slang,tran2020bayesian,lin2021tractable,lin2021structured,karakida2020understanding,daxberger2021laplace}.
	However, these approximation methods are based on FIM and focus more on the computational aspects. In contrast, in this study, we propose the use of the Nystr\"om approximation of the Hessian matrix and elucidate its theoretical properties for Newton-like algorithms.

	\noindent {\bf Newton-sketch method:} Because the cost of Hessian computation and its inverse is not affordable, some significant work has been performed on sketch-based methods. ~\cite{Lacotte2021icml} and ~\cite{pilanci2017newton} proposed the use of Newton-sketch-based methods. They have used various embedding techniques such as a subsampled randomized Hadamard transform (SRHT), sub-Gaussian, randomized orthogonal systems (ROS), and sparse Johnson-Lindenstrauss transform (SJLT). It should be noted that the approximation of the Hessian matrix and its accuracy highly depend on the embedding matrix. One of the drawbacks of the sketching method is that it needs to compute the square root of the Hessian matrix although the size of the embedding matrix is $m\times n$, which makes it impractical for both large-scale and high-dimensionality data. In contrast, the Nystr\"om approximation only needs to compute the partial Hessian matrix. It is important to note that Adaptive Newton sketch ~\cite{Lacotte2021icml} increases the $m$ size eventually, whereas the proposed method does not change the value of $m$ through out the algorithm. 

	\section{Experiments}
	In this section, we validate our proposed method for both convex and non-convex setups.
	
	\subsection{Real-world experiments (Convex setup)}
	We evaluated the performance of the proposed methods for logistic regression and $\ell_2$-svm on several benchmark datasets (See Table~\ref{tab:datasets} in the supplementary material).
	
	\begin{figure}[!t]
		\centering
		\includegraphics[width=1\linewidth,keepaspectratio]{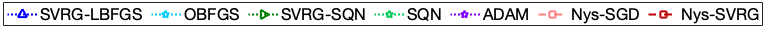}
		\includegraphics[width=.32\linewidth,height=3cm]{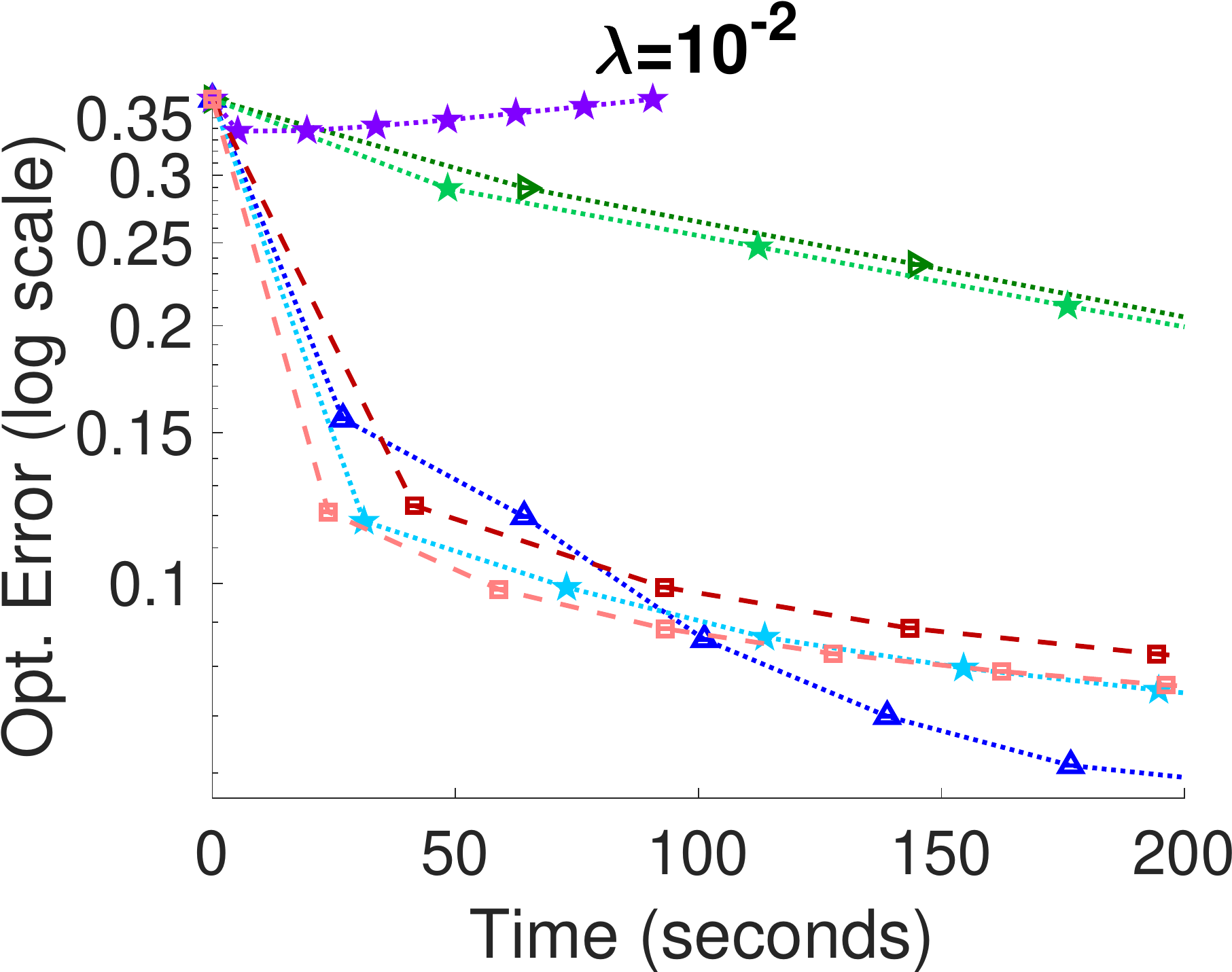}
		\includegraphics[width=.32\linewidth,height=3cm]{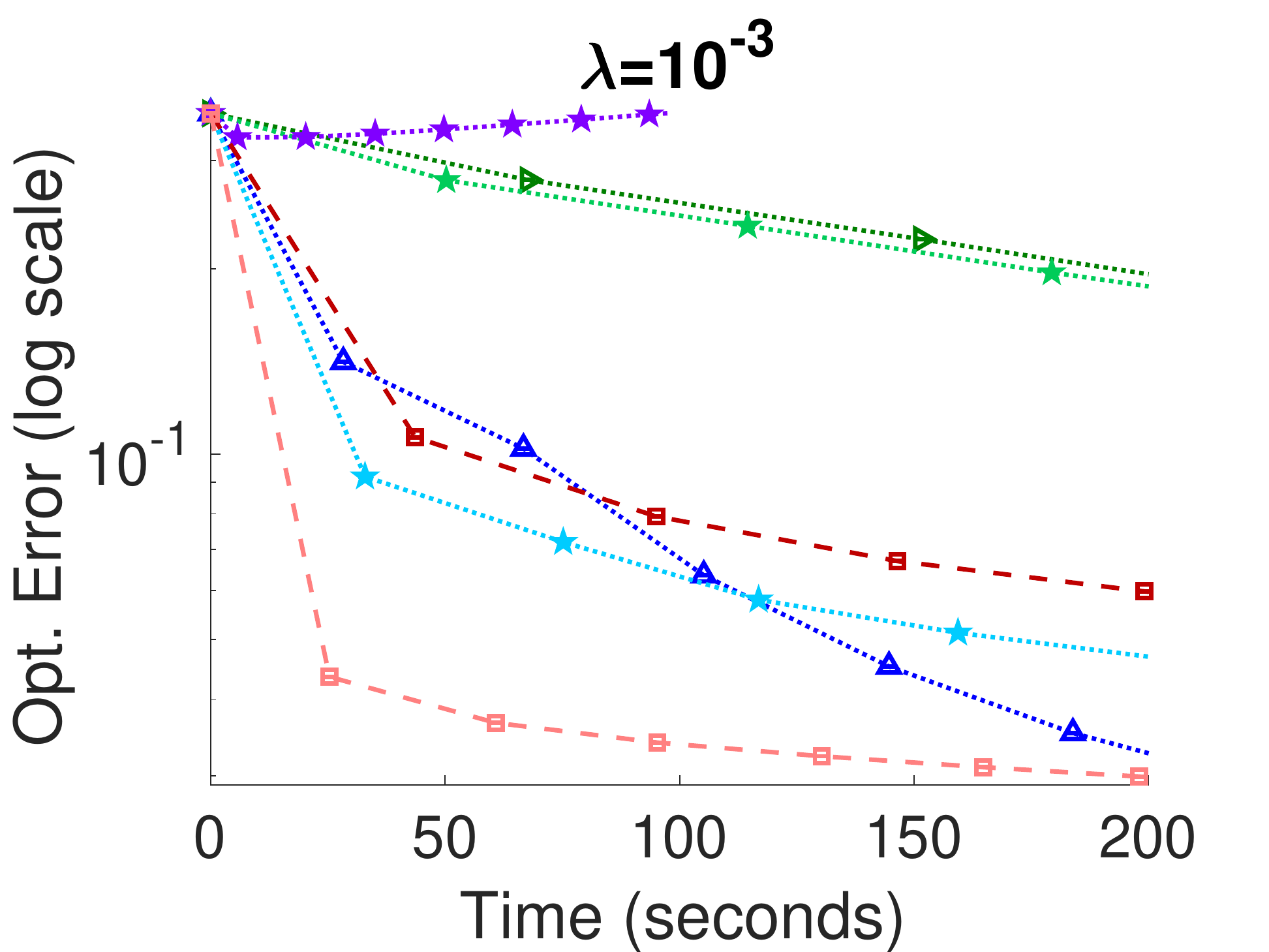}
		\includegraphics[width=.32\linewidth,height=3cm]{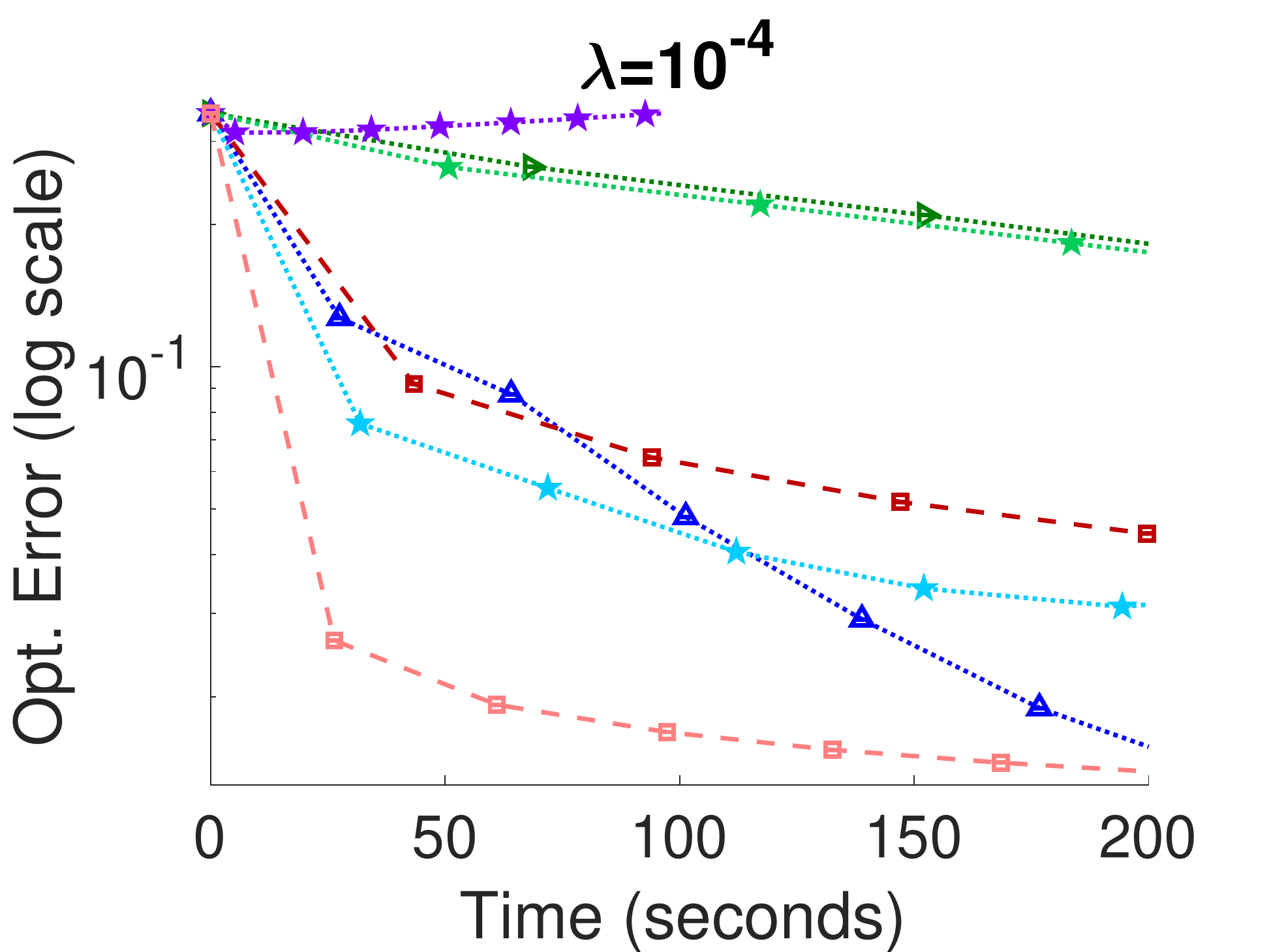} \\
		\includegraphics[width=.32\linewidth,height=3cm]{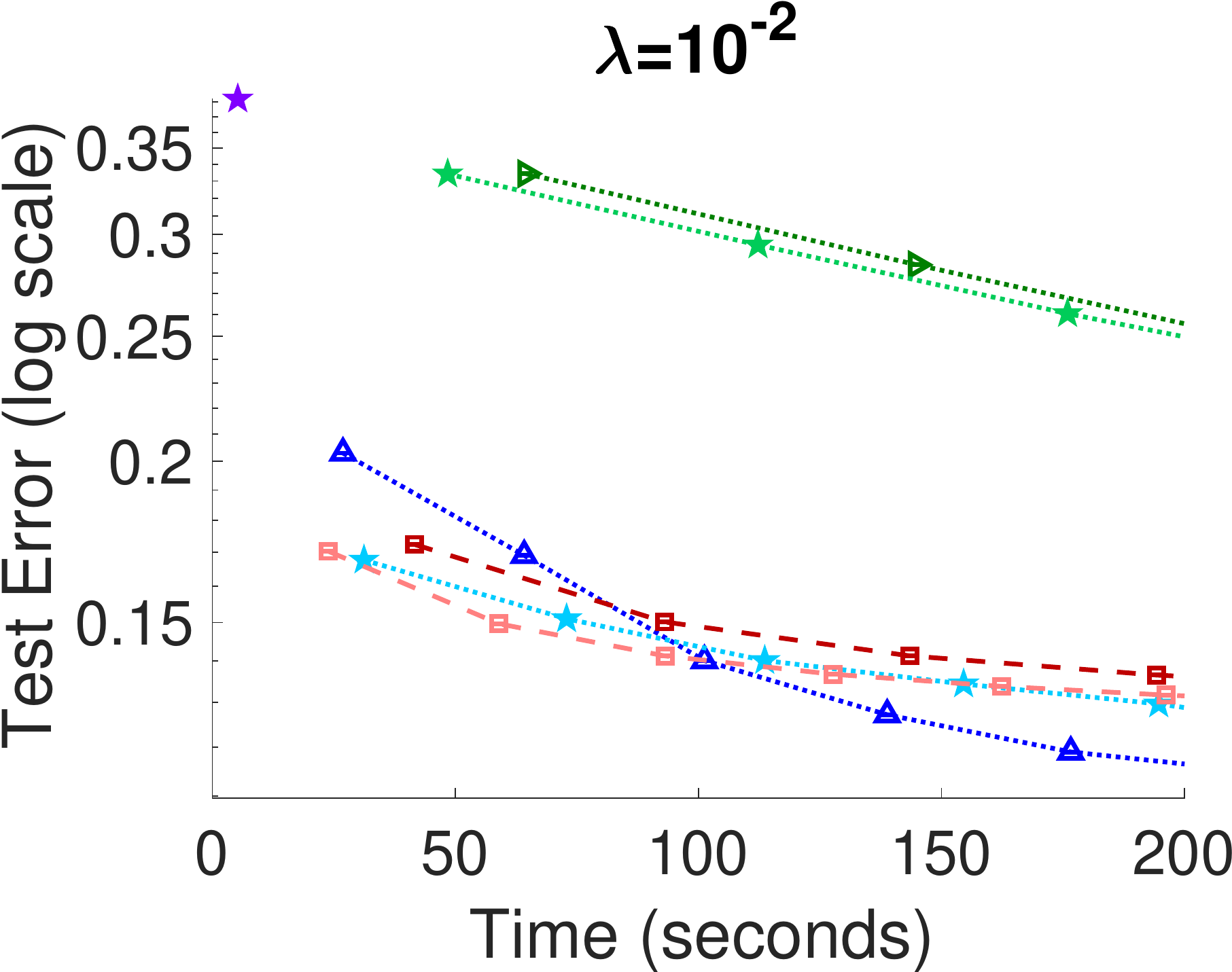}
		\includegraphics[width=.32\linewidth,height=3cm]{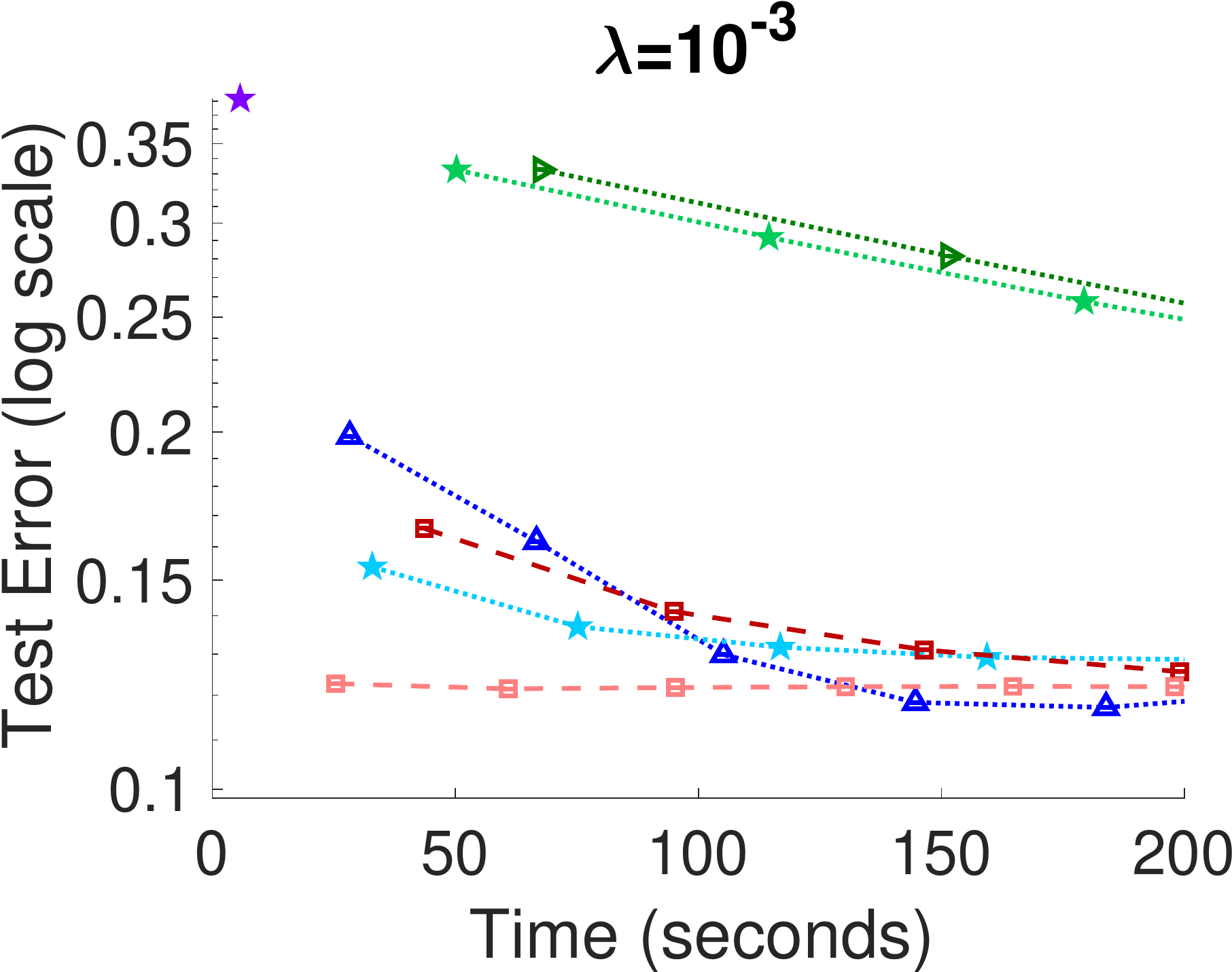}
		\includegraphics[width=.32\linewidth,height=3cm]{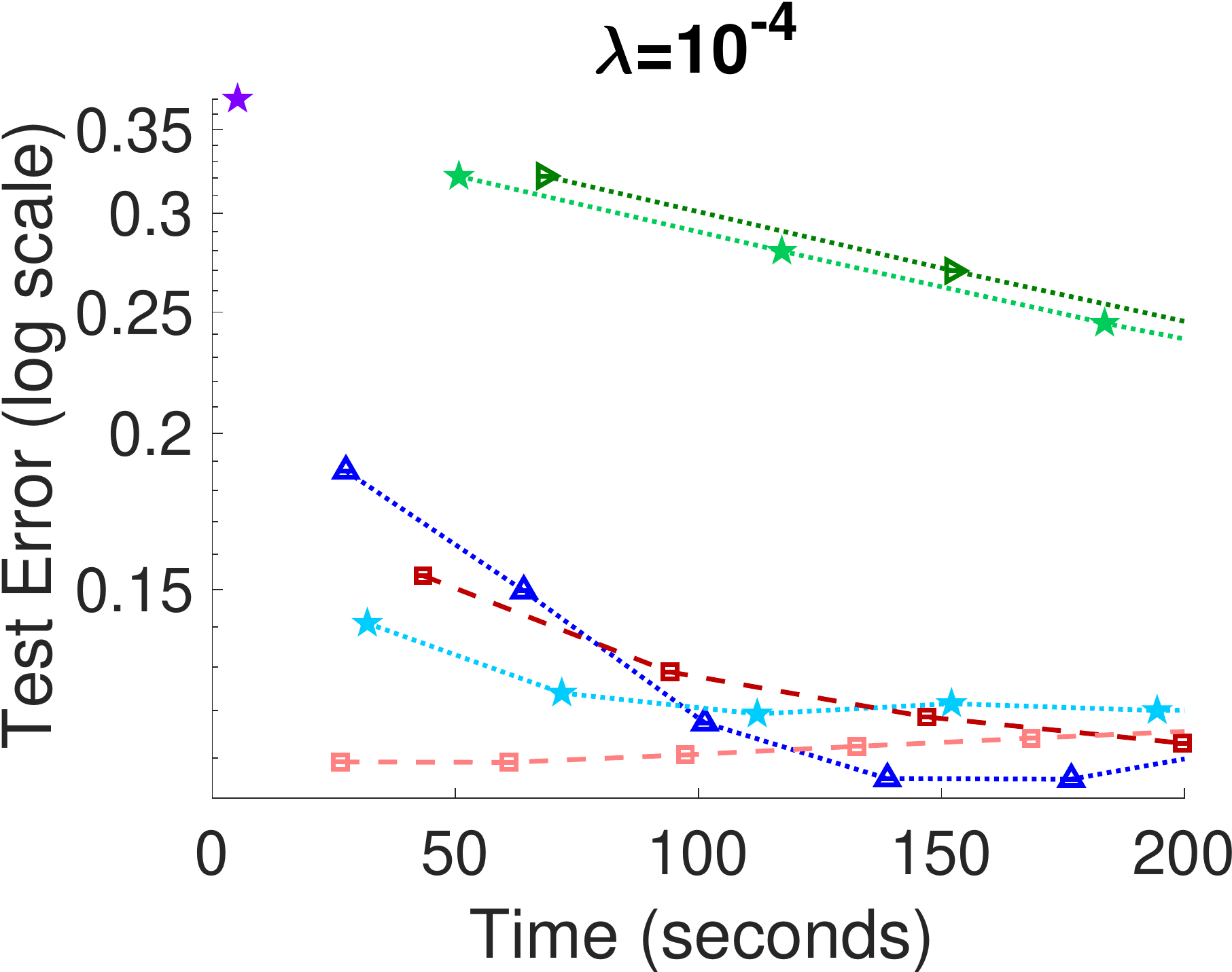} \\
		\includegraphics[width=.32\linewidth,height=3cm]{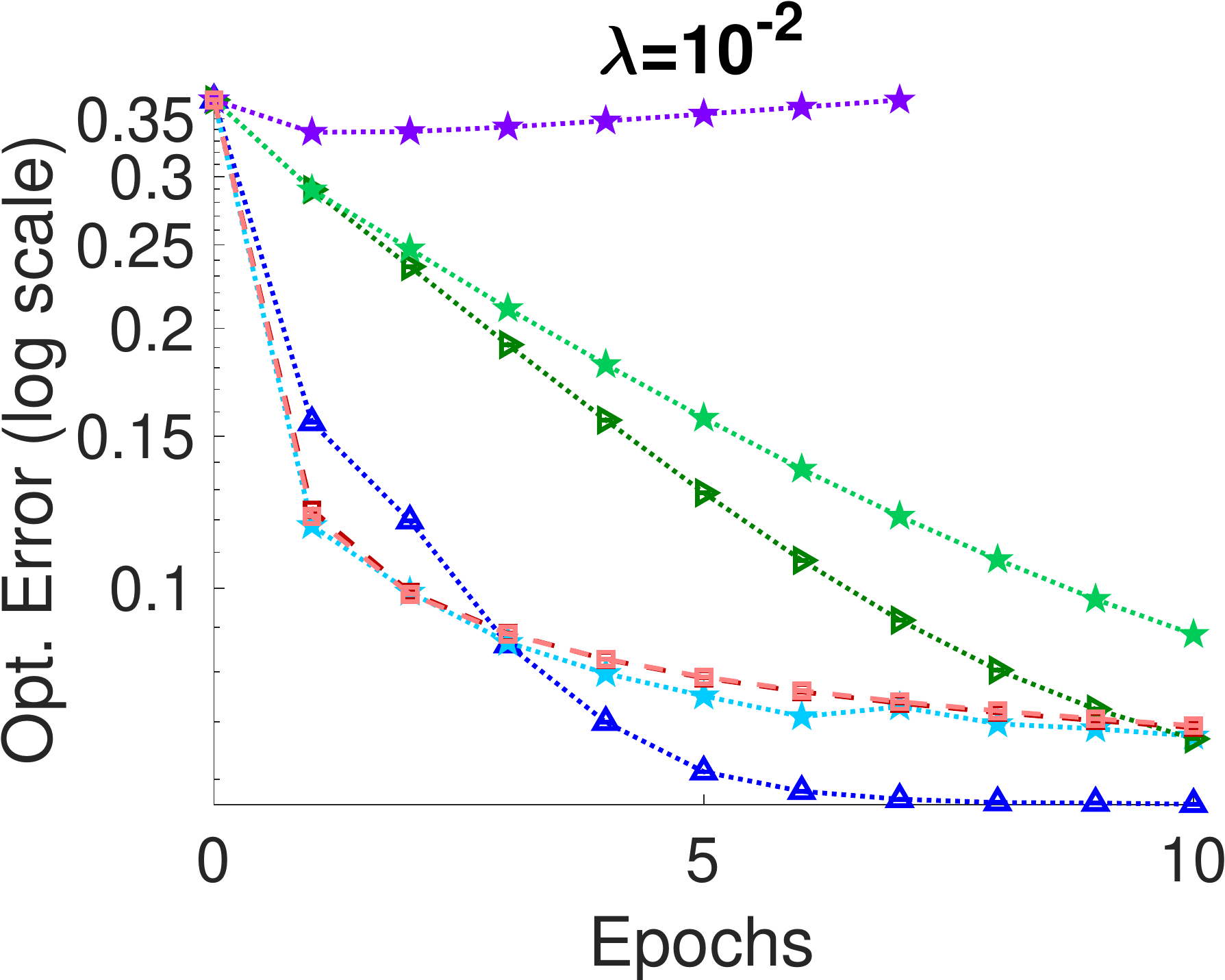}
		\includegraphics[width=.32\linewidth,height=3cm]{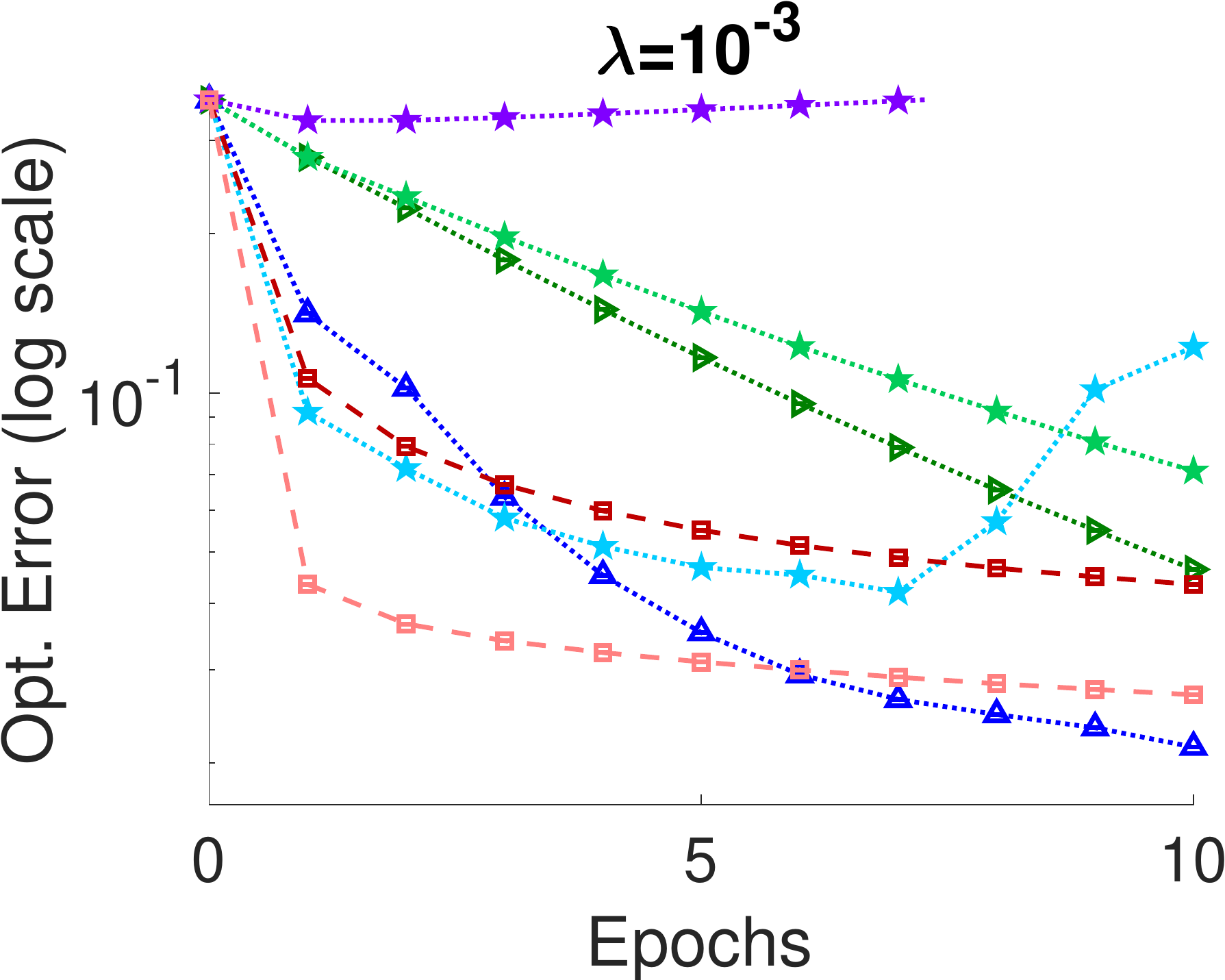}
		\includegraphics[width=.32\linewidth,height=3cm]{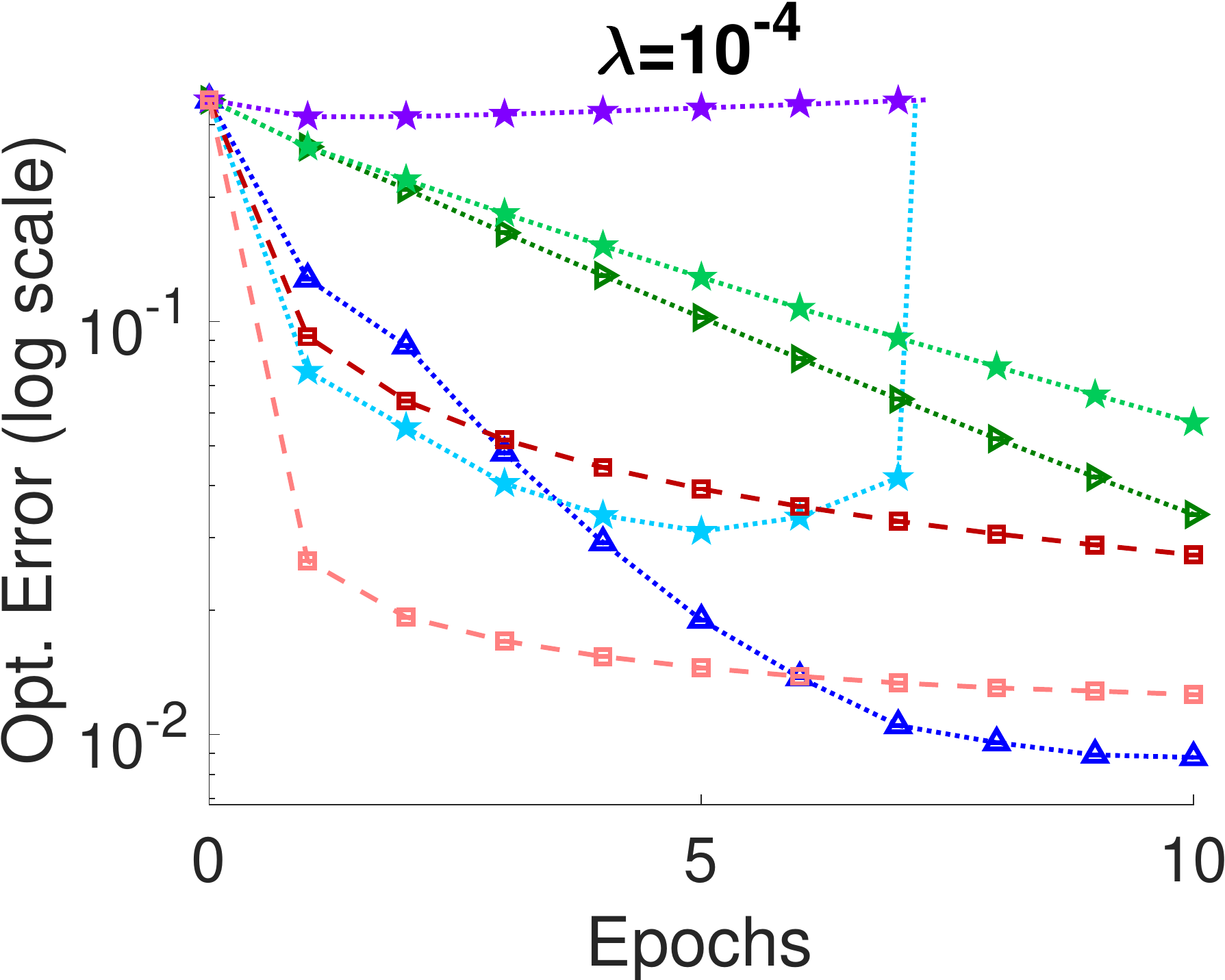} \\
		\includegraphics[width=.32\linewidth,height=3cm]{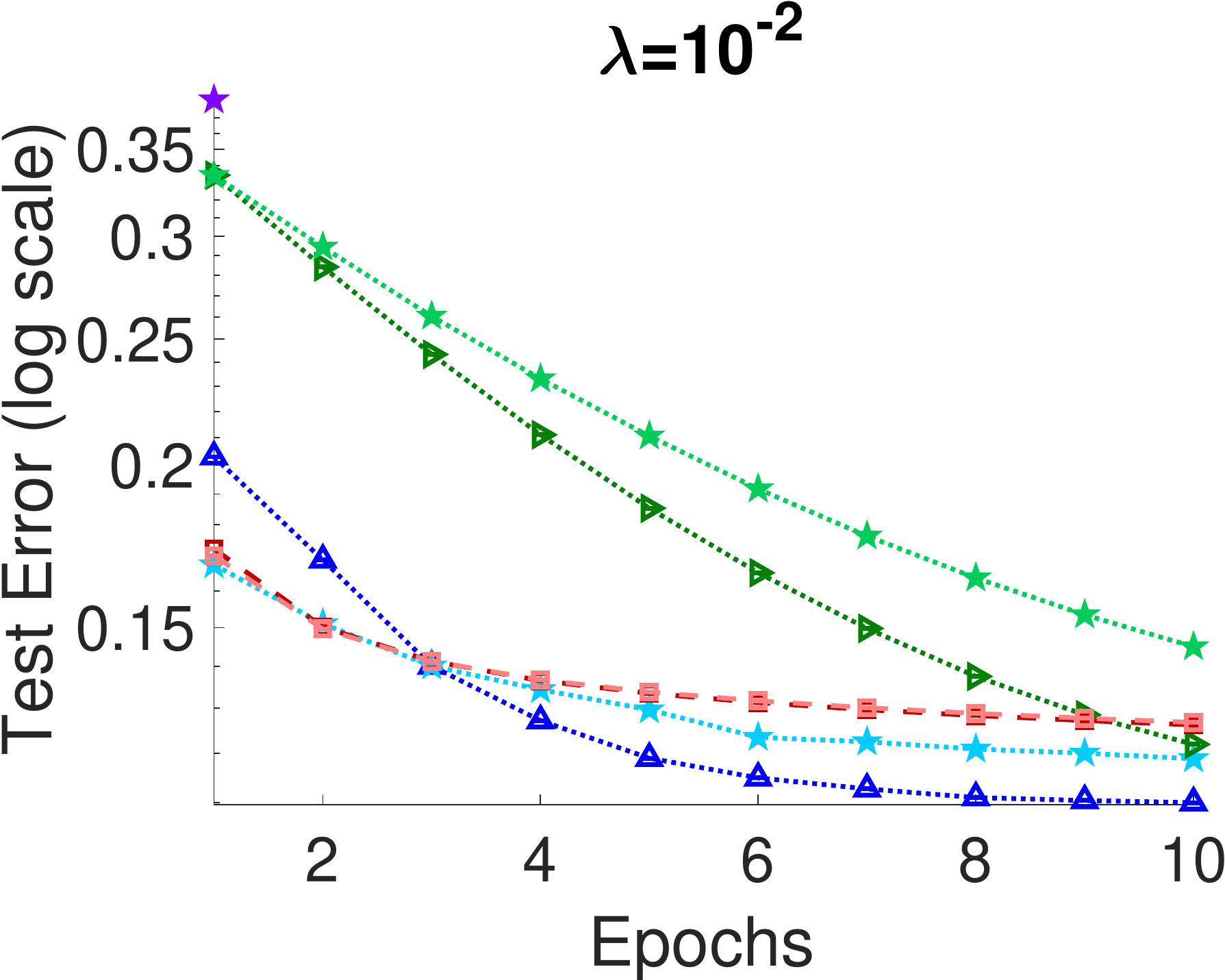}
		\includegraphics[width=.32\linewidth,height=3cm]{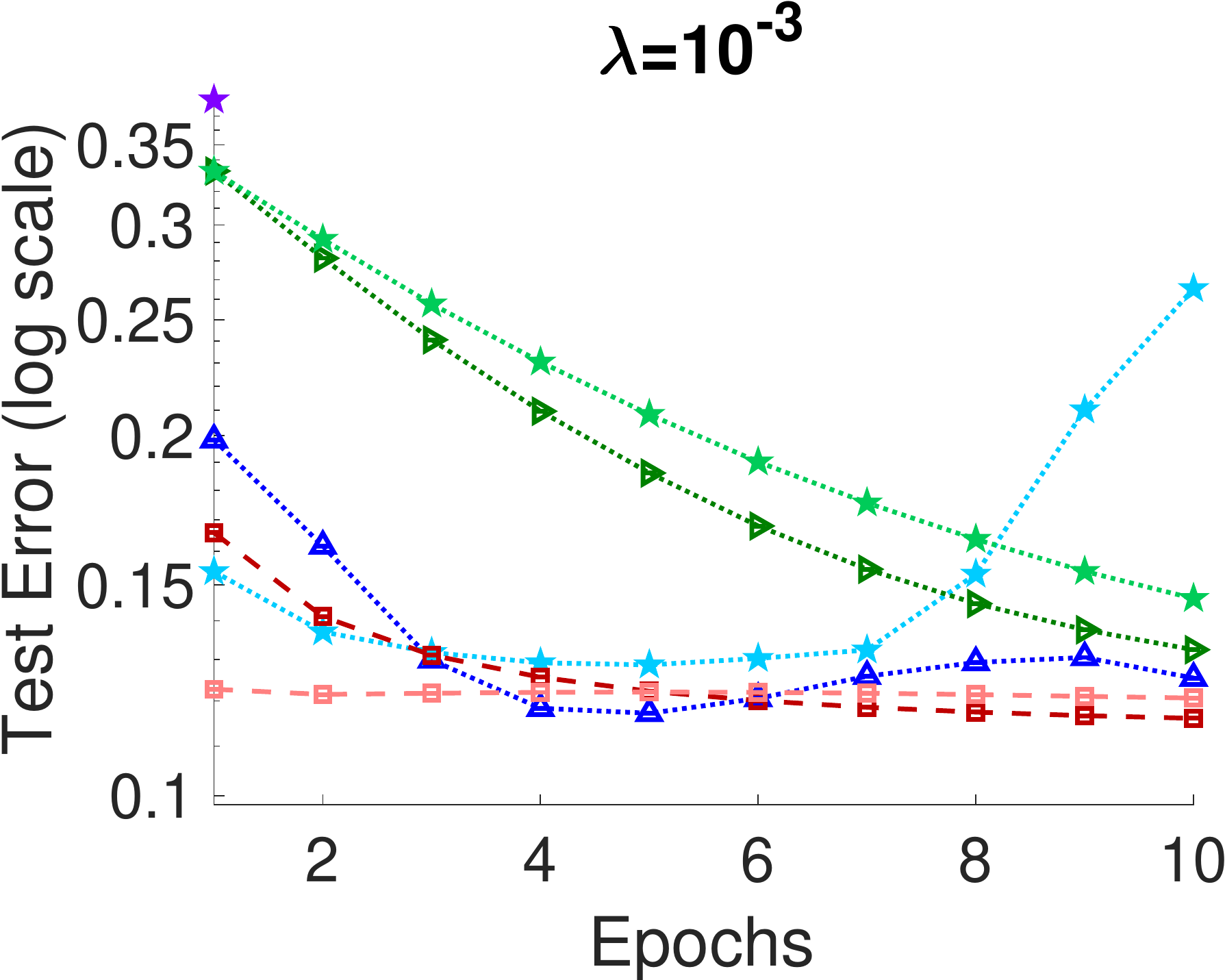}
		\includegraphics[width=.32\linewidth,height=3cm]{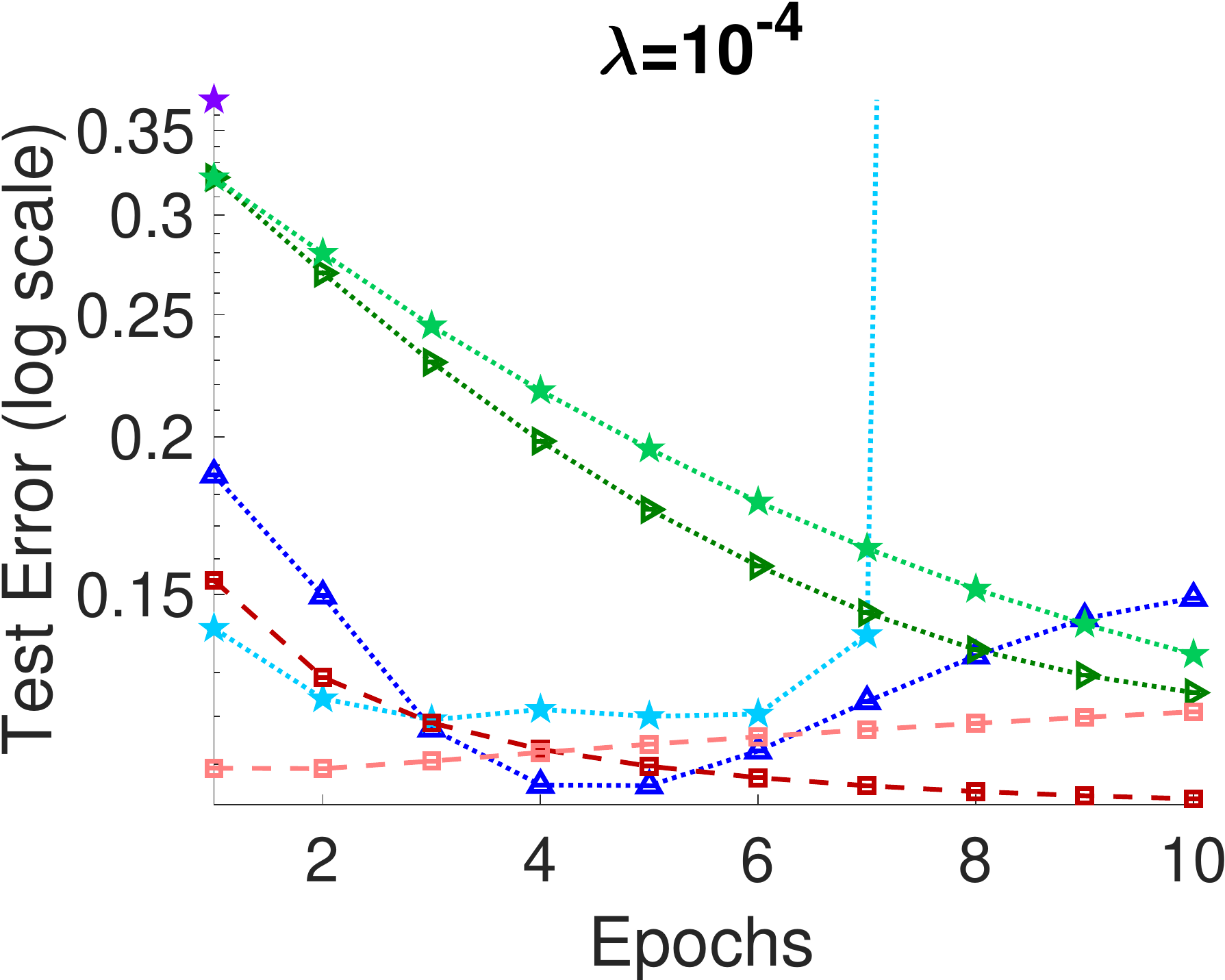}
		\caption{Comparison of the proposed method with the existing methods on \textit{real-sim} dataset for optimization error and test error with various regularizer. Top two rows are \emph{w.r.t.} time and bottom two rows are \emph{w.r.t.} epochs.}
		\label{fig:realsim}
		\vspace{-.2in}
	\end{figure}
	
		\begin{figure*}[t]
		\centering
		\begin{subfigure}[t]{0.235\textwidth}
			\centering
			\includegraphics[width=0.99\textwidth,keepaspectratio]{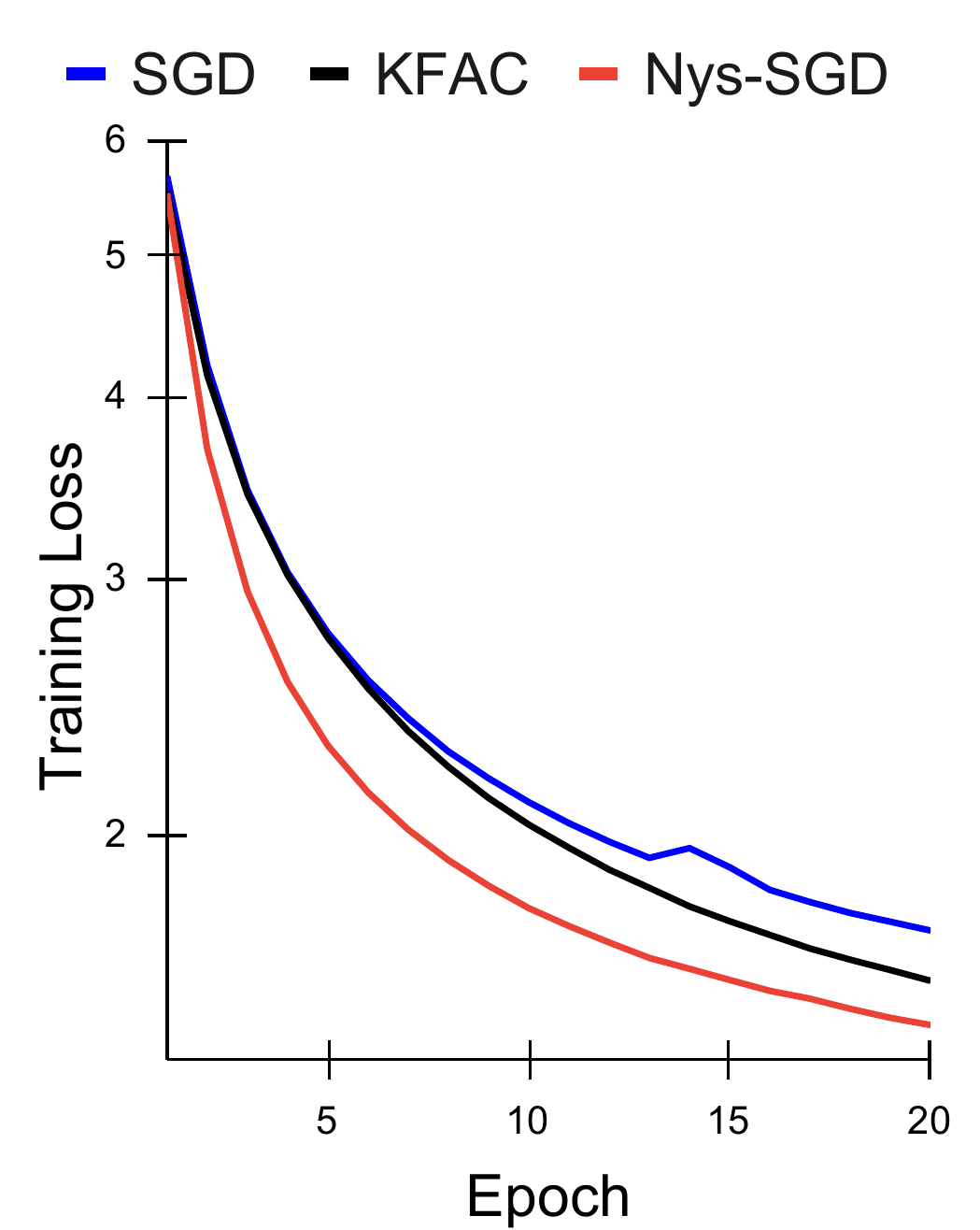}
			\caption{ResNet152 (Train loss)}
		\end{subfigure}
		\begin{subfigure}[t]{0.235\textwidth}
			\centering
			\includegraphics[width=0.99\textwidth,keepaspectratio]{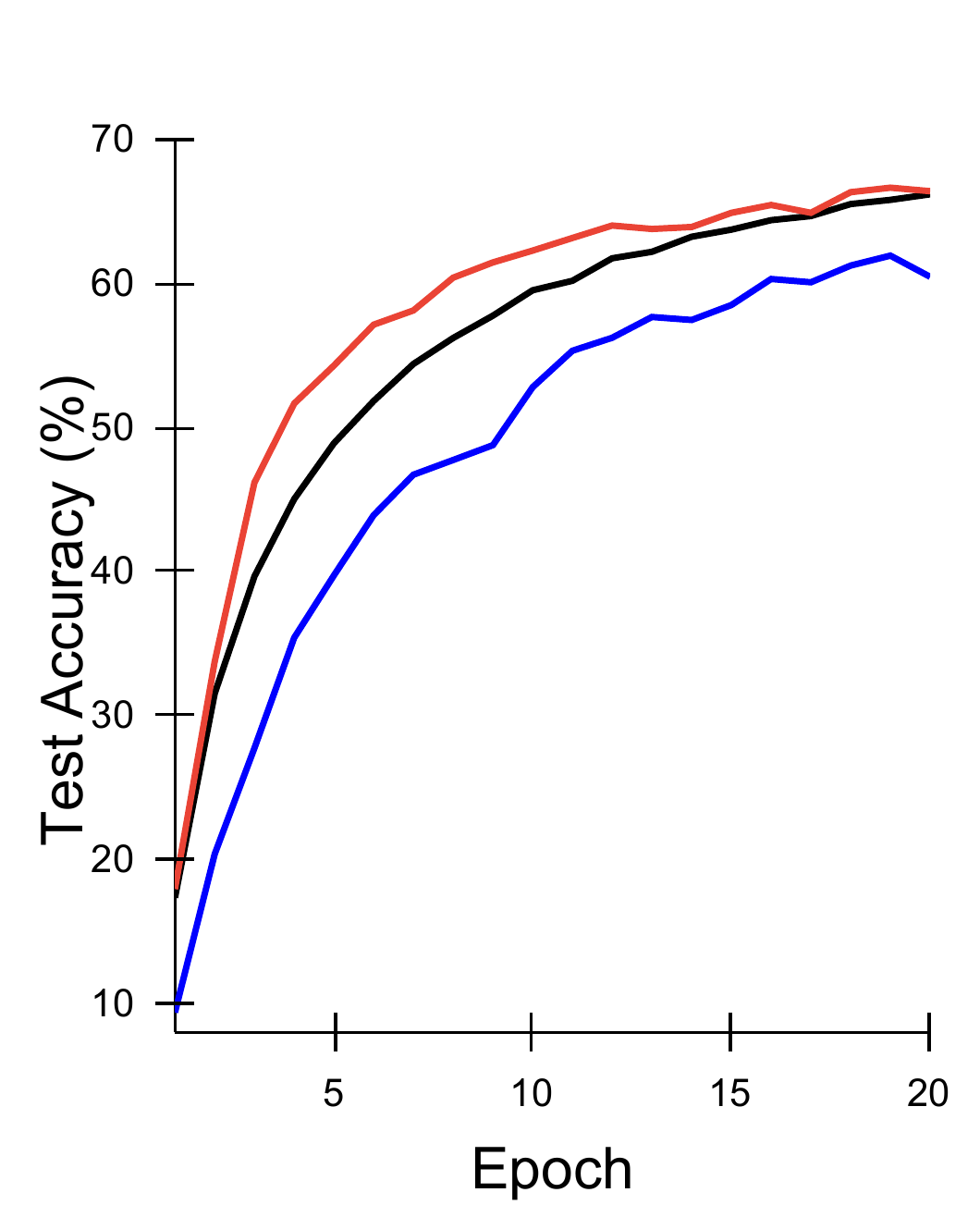}
			\caption{ResNet152 (Test Acc.)}
		\end{subfigure}
		\begin{subfigure}[t]{0.235\textwidth}
			\centering
			\includegraphics[width=0.99\textwidth,keepaspectratio]{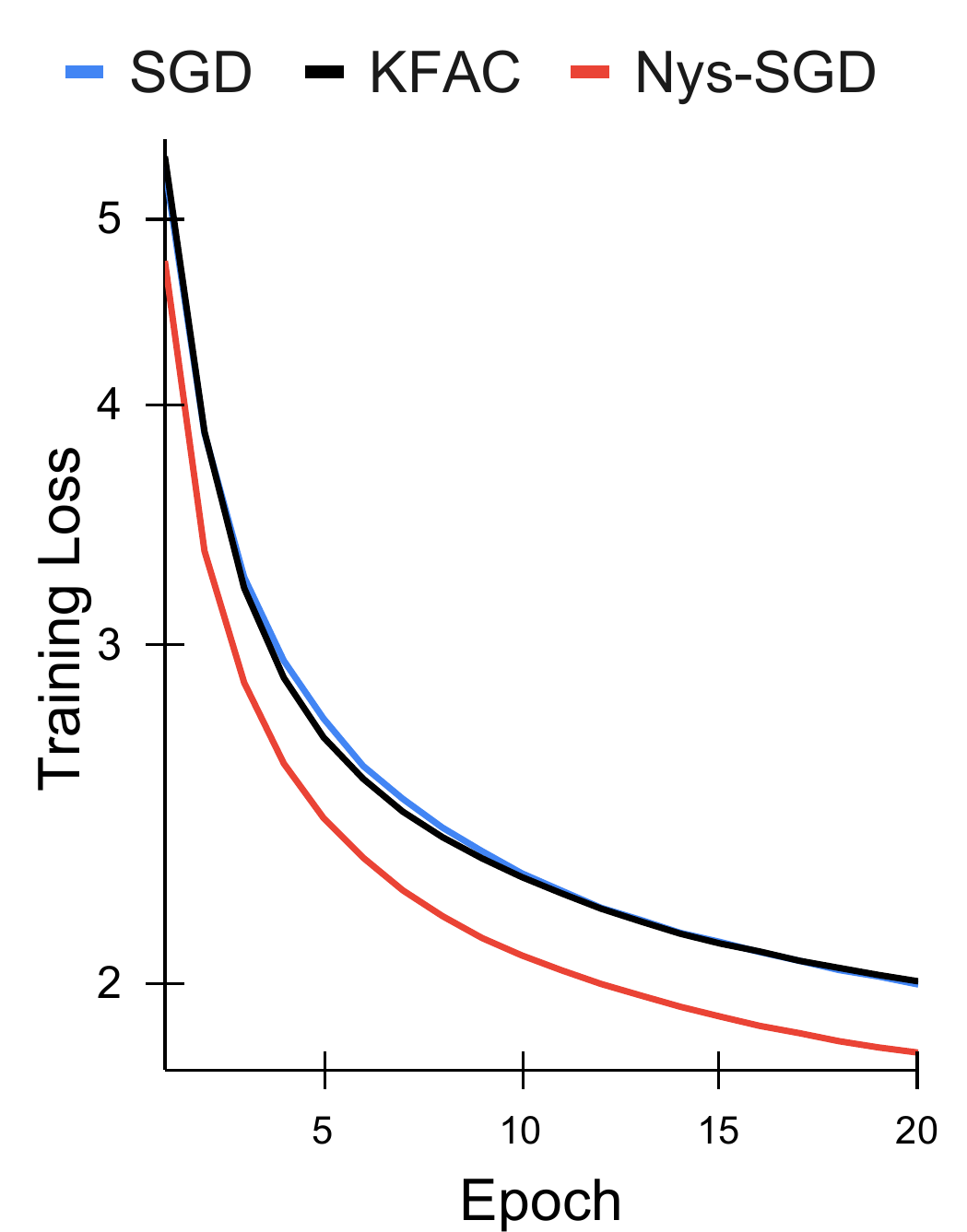}
			\caption{EfficientNet (Train loss)}
		\end{subfigure}
		\begin{subfigure}[t]{0.235\textwidth}
			\centering
			\includegraphics[width=0.99\textwidth,keepaspectratio]{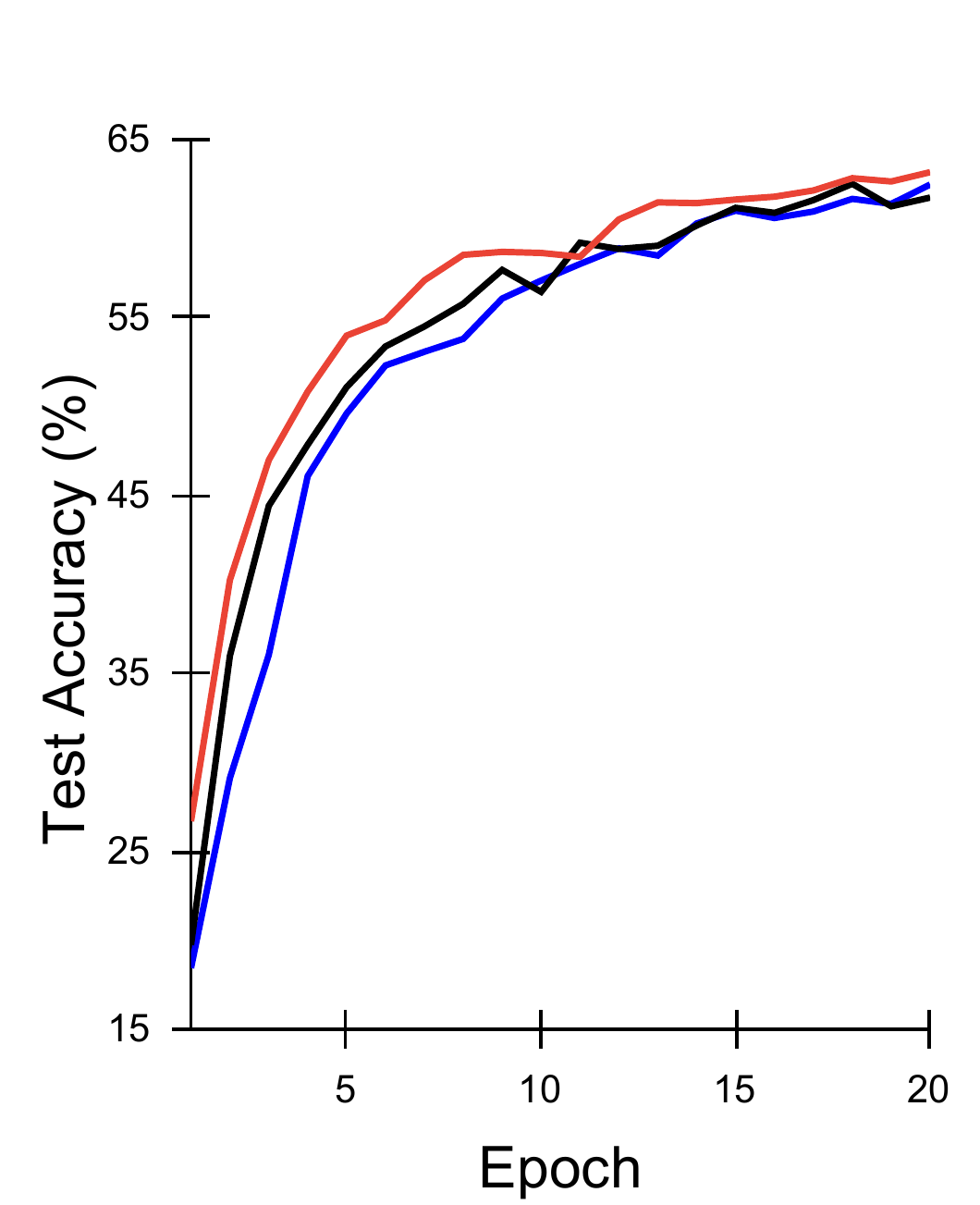}
			\caption{EfficientNet (Test Acc.)}
		\end{subfigure}
		\caption{Results on Imagenet using ResNet152 and EfficientNet, respectively. \label{fig:imagenet-resnet152}}
	\end{figure*}
	
	\noindent {\bf Experimental Setup:} We compared the proposed methods with existing \textit{state-of-the-art} first and second order optimization methods, namely, Adam, SVRG-LBFGS~\cite{Kolte2015accelerating}, oBFGS~\cite{Schraudolph2007stochastic}, SVRG-SQN~\cite{Moritz2016linearly}, and SQN~\cite{Byrd2016stochastic}. We compare various variants of the proposed methods, namely, Nystr\"om SGD (Nys-SGD) and Nystr\"om SVRG (Nys-SVRG) for optimization error and classification accuracy on both training and test data.

	\begin{wrapfigure}{r}{0.38\textwidth}
		\centering
		\includegraphics[width=0.8\linewidth,keepaspectratio]{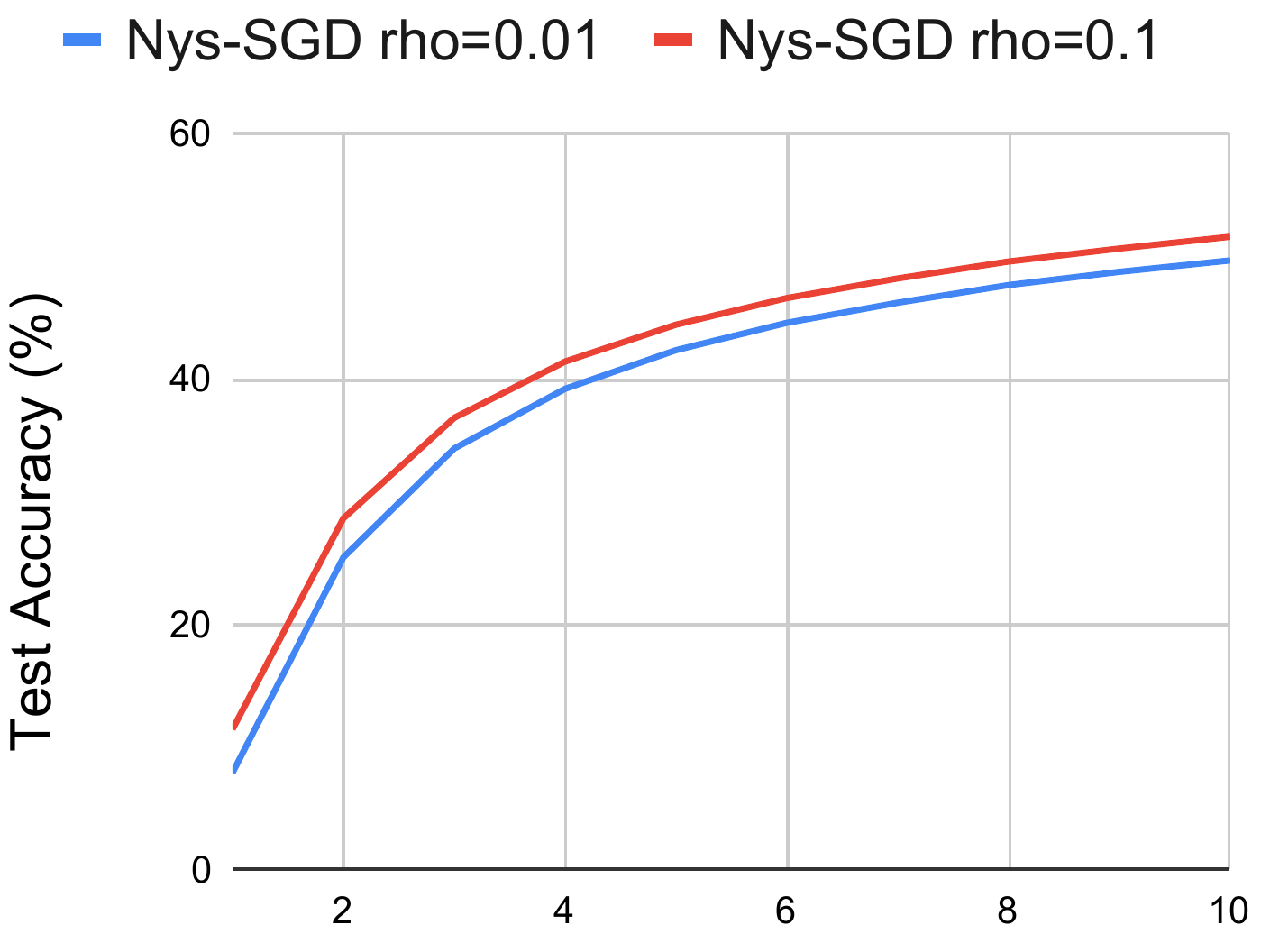}
		\caption{Effect of the $\rho$ on the test accuracy for ResNet18 on imagenet dataset. $x$-axis is the number of epochs.}
		\label{fig:diff_rho}
	\end{wrapfigure}	

	We demonstrate the performance of the proposed and existing methods on the $\ell_2$-regularized logistic regression problem with the regularizer $\lambda\in\{10^{-2},10^{-3},10^{-4}\}$. For each $\lambda$, all methods were tuned for other hyperparameters such as $\eta\in\{10^{0},10^{-1},10^{-2},10^{-3},10^{-4},10^{-5}\}$, batch size, epochs, etc. The  parameter $m=50$ and $\rho\in\{10^{0},10^{-1},10^{-2},10^{-3}\}$ were specific to the proposed variants only. 
	The memory used in the quasi-Newton method was set to $20$, which is a commonly used value~\cite{Kolte2015accelerating,Byrd2016stochastic}. 
	We report the optimization error on the training set (Opt. Error) and testing set (test error) with respect to epochs and training CPU time cost per epoch. The best-performing model was selected based on the minimum optimization error on the training set and presented its corresponding validation error. We implemented the existing and proposed methods in MATLAB using the SGDLibrary~\cite{Kasai2017jmlr}.
	We computed the results on an Intel(R) Xeon(R) CPU E5-2698v4 @ 2.20GHz with 40 cores running MATLAB R2019a.

	\begin{wraptable}{r}{0.5\textwidth}
		\centering
		\caption{Per iteration computational time (seconds). [*For KFAC on EfficientNet with batch size 128 could not fit into memory]}
		\begin{tabular}{lllll}
			\hline
			Model & Method & Batch & Update Time & Hessian  \\
			\hline
			& SGD & 128 & 1.006 & -\\
			ResNet152 & KFAC & 128 & 1.064 & -\\
			& Nys-SGD & 128 & 1.060 & 0.643\\
			\hline
			& SGD & 128 & 0.341 & -\\
			EfficientNet & KFAC & 64* & 1.173 & -\\
			& Nys-SGD & 128 & 0.347 & 2.620\\
			\hline
		\end{tabular}
		\label{tab:comput-time}
	\end{wraptable}
	
	\noindent {\bf Results:} Figure~\ref{fig:realsim} show the comparisons of the proposed methods with Adam and existing quasi-Newton methods. We present the numerical results per epoch as well CPU time. It may be observed from the results that the proposed methods performed better than or very competitively with the existing methods consistently on all the datasets with the varying number of dimensions and samples. Additionally, the performance of the proposed methods was almost stable when increasing the regularization parameters, whereas that of the existing methods changed significantly. The existing methods such as SVRG-LBFGS and OBFGS sometimes performed very competitively to the proposed methods. However, their performance was not very stable as they could fluctuate significantly even after reaching very close to the optimal point.
	It may be observed from the results that the proposed methods generated robust curvature information using only 50 columns of Hessian in the Nystr\"om approximation, compared to the existing quasi-Newton methods using traditional secant equation. As we took the initial point $\boldw_0$ to be the solution of the least squares problem, the first-order method Adam was unable to significantly reduce the optimization cost. It may be observed that existing stochastic methods require long CPU times for problems with higher dimensionality. More precisely, for the \textit{real-sim} dataset, SQN and SVRG-SQN require a considerable period of time to reach the optimal point. (See additional results for logistic regression and $\ell_2$-svm in appendix)

	\subsection{Real-world experiment (Non-convex setup)}
	We also evaluated the performance of the Nystr\"om SGD on the well-known deep models on Imagenet dataset.
		
	\noindent {\bf Experimental Setup:}
	We compared our method with the first-order methods SGD and the well-known approximate second-order method KFAC on ResNet152 \cite{he2016deep} and EfficientNet \cite{tan2019efficientnet} models.
	For Nystr\"om SGD, we used $\rho=0.1$ and fixed the $m=\log_2|\boldw|$, where $|\boldw|$ is the number of parameters in the respective model. We used a batch size of 128. We used a random sample of size of $\min\{6400,n\times 0.01\}$ to compute the partial Hessian $\boldC$ for Nystr\"om SGD. The update frequency used to re-estimate the preconditioner in KFAC and its variants is set to $200$, as used in their experiments. The ImageNet results were computed on a Quadro RTX 8000 GPU.	
	
	\noindent {\bf Results:}
	Figure~\ref{fig:imagenet-resnet152} present the results of the ResNet152 and EfficientNet on the ImageNet dataset. The proposed method outperformed both the SGD and KFAC for both the models in terms of training loss as well as test accuracy, showing the better optimization and generalization ability of the trained models. Table \ref{tab:comput-time} shows the computational time comparison of methods. The per update computational time of the Nystr\"om SGD on ResNet152 is 1.703 seconds which is slightly slower than SGD and KFAC. To further speed up the  Nystr\"om SGD is an interesting future work. Figure \ref{fig:diff_rho} shows the effect of the $\rho$ parameter for ResNet18. As can be seen, the $\rho$ parameter affects the model performance. We found setting $\rho = 0.1$ performs well in practice.

	\section{Conclusion}
	In this paper, we propose to approximate the Hessian matrix using Nystr\"om method  and to use the approximated Hessian for stochastic optimization. In terms of optimal cost, the proposed Nystr\"om based methods performed better than the first-order methods Adam, SGD, and SVRG, and performed better than or comparably with the stochastic quasi-Newton methods SVRG-LBFGS, oBFGS, SVRG-SQN, and SQN. In terms of computational time, our methods performed compares favorably with the stochastic quasi-Newton methods. Additionally, experimental results on several publicly available benchmark datasets show that the proposed method exhibited much more stable and generalized behavior in comparison to both first- and second-order methods, which can be attributed to the better preservation of the curvature information through Nystr\"om approximation of the Hessian matrix.
	
	\section*{Acknowledgement}
	We would like to thank Mohammad Emtiyaz Khan for helpful discussion.
	
	\bibliography{ref}
	\bibliographystyle{plain}

	\newpage
	\appendix
	\onecolumn
	
	\section{Proofs}
	
	\subsection{Proof of Lemma~\ref{lem:lowerboundNystrom}}
	\begin{proof}
		Recall that $|\|\boldA\| - \|\boldB\| | \leq \|\boldA - \boldB\|$ for the matrix $\boldA,\boldB$ of the same size. Let $rank(\boldH) = r$, and its singular values is given by $\sigma_1\geq \sigma_2 \geq \ldots \geq \sigma_r > \sigma_{r+1} = \ldots = \sigma_d = 0.$ Then the difference of the Hessian $\boldH$ with its best $k$-rank approximation is $\|\boldH - \boldH_k\|_2 = \sigma_{k+1}(\boldH)$. Applying these facts in Nystr\"om error bound~\eqref{eq:Nystromerrorbound}, we then get:
		\begin{align*}
			|\|\boldN\| - \|\boldH\|| & \leq \|\boldH - \boldN\| \leq \sigma_{k+1} + \epsilon \sum_{i=1}^{d} \boldH^2_{ii},\\
			\|\boldN\| & \leq \|\boldH\| + \sigma_{k+1} + \epsilon \sum_{i=1}^{d} \boldH^2_{ii},\\
			\sigma_{\max}(\boldN) & \leq \sigma_{r}(\boldH) + \sigma_{k+1}(\boldH) + \epsilon \sum_{i=1}^{d} \boldH^2_{ii}
		\end{align*}
		Let $\Gamma = (\sigma_{r}(\boldH) + \sigma_{k+1}(\boldH) + \epsilon \sum_{i=1}^{d} \boldH^2_{ii})^2,$ then $\lambda_{\max}(\boldN) \leq \Gamma.$
		Finally,
		\begin{equation}\label{asm:approx_hess_bound}
			\boldzero \preceq \boldN \preceq \Gamma \boldI.    
		\end{equation}
	\end{proof}
	\subsection{Proof of Lemma~\ref{lem:Nystrom_bound}}
	\begin{proof} \rm
		Let $\lambda_1(\boldN),\ldots,\lambda_d(\boldN)$\footnote[2]{For a matrix $X$, $\lambda(X) $ denotes an eigenvalue of $X$.}
		are the eigenvalues of the Nystrom approximation 
		$\boldZ\boldZ^\top = \boldN$ 
		in the decreasing order; i.e., 
		$\lambda_{i}(\boldN) \geq \lambda_{i+1}(\boldN)$, for 
		$i \in \{1,\ldots,d-1\}$. 
		From the~\eqref{asm:approx_hess_bound}, it is clear that 
		$0 = \lambda_{min}(\boldN)$ and $\lambda_{max} \leq \Gamma $.
		
		Now, using Weyl's theorem~\cite[Theorem 3.2.1]{Bhatia2013matrix} for the eigenvalues of two symmetric matrices, we prove bound on $\boldN + \rho \boldI $, for  $\rho > 0$.
		\begin{align*}
			\lambda_i(\boldN) + \lambda_d(\rho \boldI) \leq
			\lambda_i(\boldN + \rho \boldI) \leq \lambda_i(\boldN) + 
			\lambda_1(\rho \boldI), \text{~ for } i \in \{1,\ldots,d\},
		\end{align*}
		which implies,
		\begin{equation*}
			\lambda_{min}(\boldN) + \rho \leq \lambda_i(\boldN + \rho I)
			\leq \lambda_{max}(\boldN) + \rho, \text{~ for } i \in \{1,\ldots,d\}.
		\end{equation*}
		Therefore, 
		\begin{equation*}
			\rho \boldI \preceq (\boldN +\rho \boldI) 
			\preceq (\Gamma + \rho ) \boldI.
		\end{equation*}
	\end{proof}
	
	\subsection{Proof of Lemma~\ref{lem:InverseNB}}
	\begin{proof}
		From the Lemma~\ref{lem:Nystrom_bound}, it is easy to see that,
		\begin{equation*}
			\frac{1}{( \Gamma + \rho )}\boldI
			\preceq (\boldN +\rho \boldI)^{-1} 
			\preceq \frac{1}{ \rho } \boldI.
		\end{equation*}
		Therefore, 
		\begin{equation*}
			\frac{1}{( \Gamma + \rho )}\boldI
			\preceq \boldB 
			\preceq \frac{1}{ \rho } \boldI.
		\end{equation*}
	\end{proof}
	\vfill
	
	\subsection{Proof of Lemma~\ref{lem:gradient_bound}}
	\begin{proof}
		From the strong convexity of $f$,\\
		\begin{align}
			f(\boldz)\geq & ~f(\boldw) + \nabla f(\boldw)^\top(\boldz - \boldw)+\frac{\mu}{2} \|\boldz - \boldw\|^2 \nonumber\\
			\geq & ~f(\boldw) + \nabla f(\boldw)^\top\left(-\frac{1}{\mu}\nabla f(\boldw)\right) + \frac{\mu}{2}\left\|\frac{1}{\mu}\nabla f(\boldw)\right\|^2\nonumber \\
			= & ~f(\boldw) - \frac{1}{2\mu} \|\nabla f(\boldw) \|^2, \label{eq:grad_bound}
		\end{align}
		where the last inequality holds from the minimizer $\boldz = \boldw - \frac{1}{\mu} \nabla f(\boldw)$ of quadratic model:
		\begin{equation*}
			q(\boldz) = f(\boldw) + \nabla f(\boldw)^\top(\boldz -\boldw) + \frac{\mu}{2}\| \boldz -\boldw\|^2.
		\end{equation*}
		By setting $\boldz = \boldw_* $ in \eqref{eq:grad_bound} gives
		\begin{equation*}
			\|\nabla f(\boldw) \|^2 \geq 2\mu (f(\boldw) - f(\boldw_*)).
		\end{equation*}
	\end{proof}
	
	\subsection{Proof of Theorem~\ref{thmNys-SVRG}}
	\begin{proof}
		Using Assumption $\ref{asm:Lipschitz}$, we get
		\begin{align*}
			f(\boldw_t) \leq &~ f(\boldw_{t-1}) + \nabla f(\boldw_{t-1})^\top(\boldw_t - \boldw_{t-1}) + 
			\frac{\Lambda}{2} \| \boldw_t - \boldw_{t-1} \|^2 \\
			= &~ f(\boldw_{t-1}) - \eta \nabla f(\boldw_{t-1})^\top \boldB_{\tau} \boldv_{t-1} + \frac{\eta^2 \Lambda}{2} 
			\|\boldB_{\tau} \boldv_{t-1} \|^2. \\
		\end{align*}
		Taking expectation on $f$,
		\begin{equation*}
			\mathbb{E} [f(\boldw_t)] \leq f(\boldw_{t-1}) - \eta \nabla 
			f(\boldw_{t-1})^\top \boldB_{\tau} \nabla f(\boldw_{t-1}) + \frac{\eta^2 \Lambda }{2} \mathbb{E}\|\boldB_{\tau} \boldv_{t-1} \|^2
		\end{equation*}
		From Lemma~\ref{lem:InverseNB}, we get
		\begin{align}
			\mathbb{E} [f(\boldw_t)] \leq &~ f(\boldw_{t-1}) - \frac{\eta}{(\Gamma + \rho)} \|\nabla 
			f(\boldw_{t-1})\|^2 + \frac{\eta^2 \Lambda }{2(\rho)^2} \mathbb{E}\|\boldv_{t-1} \|^2\\
			\leq &~ f(\boldw_{t-1}) - \frac{2\eta \mu}{(\Gamma + \rho)}(f(\boldw_{t-1}) - f(\boldw_*)) \\
			& \qquad \qquad \qquad \qquad +\frac{2 \eta^2 \Lambda^2 }{\rho^2}
			(f(\boldw_{t-1}) - f(\boldw_*) + f(\boldw_{\tau -1}) - f(\boldw_*)) \nonumber
		\end{align}
		where the last inequality follows from the Lemma~\ref{lem:gradient_bound} and Lemma~\ref{lem:expectationbound} for the second and the third term, respectively.
		Now, taking sum over $t = 1,\ldots,\ell$, we get
		\begin{align*}
			\mathbb{E}[f(\boldw_{\ell})] \leq &~ \mathbb{E}[f(\boldw_0)] + 
			\frac{2\ell \eta^2 \Lambda^2}{(\rho )^2} \mathbb{E}[f(\boldw_{\tau -1}) - f(\boldw_*)]\\
			&~ \qquad \qquad \qquad \qquad - 2\eta \left( \frac{\mu}{(\Gamma + \rho)} - 
			\frac{\eta \Lambda^2}{\rho^2}\right) 
			\left[ \sum_{t = 1}^{\ell} f(\boldw_{t-1}) - \ell f(\boldw_*)\right]\\
			\leq &~ \mathbb{E}[f(\boldw_{\tau -1})] + \frac{2\ell \eta^2 \Lambda^2}{\rho^2}
			\mathbb{E}[f(\boldw_{\tau -1}) - f(\boldw_*)] \\
			&~\qquad \qquad \qquad \qquad - 2\eta \ell \left( \frac{\mu}{(\Gamma + \rho)} - 
			\frac{\eta \Lambda^2}{\rho^2}\right) 
			\mathbb{E}[ f(\boldw_{\tau}) -  f(\boldw_*)]
		\end{align*}
		Rearranging above inequality gives
		\begin{align*}
			0 \leq &~ \mathbb{E}[f(\boldw_{\tau -1 }) - f(\boldw_{\ell})] +    \frac{2\ell \eta^2 
				\Lambda^2}{\rho^2} \mathbb{E}[f(\boldw_{\tau -1}) - f(\boldw_*)] \\
			&~\qquad \qquad \qquad \qquad - 2\eta \ell \left( \frac{\mu}{(\Gamma + \rho)} - 
			\frac{\eta \Lambda^2}{\rho^2}\right) \mathbb{E}[ f(\boldw_{\tau}) -  f(\boldw_*)]\\
			\leq &~ \mathbb{E}[f(\boldw_{\tau -1 }) - f(\boldw_*)] +    \frac{2\ell \eta^2 
				\Lambda^2}{(\rho)^2} \mathbb{E}[f(\boldw_{\tau -1}) - f(\boldw_*)] \\
			&~\qquad \qquad \qquad \qquad - 2\eta \ell \left( \frac{\mu}{(\Gamma + \rho)} - 
			\frac{\eta \Lambda^2}{\rho^2}\right) \mathbb{E}[ f(\boldw_{\tau}) -  f(\boldw_*)]\\
			\leq &~ \left( 1 + \frac{2\ell \eta^2 
				\Lambda^2}{(\rho)^2} \right) \mathbb{E}[f(\boldw_{\tau -1}) - f(\boldw_*)] \\
			&~\qquad \qquad \qquad \qquad - 2\eta \ell \left( \frac{\mu}{(\Gamma + \rho)} - 
			\frac{\eta \Lambda^2}{\rho^2}\right) \mathbb{E}[ f(\boldw_{\tau}) -  f(\boldw_*)]
		\end{align*}
		The second inequality follows from the fact that $f(\boldw_*) \leq f(\boldw_{\ell})$. From the assumption
		of $\eta < \mu \Delta /2\Lambda^2 \delta^2$, which implies
		\begin{equation}
			\mathbb{E}[f(\boldw_{\tau}) - f(\boldw_*)] \leq \frac{1 + 2 \ell \eta^2 \Lambda^2 \delta^2}
			{2\ell \eta (\mu \Delta - \eta \Lambda^2 \delta^2)} \mathbb{E}[f(\boldw_{\tau -1 }) - f(\boldw_*)]
		\end{equation}
		where,
		$$\Delta = \frac{1}{(\Gamma + \rho )},\quad \delta = \frac{1}{ \rho} .$$
		Choosing $\ell$ and $\eta$ to satisfy \eqref{eq:l_bound}, it follows that the $\alpha < 1$.
	\end{proof}
	\subsection{Proof of Theorem~\ref{Th:New_linearconvergence_svrg}}
	\begin{proof}
		\begin{align}
			& \E \|\boldw_{t} -   \boldw_*\|^2 \nonumber\\
			& = \E \|\boldw_{t-1} - \eta \boldB_{\tau}\boldv_{t-1} - \boldw_*\|^2 \nonumber\\
			& = \|\boldw_{t-1} - \boldw_*\|^2 - 2 \eta~\E((\boldw_{t-1} - \boldw_*)^\top \boldB_{\tau}\boldv_{t-1}) + \eta^2~\E\|\boldB_{\tau}\boldv_{t-1}\|^2 
			\nonumber\\
			& = \|\boldw_{t-1} - \boldw_*\|^2 - 2 \eta (\boldw_{t-1} - \boldw_*)^\top\boldB_{\tau} \nabla f(\boldw_{t-1}) + \eta^2 \E\|\boldB_{\tau}\boldv_{t-1}\|^2 \nonumber\\
			& \leq \|\boldw_{t-1} - \boldw_*\|^2  - \frac{2\eta \mu}{(\Gamma+\rho)} \|\boldw_{t-1} - \boldw_*\|^2   + \frac{\eta^2}{\rho}\mathbb{E}\| \boldv_{t-1}\|^2\nonumber\\
			& \leq \|\boldw_{t-1} - \boldw_*\|^2  - \frac{2\eta\mu}{(\Gamma+\rho)} \|\boldw_{t-1} - \boldw_*\|^2 + \frac{4\eta^2 \Lambda}{\rho} (f(\boldw_{t-1}) - f(\boldw_*) + f(\boldw_{\tau -1}) - f(\boldw_*)).\nonumber\\
			& \leq \|\boldw_{t-1} - \boldw_*\|^2 - 2\eta\mu\Delta \|\boldw_{t-1} - \boldw_*\|^2 + 4\eta^2\Lambda\delta \left[ \frac{\Lambda}{2} \|\boldw_{t-1} - \boldw_*\|^2 + \frac{\Lambda}{2} \|\boldw_{\tau-1} \boldw_*\|^2 \right] \nonumber\\
			& = \|\boldw_{t-1} - \boldw_*\|^2 - 2\eta\mu\Delta \|\boldw_{t-1} - \boldw_*\|^2 + 2\eta^2\Lambda^2\delta \left[  \|\boldw_{t-1} - \boldw_*\|^2 + \|\boldw_{\tau-1} \boldw_*\|^2 \right] \nonumber\\
			& = \|\boldw_{t-1} - \boldw_*\|^2 - 2\eta\mu\Delta \|\boldw_{t-1} - \boldw_*\|^2 + 2\eta^2\Lambda^2\delta \left[  \|\boldw_{t-1} - \boldw_*\|^2 + \|\boldw_{\tau-1} \boldw_*\|^2 \right] \nonumber\\
			& = (1 - 2 \eta (\mu \Delta - \eta \Lambda^2 \delta))\|\boldw_{t-1} - \boldw_*\|^2 + 2\eta^2\Lambda^2\delta  \|\boldw_{\tau-1} \boldw_*\|^2  \nonumber
		\end{align}
		where first inequality comes from the bound on $\boldB_{\tau}$ and second inequality obtained by the strong convexity of $f$ and substituting the upper bound of $\E\|v_{t-1}\|^2$ 
		from~Lemma~\ref{lem:expectationbound} and
		$$\Delta = \frac{1}{(\Gamma + \rho )},\quad \delta = \frac{1}{ \rho} .$$
		The third inequality comes from the Lipschtiz continuity of $f$.
		Now applying the above inequality over the $t$ and
		${\boldw}_{\tau-1} = \boldw_0$ and ${\boldw}_{\tau} = \boldw_{\ell}$,
		\begin{align*}
			\E \|\boldw_{\tau} - \boldw_*\|^2 & \leq \left[ 1 - 2 \eta (\mu \Delta - \eta \Lambda^2 \delta)\right]^{\ell} \|\boldw_{\tau-1} - \boldw_*\|^2 \\
			& \qquad \qquad \qquad \qquad + {2\eta^2\Lambda^2\delta} ~\sum_{j=1}^{\ell} \left[ 1 - 2 \eta (\mu \Delta - \eta \Lambda^2 \delta) \right]^{j} \|\boldw_{\tau - 1} - \boldw_*\|^2\\
			& < \left[ 1 - 2 \eta (\mu \Delta - \eta \Lambda^2 \delta)\right]^{\ell} \|\boldw_{\tau-1} - \boldw_*\|^2 + \frac{\eta \Lambda^2 \delta}{(\mu \Delta - \eta \Lambda^2 \delta)} \|\boldw_{\tau - 1} - \boldw_*\|^2
		\end{align*}
		\begin{equation}
			\E \|\boldw_{\tau} - \boldw_*\|^2  <  \zeta ~\|\boldw_{\tau -1} - \boldw_*\|^2.
		\end{equation}
		where 
		\begin{equation}
			\zeta = \left[ \left( 1 - 2 \eta (\mu \Delta - \eta \Lambda^2 \delta)\right)^{\ell}
			+ \frac{\eta \Lambda^2 \delta}{ (\mu \Delta - 2\eta \Lambda^2 \delta)} \right].
		\end{equation}
	\end{proof}
	
	\clearpage
	
	\section{Relation to Nystr\"om logistic regression}
	\label{sec:nys-lr}
	Here, we describe the relation to the Nystr\"om logistic regression \cite[Section 3.2.1]{talwalkar2010matrix}. 
	
	Let $\{(\boldx_i,y_i)\}_{i=1}^n$ be given $n$ training samples, where $\boldx_i \in \mathbbR^d$ and $y \in \{0,1\}$. Then, the optimization problem of the regularized logistic regression is given as
	\begin{align*}
		\min_{\boldw} & \hspace{.3cm} -\sum_{i = 1}^{n}y_i \log \sigma(\boldw^\top \boldx_i) + (1-y_i)\log (1 - \sigma(\boldw^\top \boldx_i)) + \frac{\lambda}{2}\|\boldw\|_2^2,
	\end{align*}
	where $\sigma(a) = \frac{1}{1 + \exp(-a)}$ is the sigmoid function. The Hessian of the regularized logistic regression can be given as
	\begin{align*}
		\boldH = \boldX \boldD \boldX^\top + \lambda \boldI_d,
	\end{align*}
	where $\boldX = [\boldx_1, \boldx_2, \ldots, \boldx_n] \in \mathbbR^{d \times n}$ and $\boldD \in \mathbbR^{n \times n}$ is the diagonal matrix whose $i$-th diagonal element is $\sigma(\boldw^\top \boldx_i)(1-\sigma(\boldw^\top \boldx_i))$. Thus, the Newton step of the regularized logistic regression is given as
	\begin{align*}
		\boldw_\tau = \boldw_{\tau-1} - (\boldX \boldD \boldX^\top + \lambda \boldI_d)^{-1}\boldv_{\tau-1},
	\end{align*}
	where $\boldv_t$ denotes the gradient at $\tau$-th epoch.
	
	For the kernel logistic regression, the optimization problem can be given as
	\begin{align*}
		\min_{\boldw} & \hspace{.3cm} -\sum_{i = 1}^{n}y_i \log \sigma(\boldw^\top \boldk(\boldx_i)) + (1-y_i)\log (1 - \sigma(\boldw^\top \boldk(\boldx_i))) + \frac{\lambda}{2}\|\boldw\|_2^2,
	\end{align*}
	where $\boldk(\boldx) = (K(\boldx_1, \boldx), K(\boldx_2, \boldx), \ldots, K(\boldx_n, \boldx))^\top \in \mathbbR^{n}$ and $\boldw \in \mathbbR^{n}$. Note that we use the $\ell_2$-regularization of $\boldw$. Thus, the Hessian of kernel logistic regression can be given as
	\begin{align*}
		\boldH = \boldK \boldD \boldK^\top + \lambda \boldI_n,
	\end{align*}
	where $\boldK = [\boldk(\boldx_1), \boldk(\boldx_2), \ldots, \boldk(\boldx_n)] \in \mathbbR^{n \times n}$ is the Gram matrix. 
	
	Then, the Newton step of kernel logistic regression is given as
	\begin{align*}
		\boldw_\tau = \boldw_{\tau-1} - (\boldK \boldD \boldK^\top + \lambda \boldI_n)^{-1}\boldv_{\tau-1}.
	\end{align*}
	\cite{talwalkar2010matrix} proposed the Nystr\"om logistic regression algorithm, where the Nystr\"om method is used to approximate the Hessian of the regularized logistic regression. Thus, it can be regarded as a variant of Nystr\"om-SGD. However, \cite{talwalkar2010matrix} only considered the regularized logistic regression, in which the Hessian can be explicitly obtained, with deterministic optimization. In contrast, we propose the Nystr\"om method for a general Hessian matrix for stochastic optimization and show its efficacy for deep learning models. Moreover, we elucidate the theoretical properties of the Nystr\"om-SGD algorithm for convex functions.
	
	\clearpage
	\section{Experimental setup}

	\begin{table}[!h]
		\small
		\centering
		\caption{Details of the datasets used in the experiments}
		\begin{tabular}{lrrr}
			\hline
			\textbf{Dataset} & \textbf{Dim} & \textbf{Train} & \textbf{Test} \\
			\hline
			\textit{adult}\footnotemark[1] & 123 & 32,561 & 16,281 \\
			\textit{mnist}\footnotemark[2]     & 784 & 60,000 & 10,000 \\
			\textit{cifar10}\footnotemark[3]     & 3,072 & 50,000 & 10,000 \\
			\textit{real-sim}\footnotemark[1]    & 20,958 & 57,909 & 14,400 \\
			\textit{w8a}\footnotemark[1]    & 300 & 49,749 & 14,951 \\
			\hline
		\end{tabular}
		\label{tab:datasets}
		\vspace{-.2in}
	\end{table}
	\footnotetext[1]{Available at LIBSVM~\cite{chang2011libsvm} \url{https://www.csie.ntu.edu.tw/cjlin/libsvm/}}
	\footnotetext[2]{The original 10 classes in the MNIST are converted to the binary classes based on the round (0,3,6,8,9) vs. non-round (1,2,4,5,7) digits, similar to~\cite{Lan2019tnnls,Zhang2012aistats}}
	\footnotetext[3]{The original 10 classes in the CIFAR10 are converted to the binary classes based on the natural (bird, cat, deer, dog, frog, horse) vs. man-made (airplane, automobile, ship, truck) objects, similar to~\cite{Hsieh2014icml}}
	
	\begin{table}[!h]
		\small
		\centering
		\caption{Details of the datasets used in the deep learning experiments}
		\begin{tabular}{lrrrr}
			\hline
			\textbf{Dataset} & \textbf{Dim} & \textbf{Train} & \textbf{Test} \\
			\hline
			\textit{imagenet}      & $224\times 224\times 3$ & 1.2M & 50,000 & 1000\\
			\hline
		\end{tabular}
		\label{tab:deepdatasets}
		\vspace{-.2in}
	\end{table}	
	
	\clearpage
	
	\section{Additional numerical experiments on \textit{mnist}}
	\begin{figure}[!h]
		\includegraphics[width=1\linewidth,keepaspectratio]{ICML_convex_figures/LEGEND_CONVEX.png}\\
		\includegraphics[width=.32\linewidth,height=4cm]{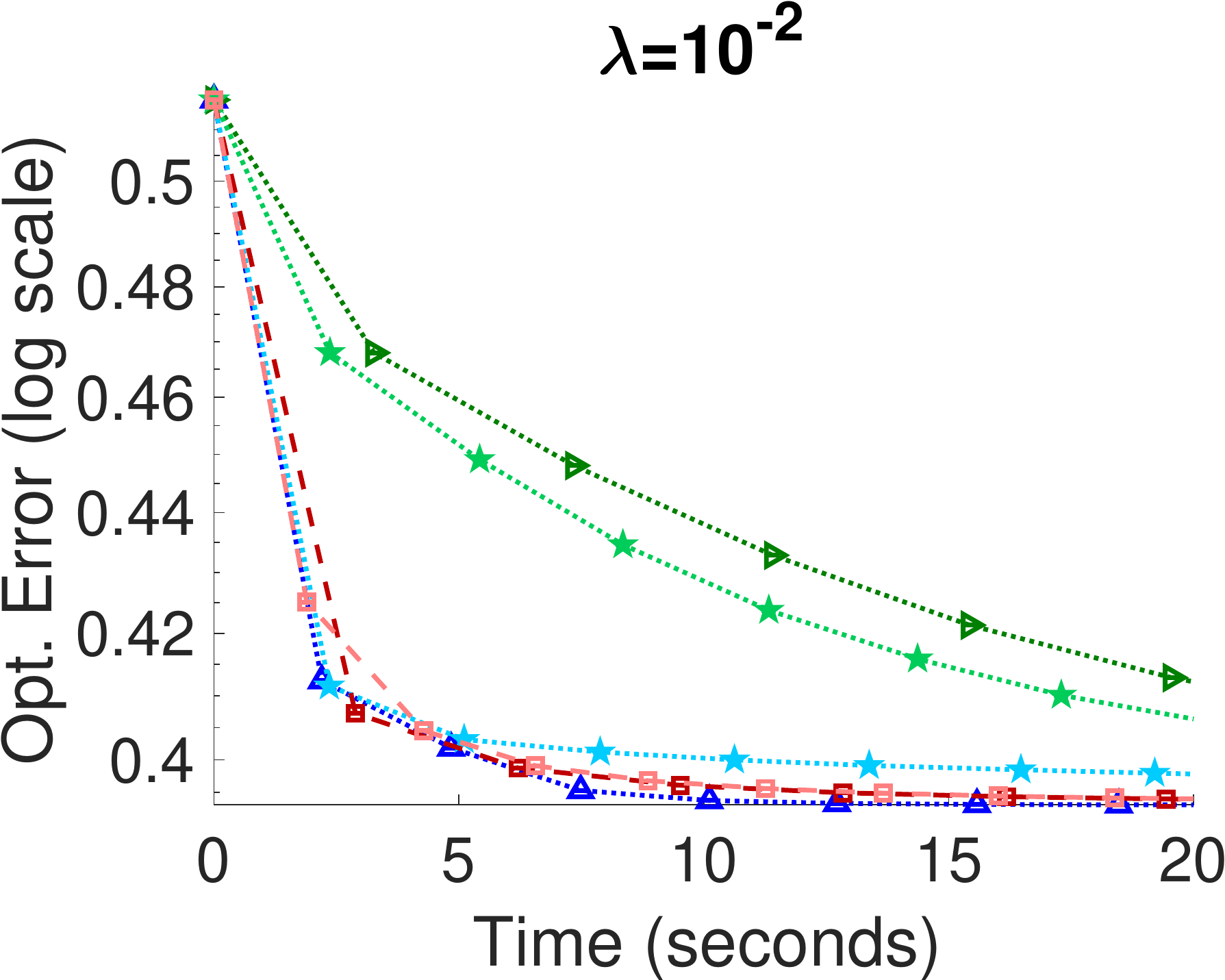} 
		\includegraphics[width=.32\linewidth,height=4cm]{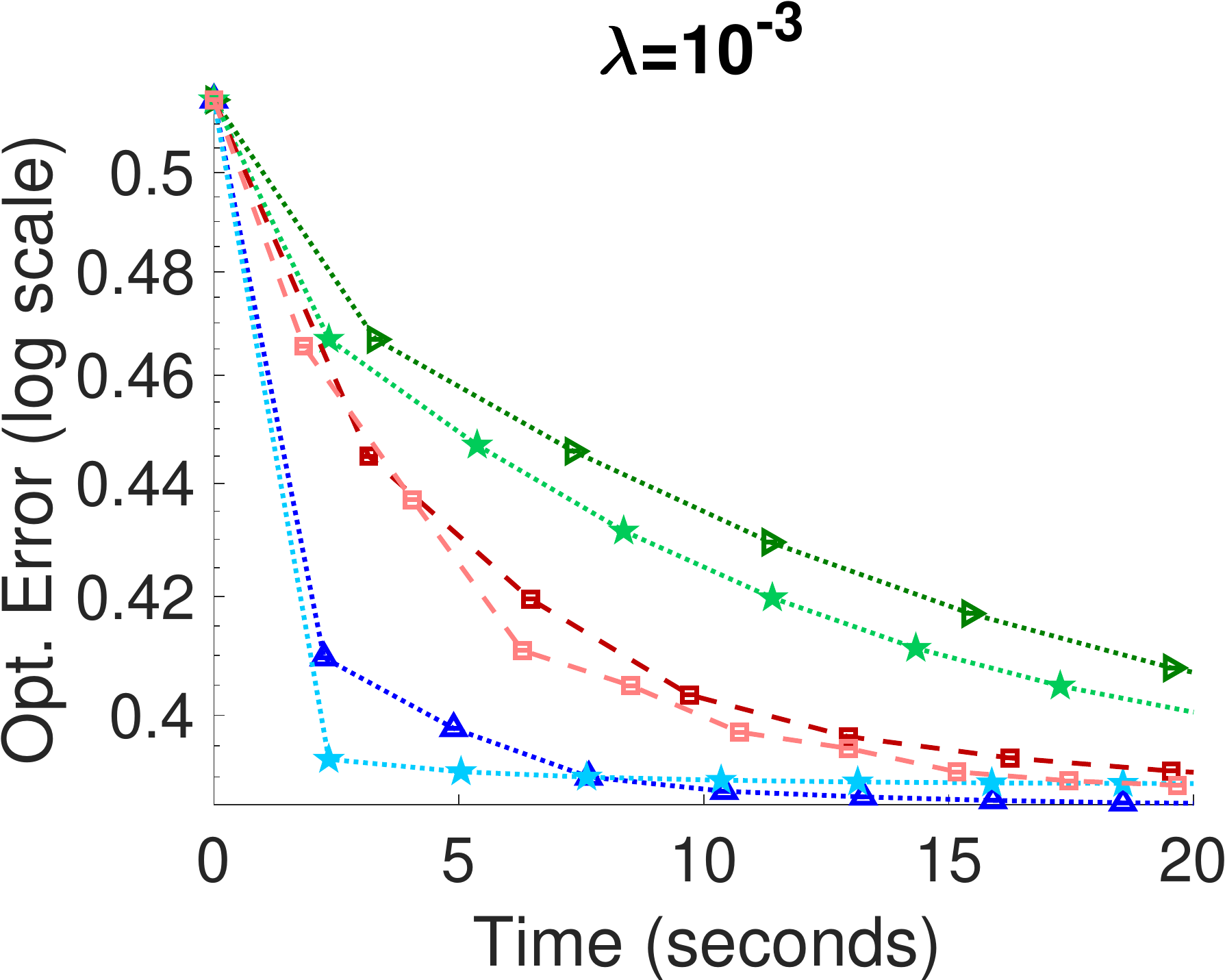} 
		\includegraphics[width=0.32\linewidth,height=4cm]{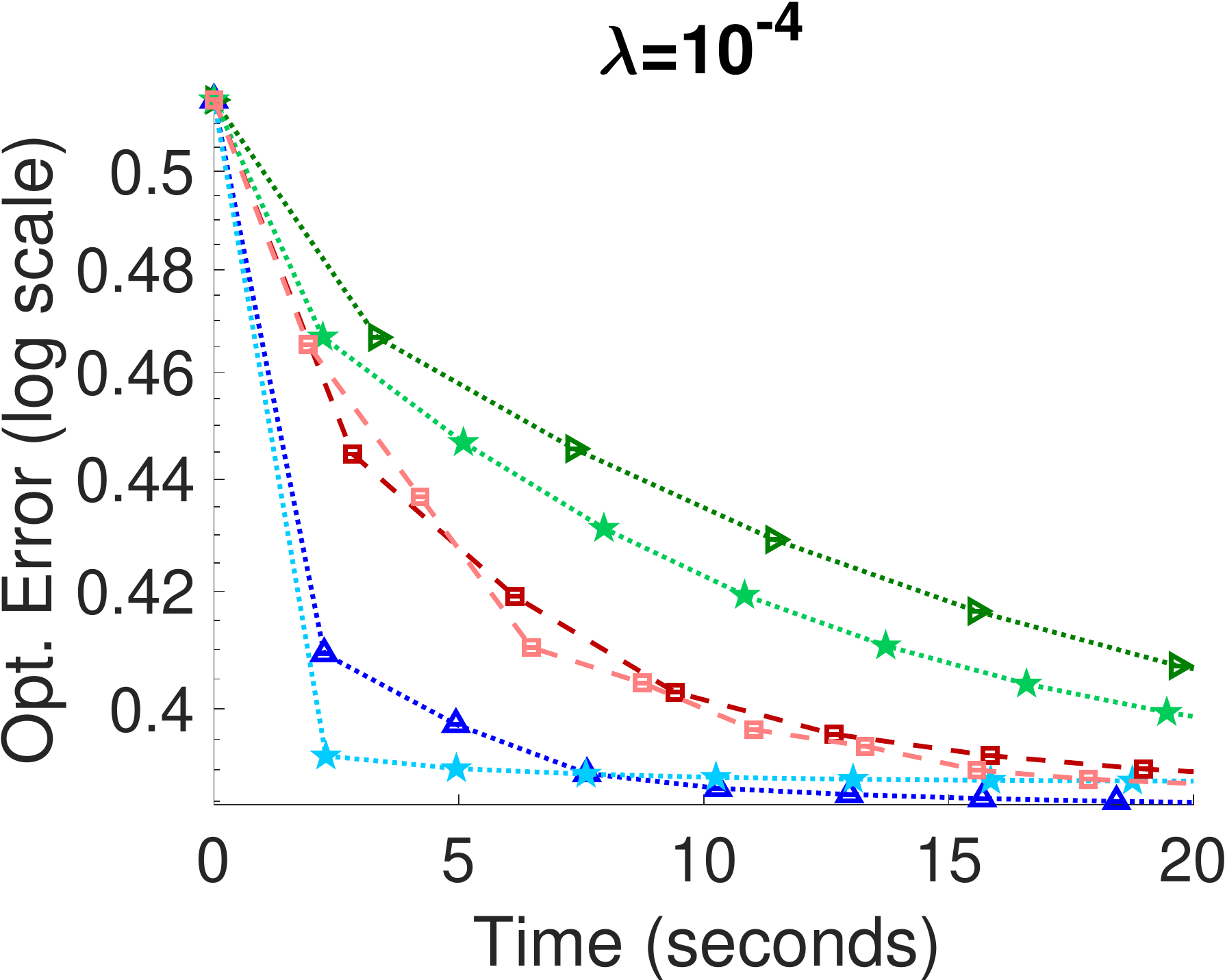} \\
		\includegraphics[width=.32\linewidth,height=4cm]{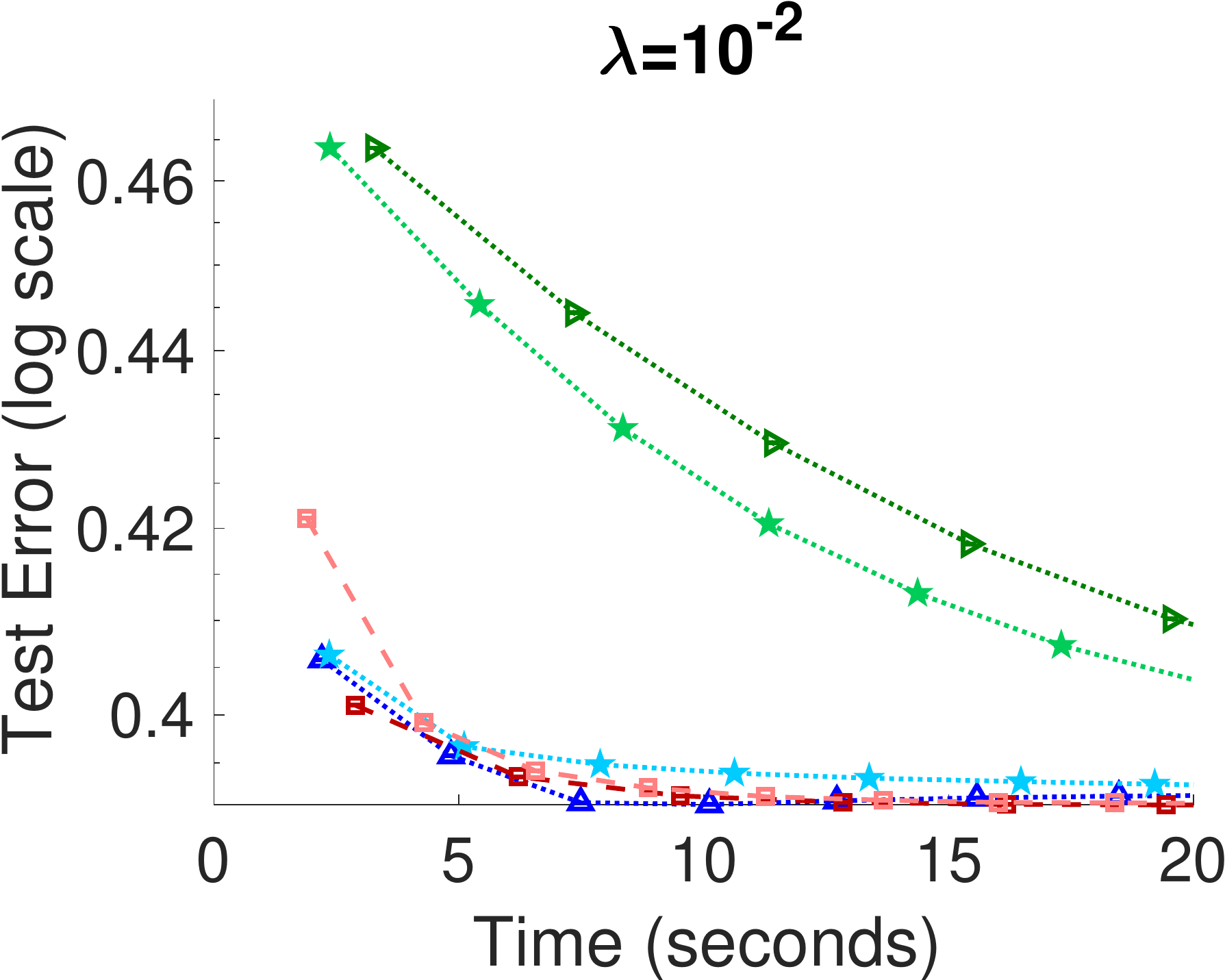} 
		\includegraphics[width=.32\linewidth,height=4cm]{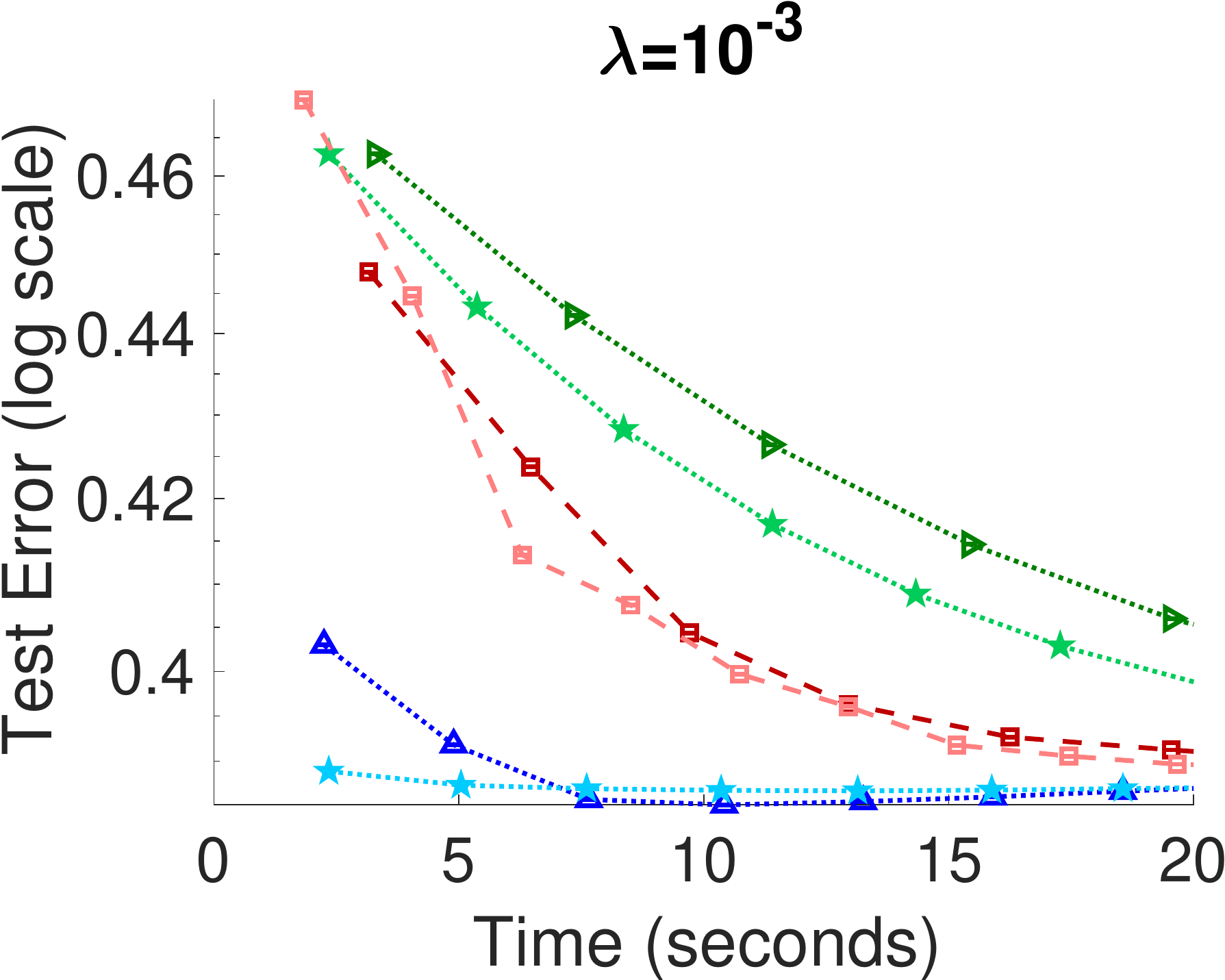} 
		\includegraphics[width=0.32\linewidth,height=4cm]{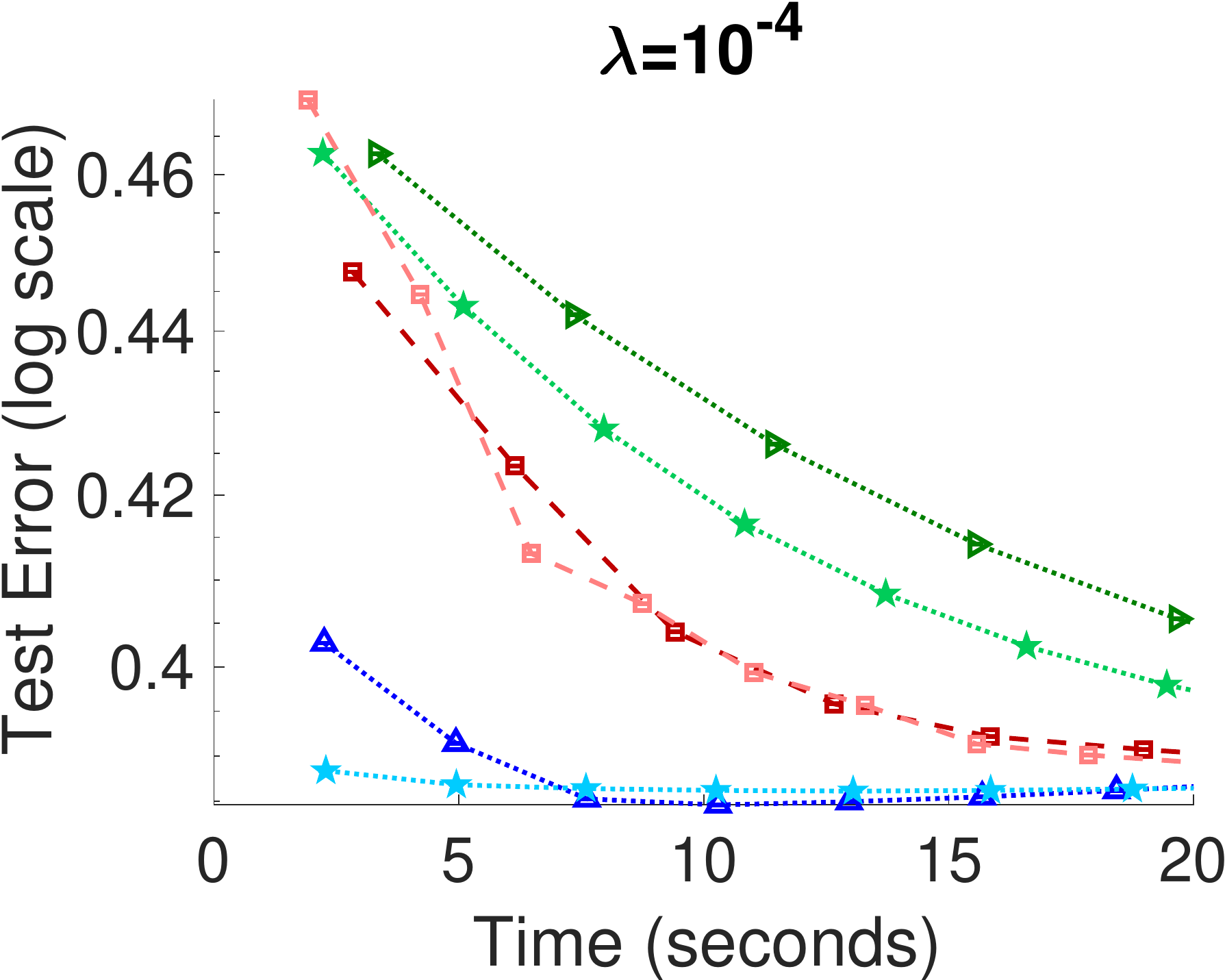} \\
		\includegraphics[width=.32\linewidth,height=4cm]{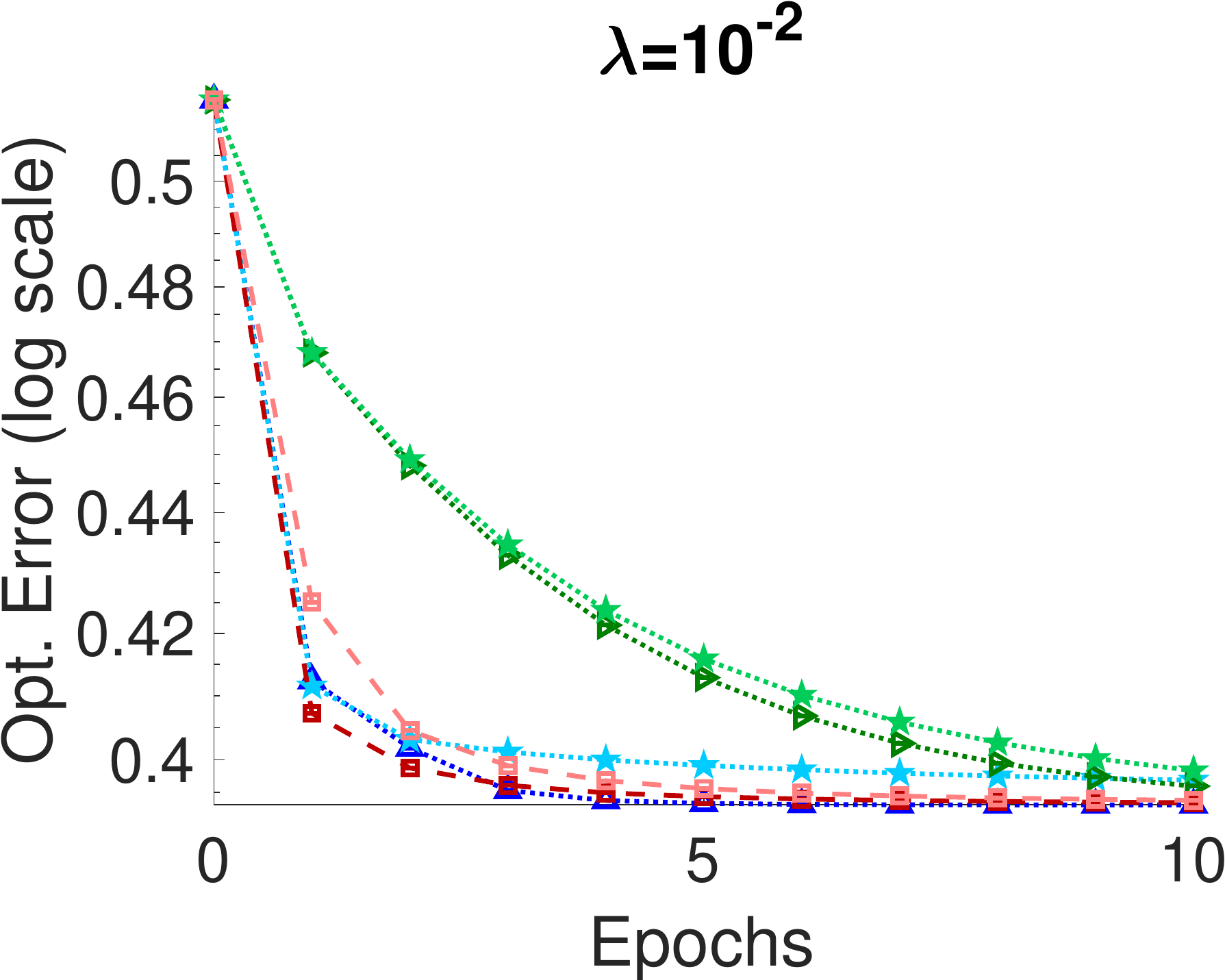} 
		\includegraphics[width=.32\linewidth,height=4cm]{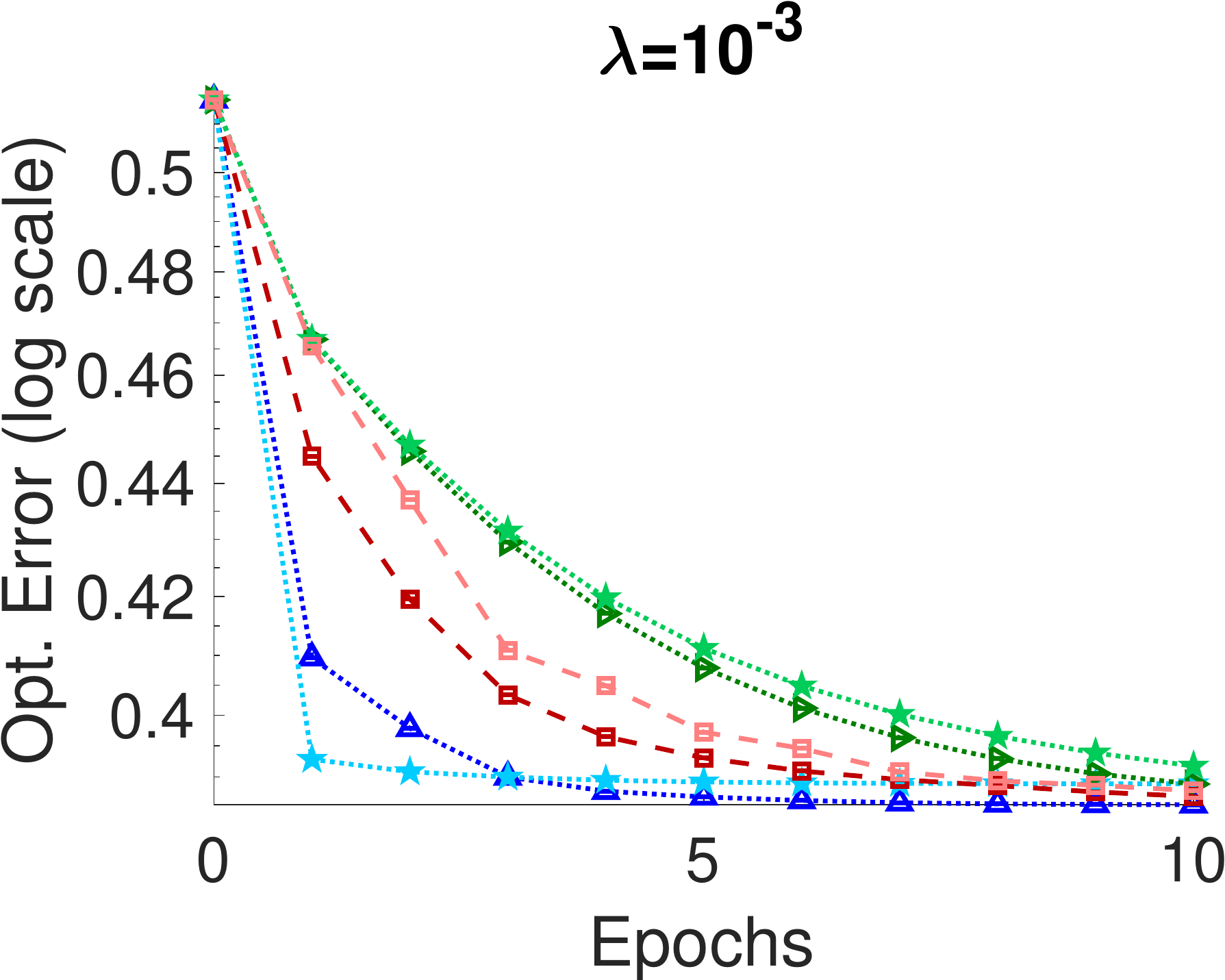} 
		\includegraphics[width=0.32\linewidth,height=4cm]{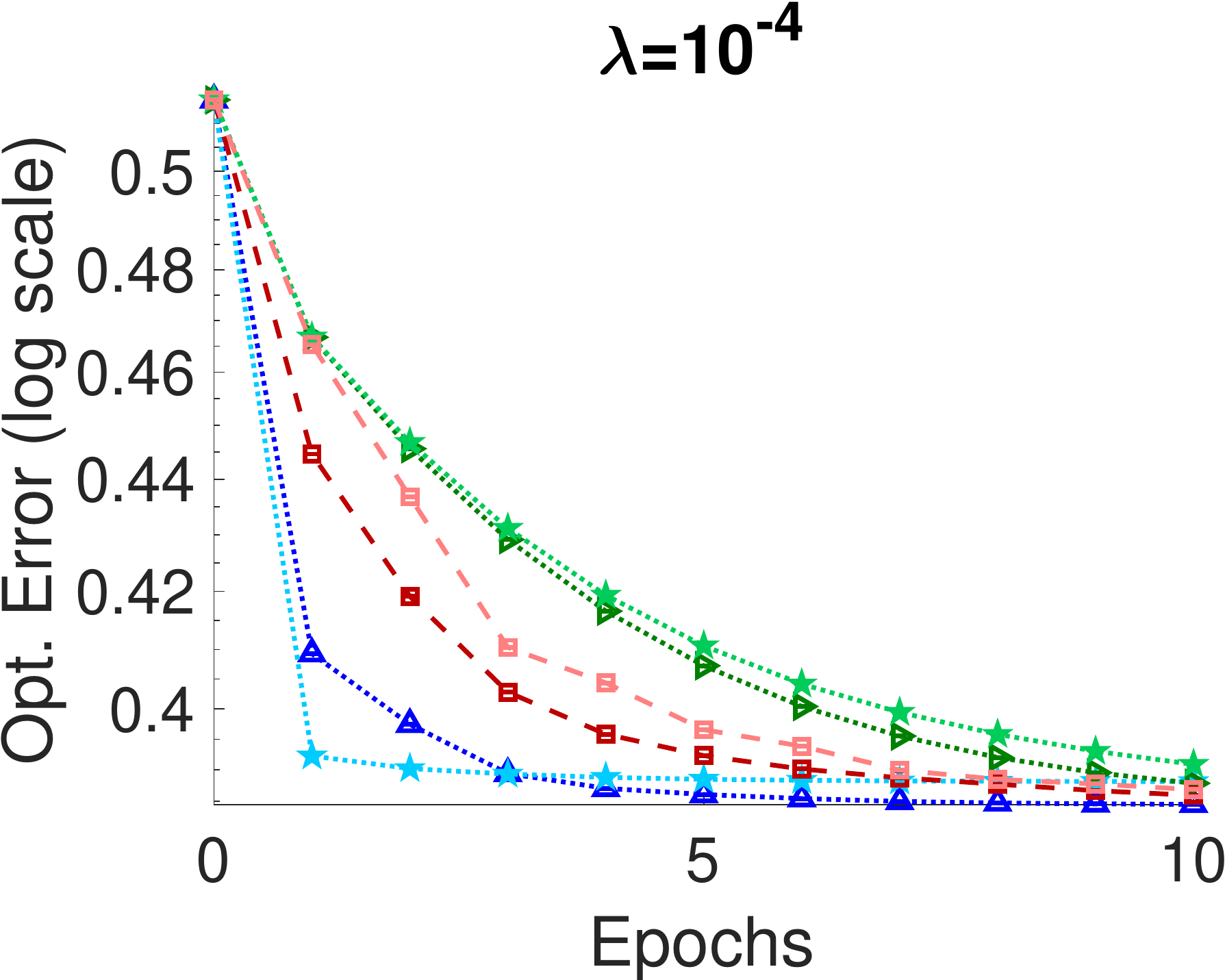} \\
		\includegraphics[width=.32\linewidth,height=4cm]{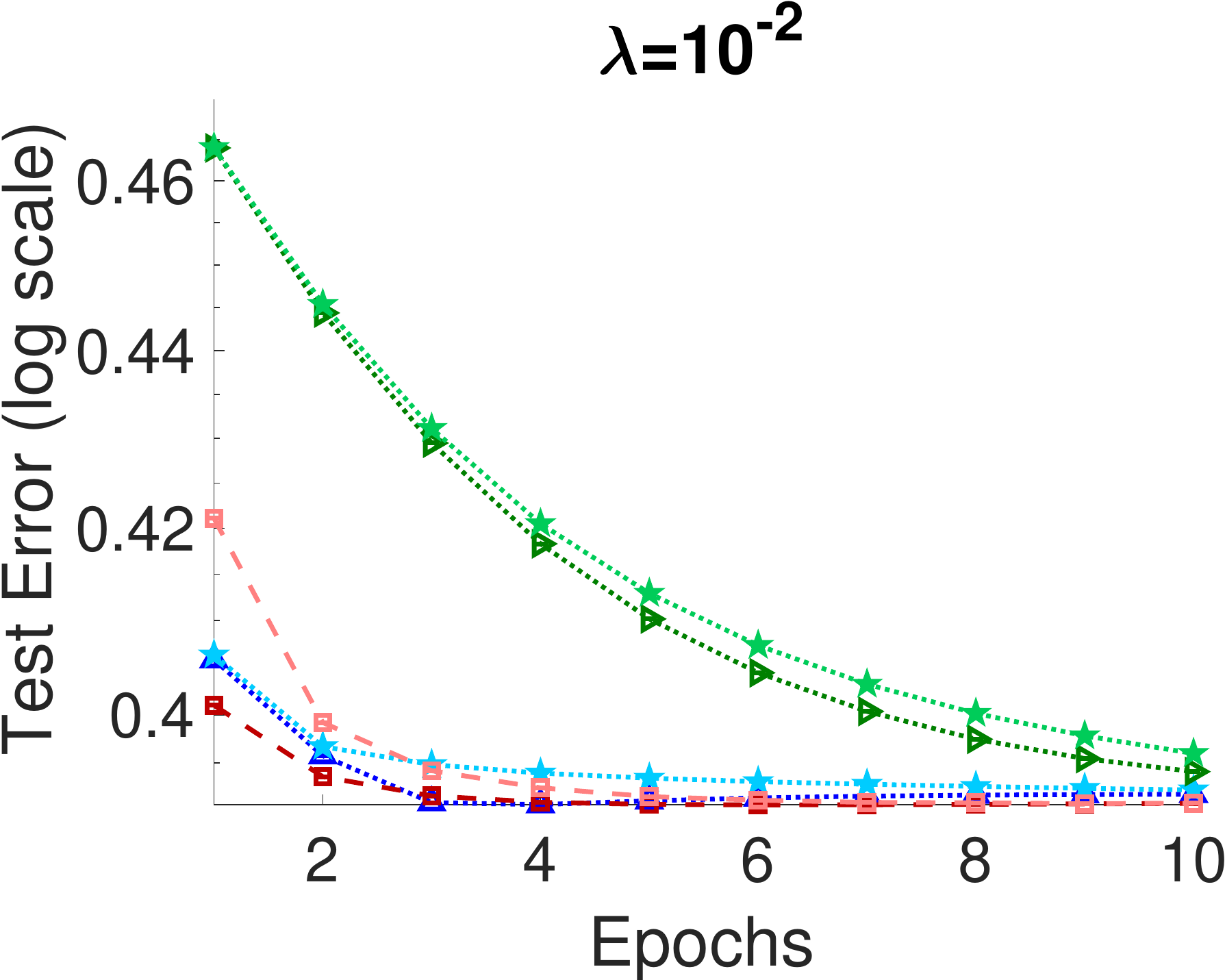} 
		\includegraphics[width=.32\linewidth,height=4cm]{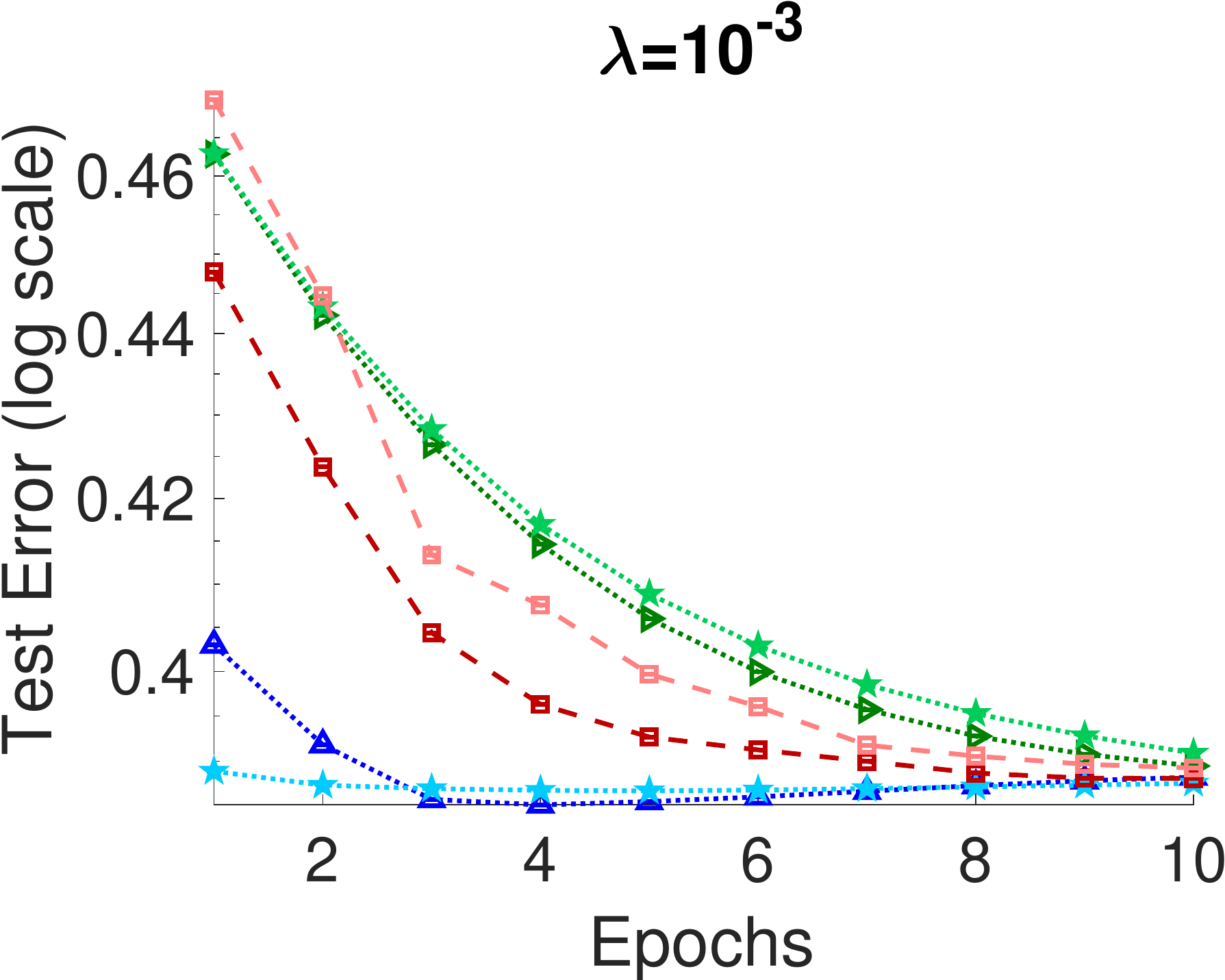} 
		\includegraphics[width=0.32\linewidth,height=4cm]{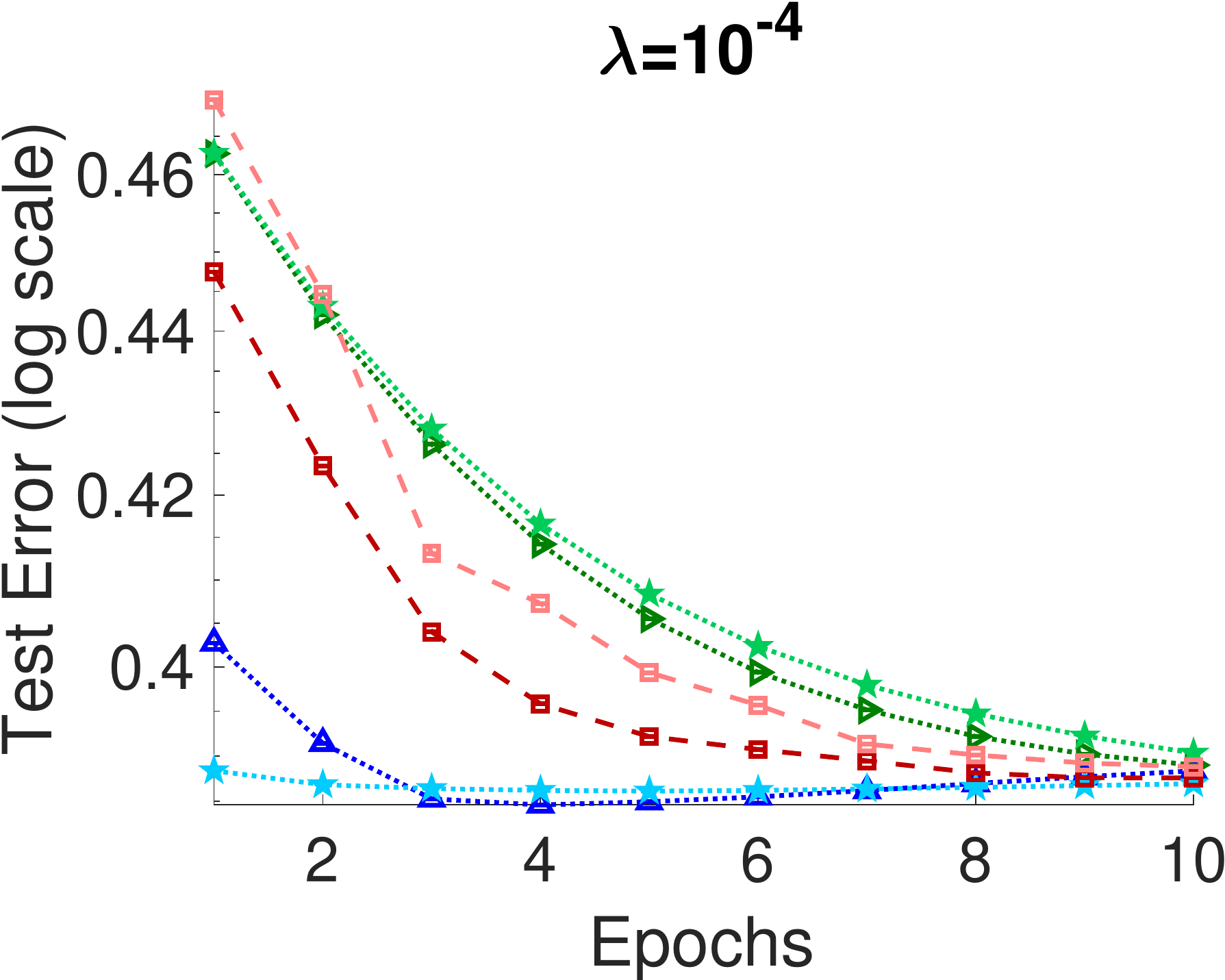}
		\caption{Comparison of proposed methods with the existing optimization methods for optimization performances achieved within 20 seconds on \textit{mnist} dataset with various regularizer $\lambda$ for \textit{logistic regression}.}
		\label{fig:mnist}
	\end{figure}
	
	\clearpage
	
	\section{Additional numerical experiments on \textit{cifer}}
	
	\begin{figure}[!h]
		\includegraphics[width=1\linewidth,keepaspectratio]{ICML_convex_figures/LEGEND_CONVEX.png}\\
		\includegraphics[width=.32\linewidth,height=4cm]{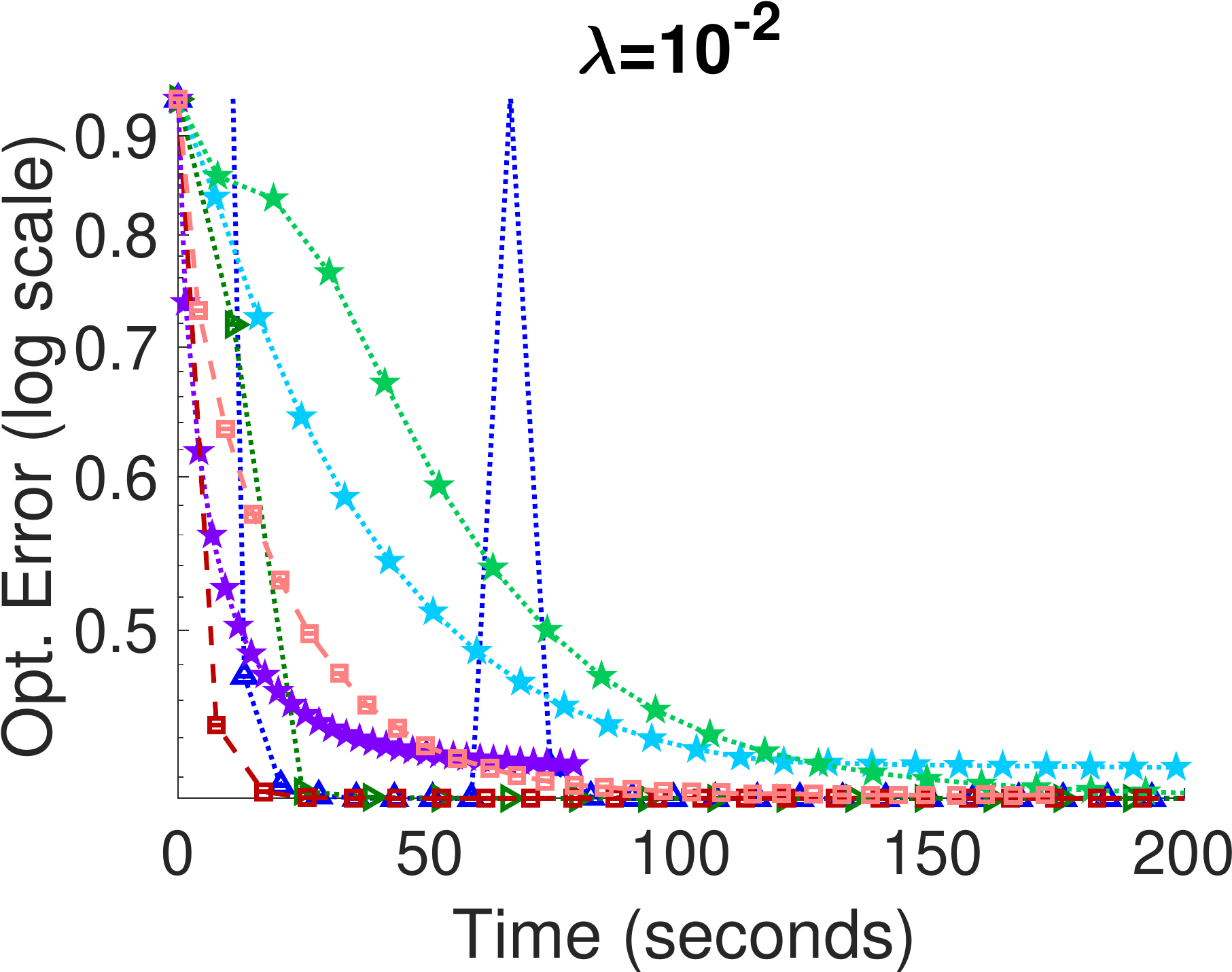} 
		\includegraphics[width=.32\linewidth,height=4cm]{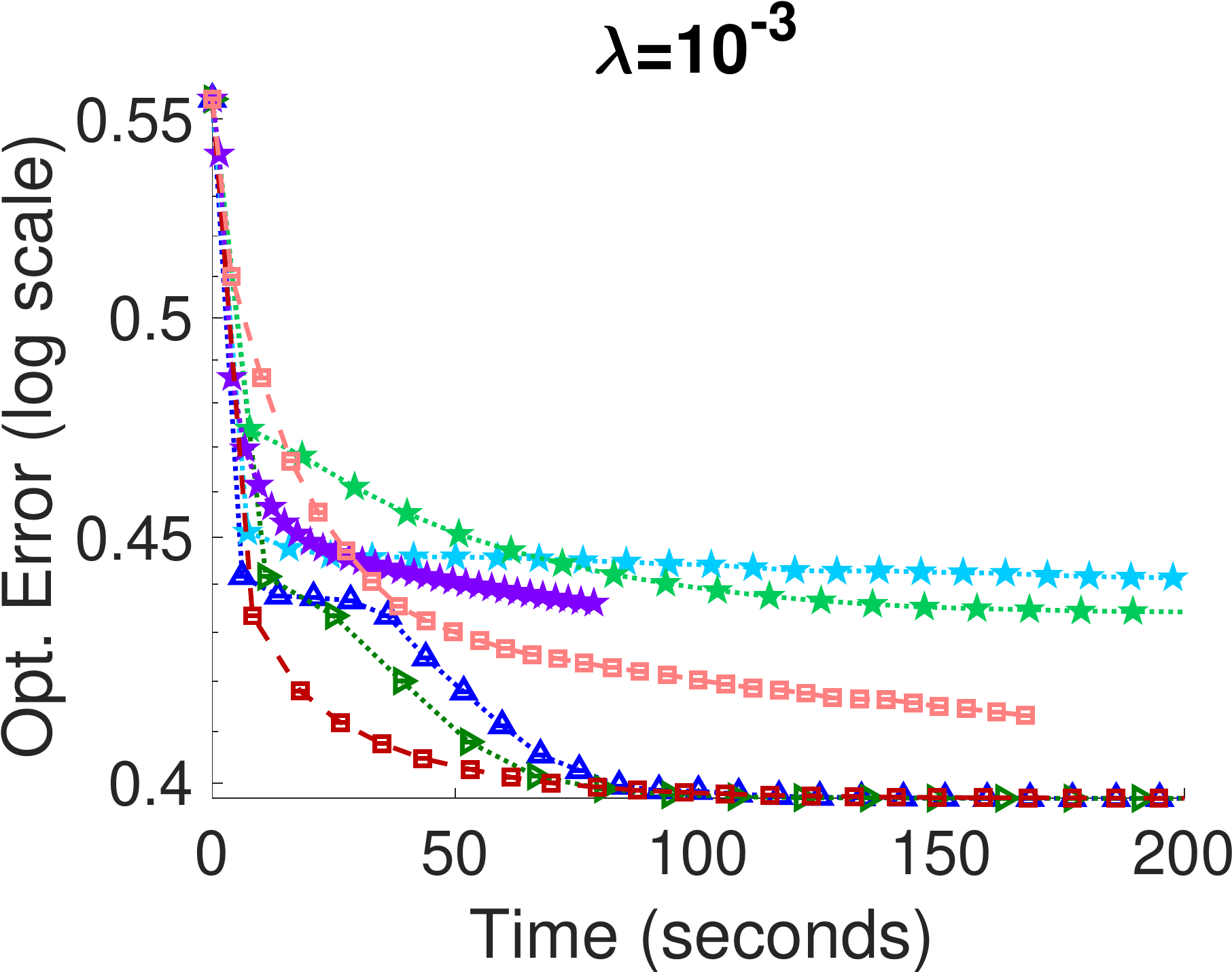} 
		\includegraphics[width=.32\linewidth,height=4cm]{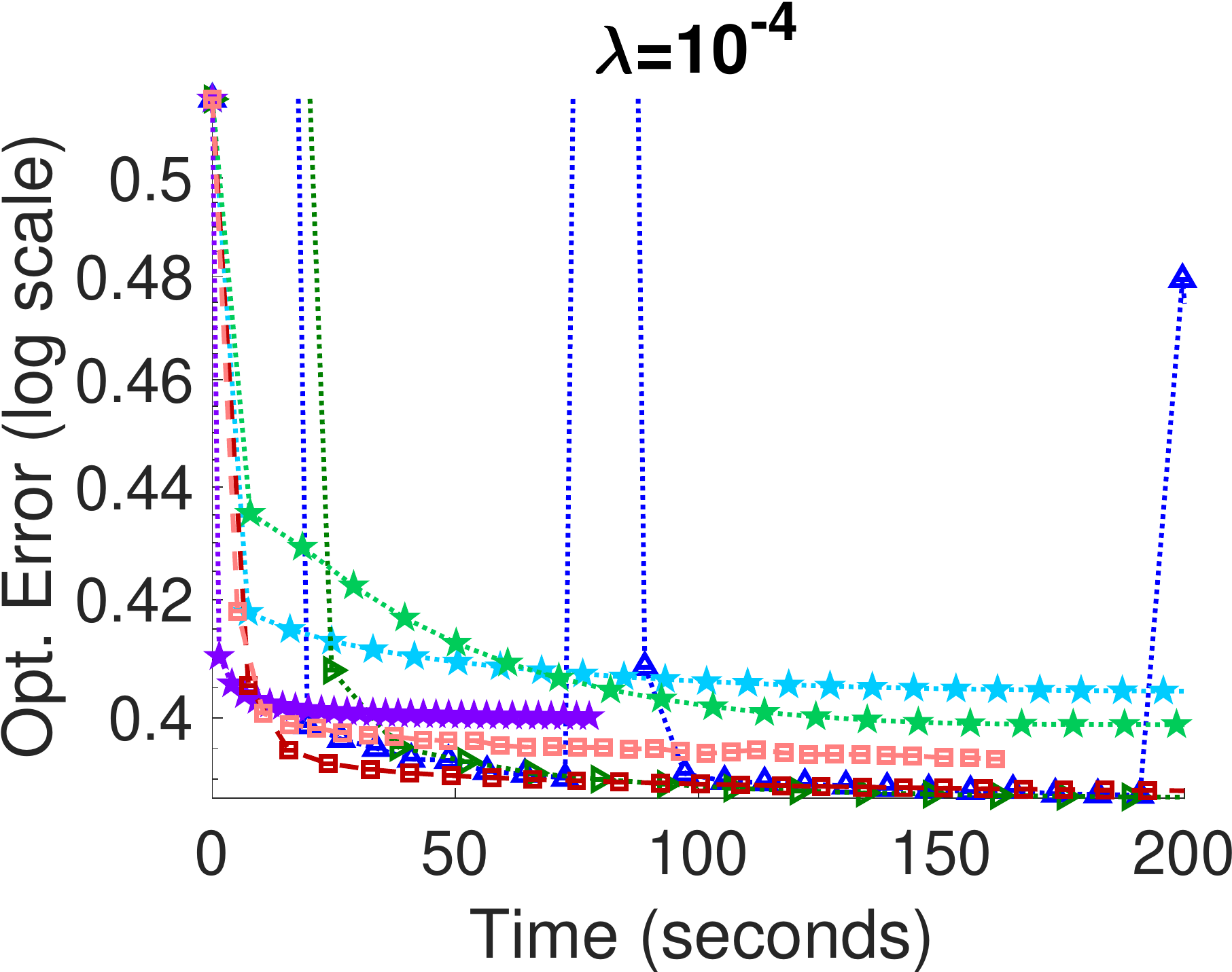} \\
		\includegraphics[width=.32\linewidth,height=4cm]{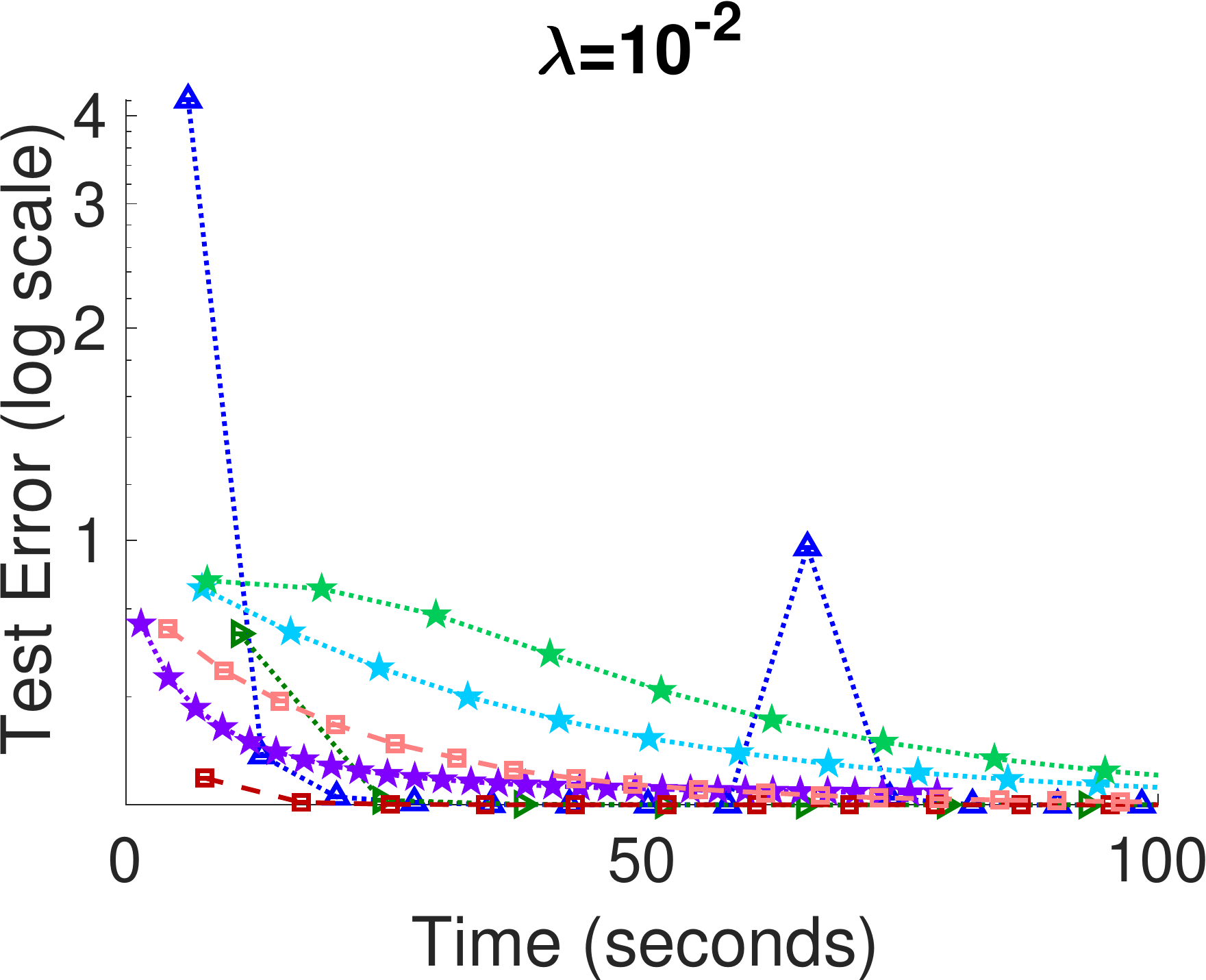} 
		\includegraphics[width=.32\linewidth,height=4cm]{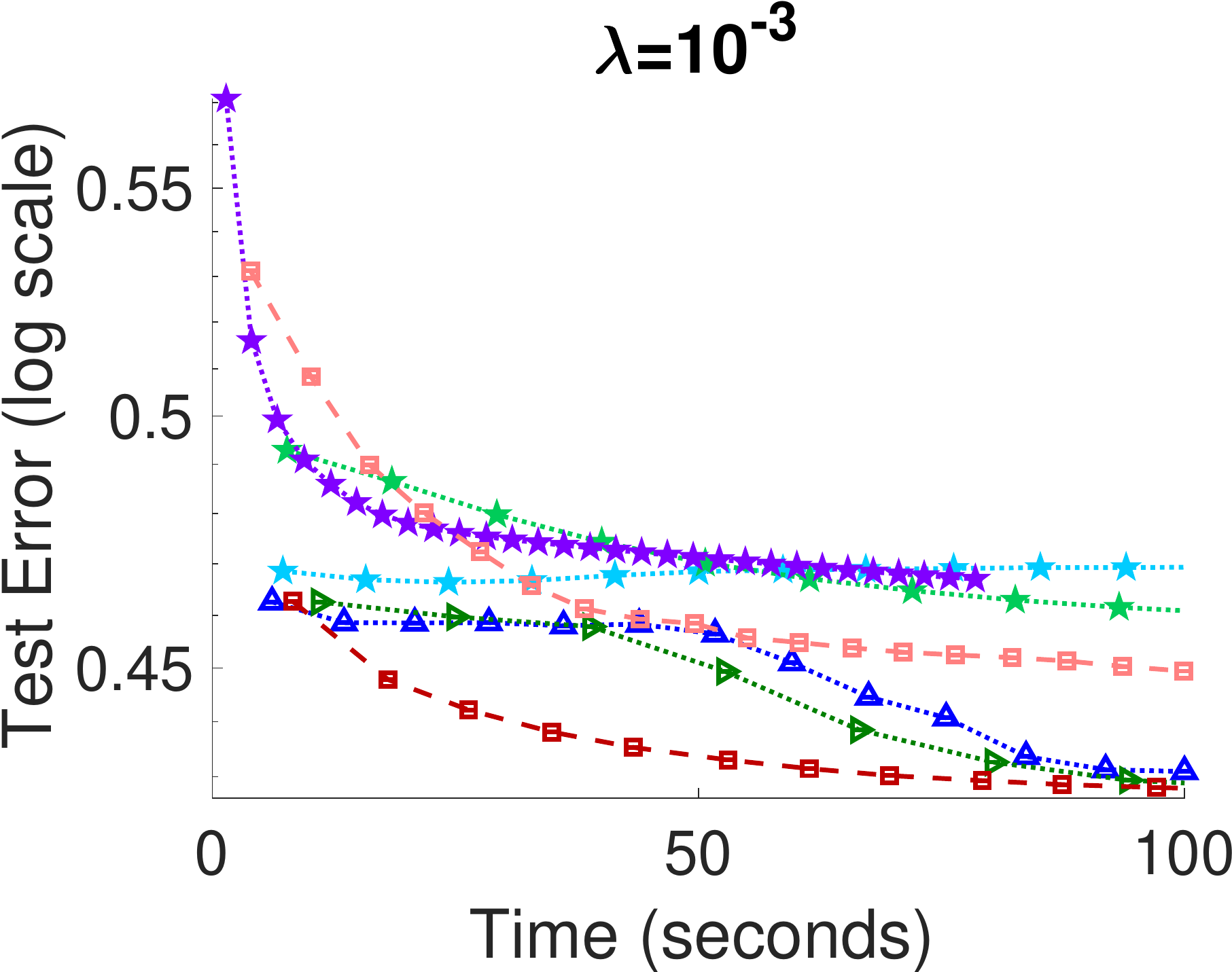} 
		\includegraphics[width=.32\linewidth,height=4cm]{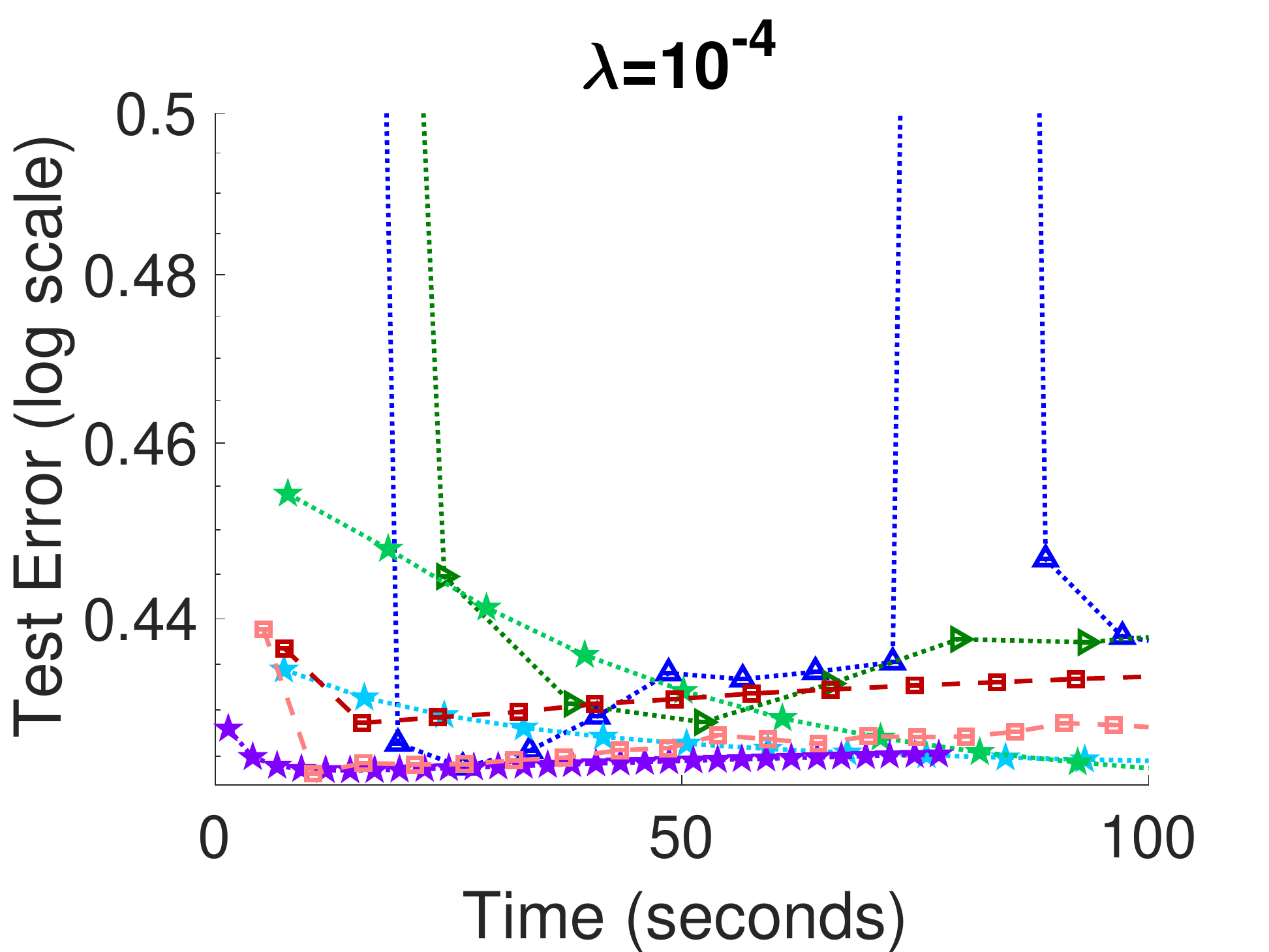} \\
		\includegraphics[width=.32\linewidth,height=4cm]{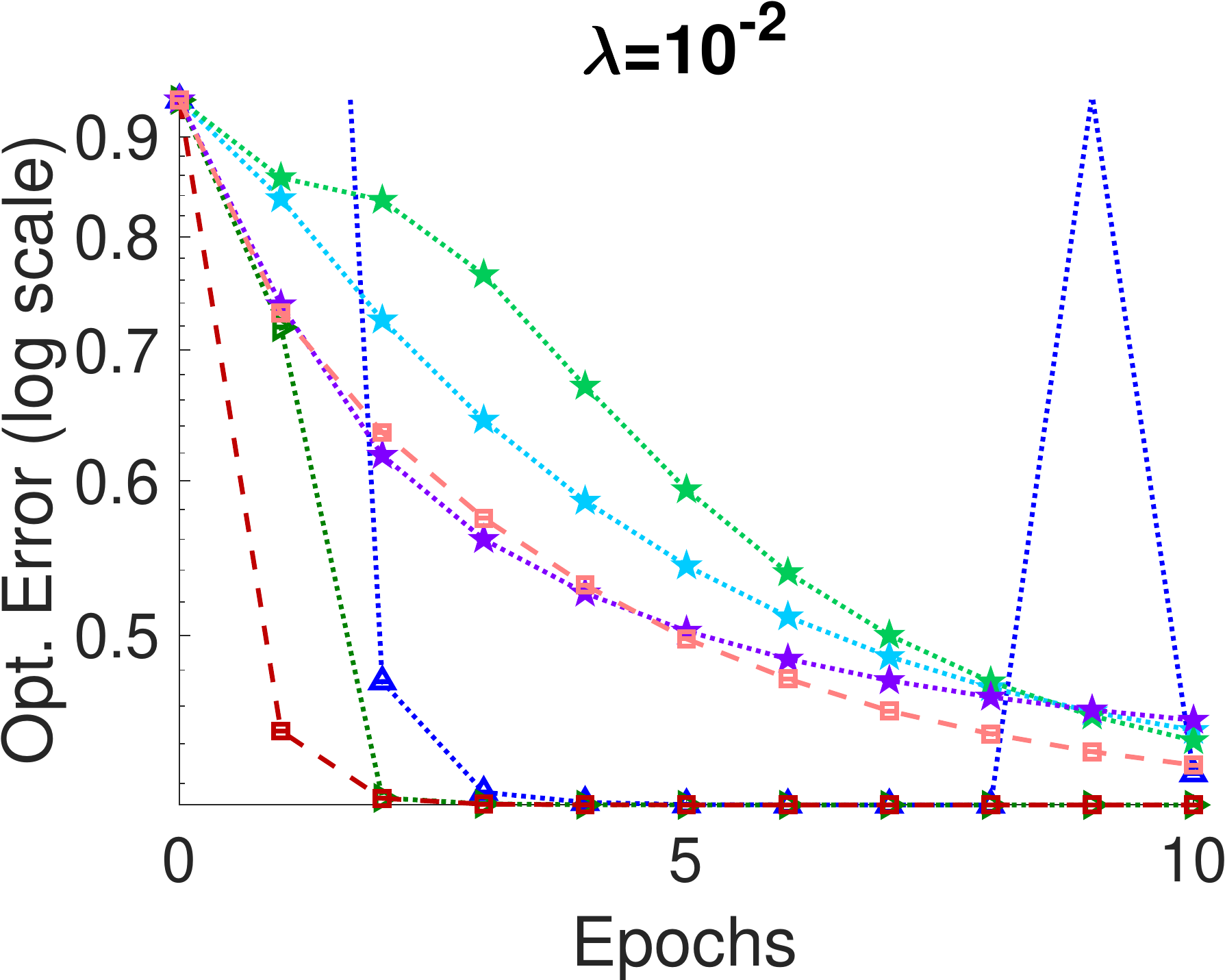} 
		\includegraphics[width=.32\linewidth,height=4cm]{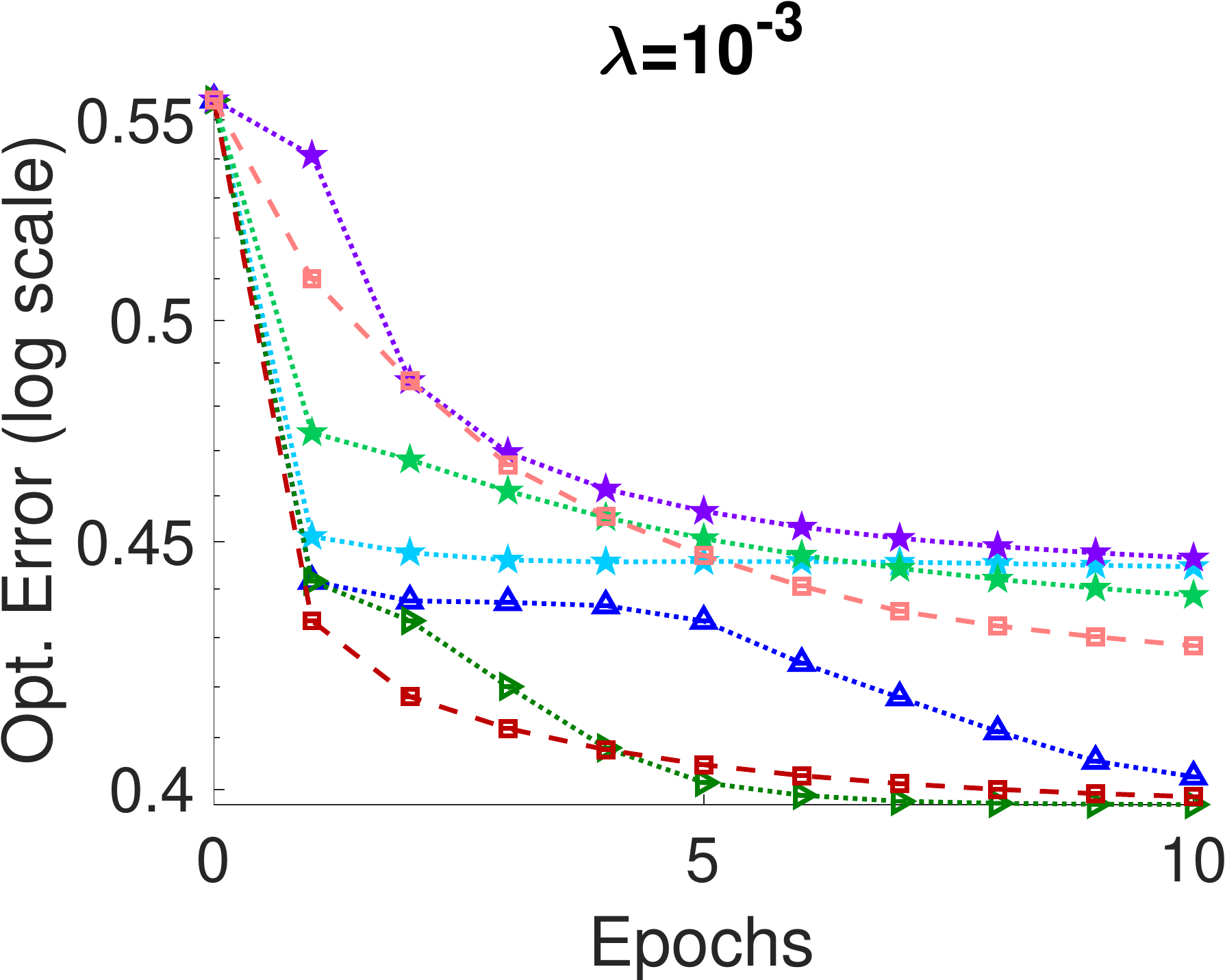} 
		\includegraphics[width=.32\linewidth,height=4cm]{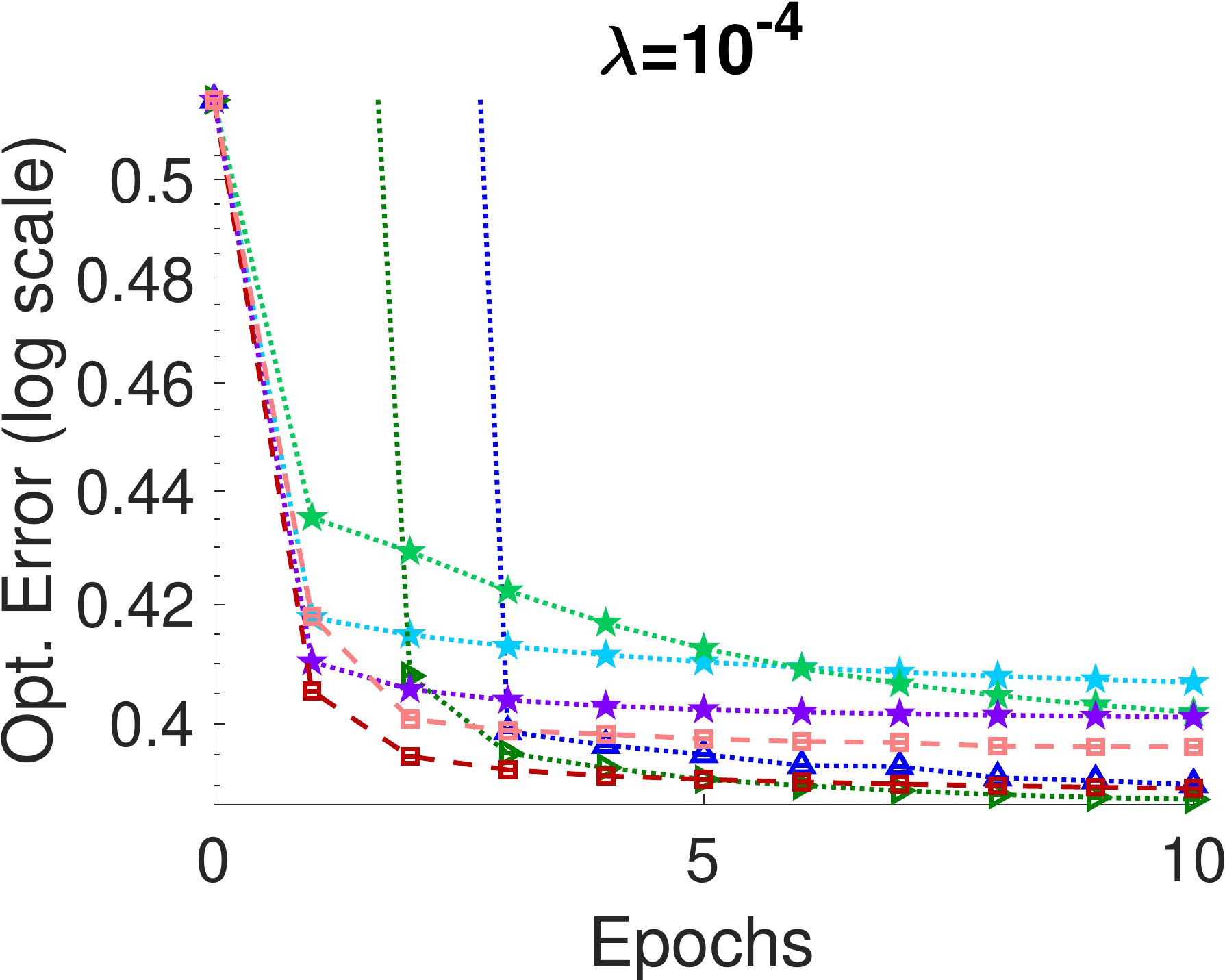} \\
		\includegraphics[width=.32\linewidth,height=4cm]{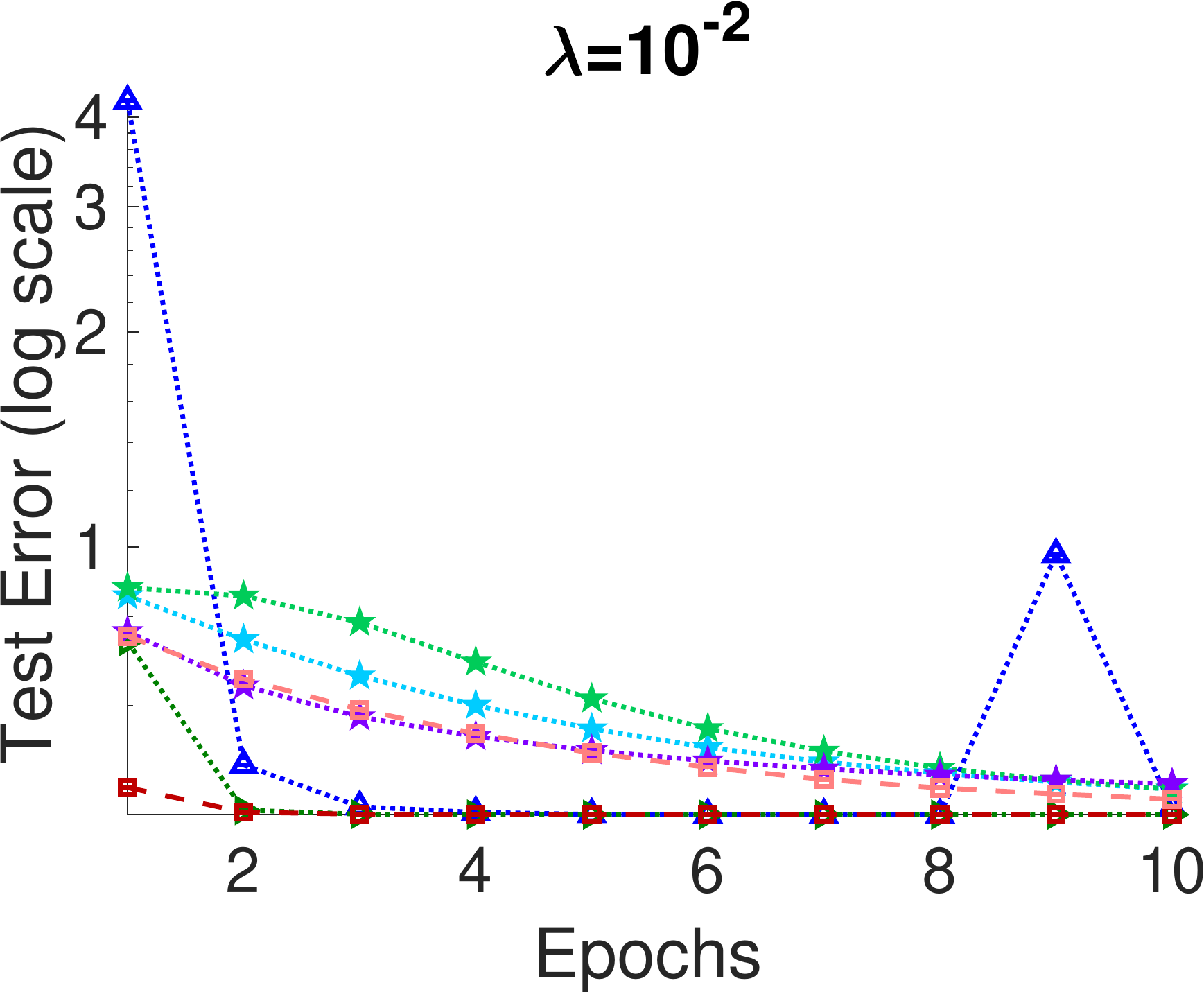} 
		\includegraphics[width=.32\linewidth,height=4cm]{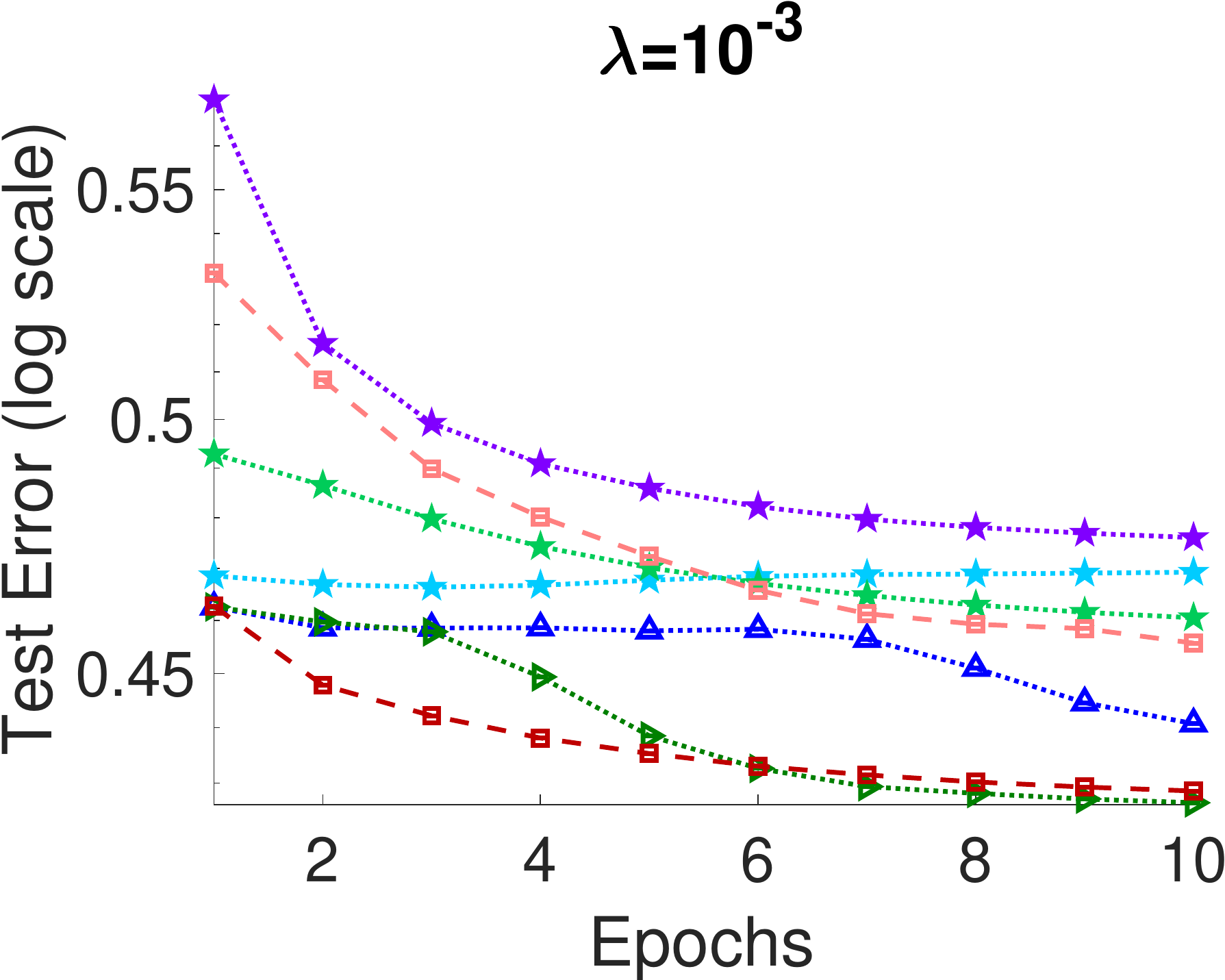} 
		\includegraphics[width=.32\linewidth,height=4cm]{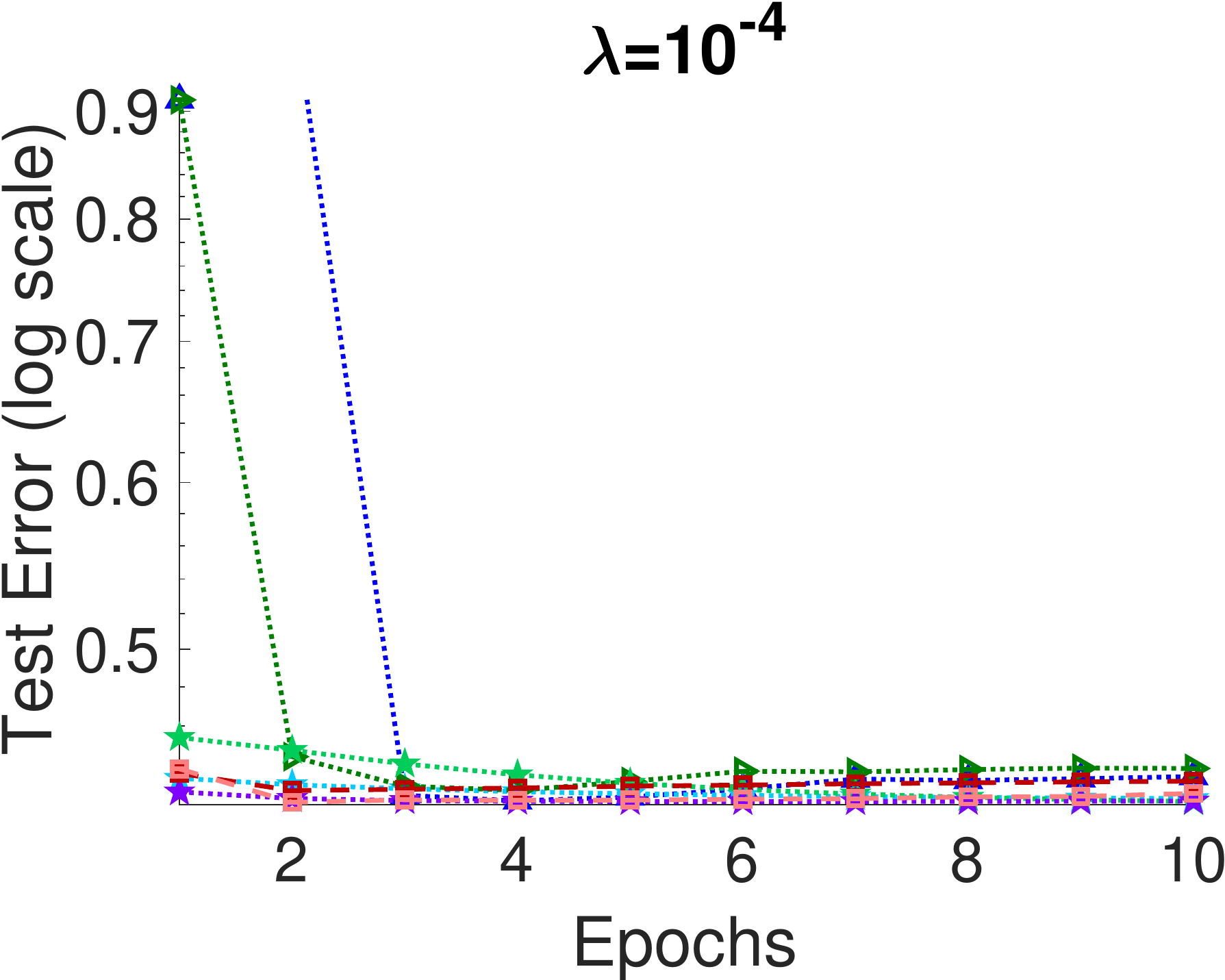}
		\caption{Comparison of proposed methods with the existing optimization methods for optimization performances achieved within 200 seconds on \textit{cifar10} dataset with various regularizer $\lambda$  for \textit{logistic regression}. }
		\label{fig:cifar}
	\end{figure}
	
	\clearpage
	
	\section{Additional numerical experiments on \textit{adult}}
	\begin{figure}[!h]
		\centering
		\includegraphics[width=1\linewidth,keepaspectratio]{ICML_convex_figures/LEGEND_CONVEX.png}\\
		\includegraphics[width=.32\linewidth,height=4cm]{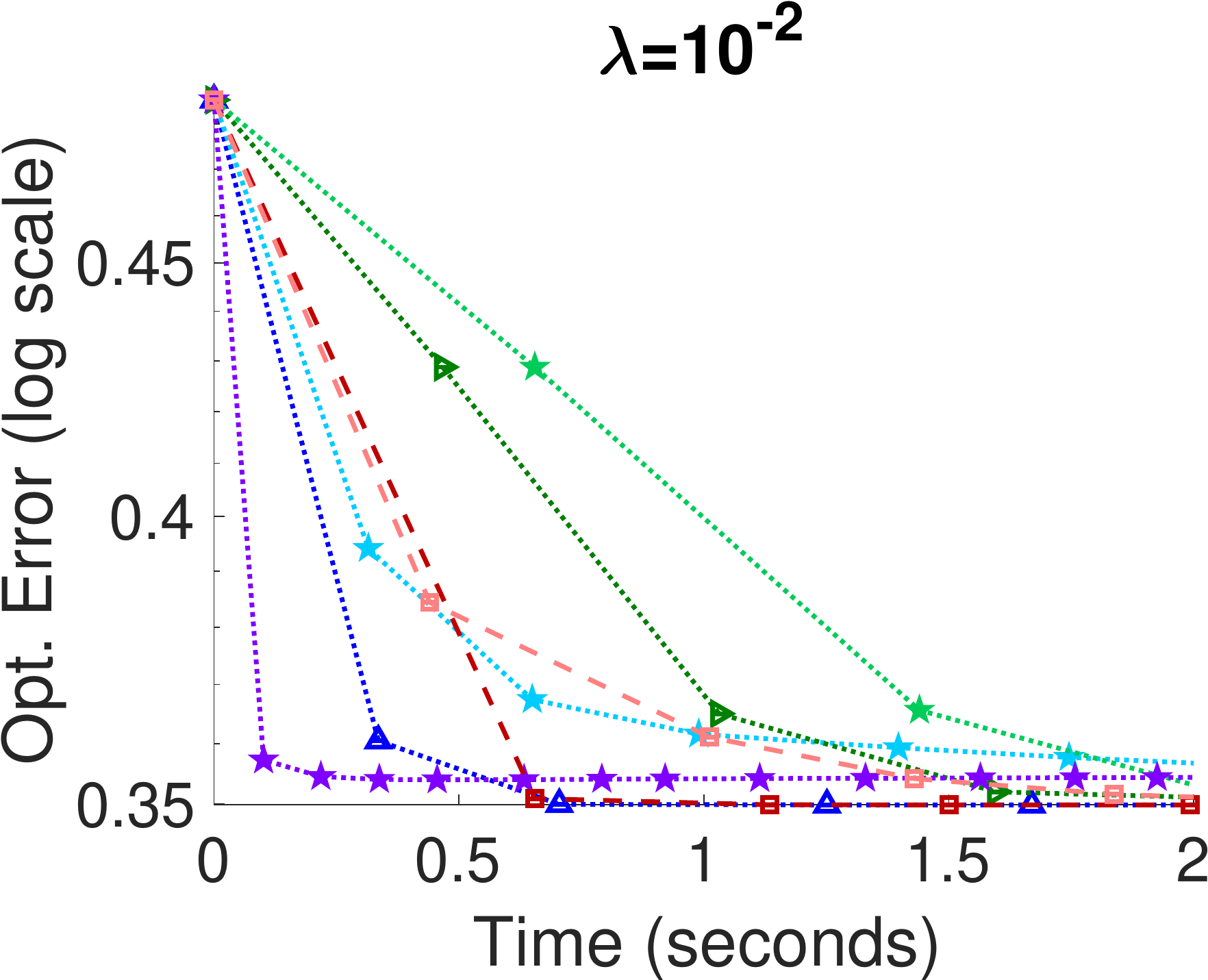}
		\includegraphics[width=.32\linewidth,height=4cm]{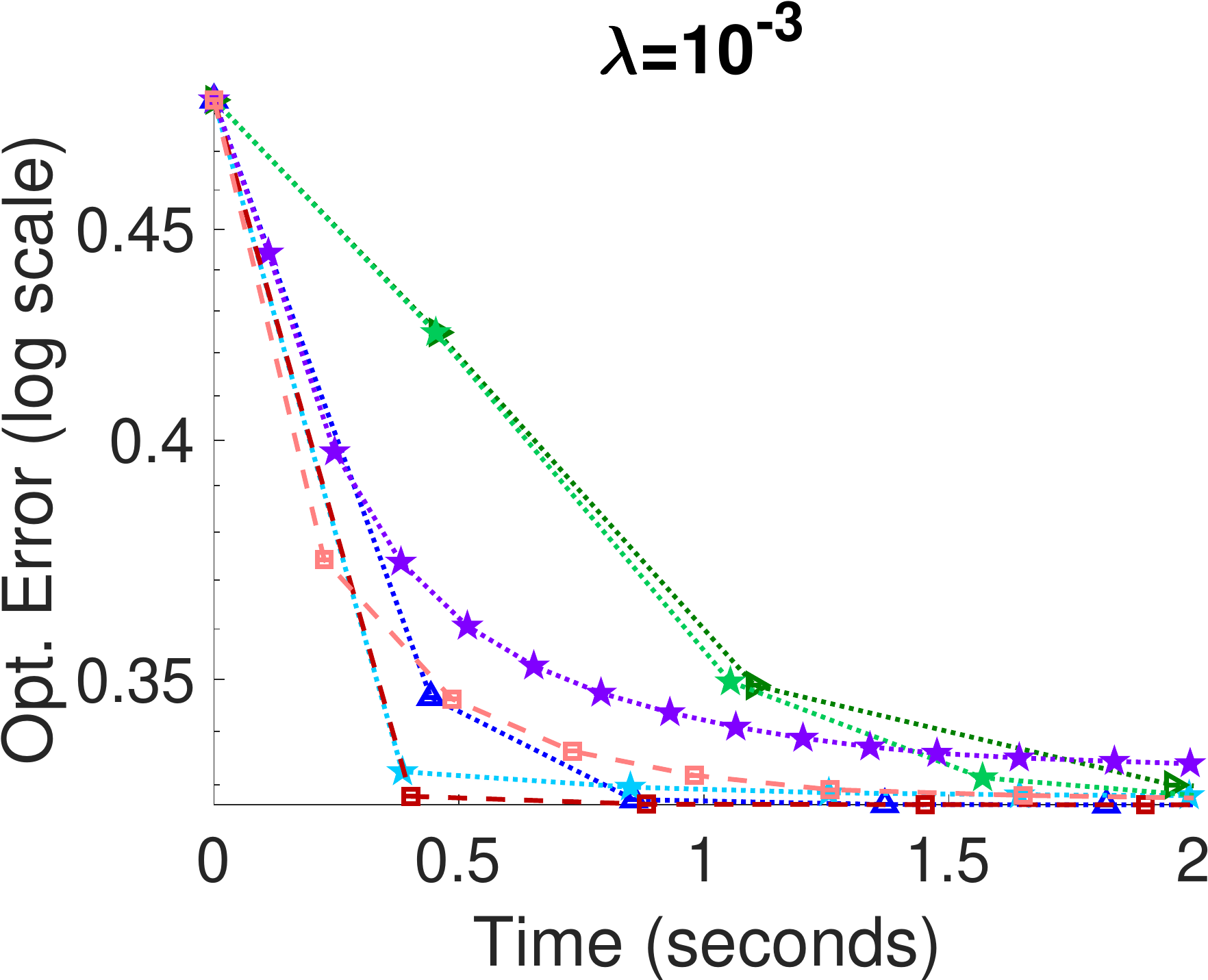}
		\includegraphics[width=.32\linewidth,height=4cm]{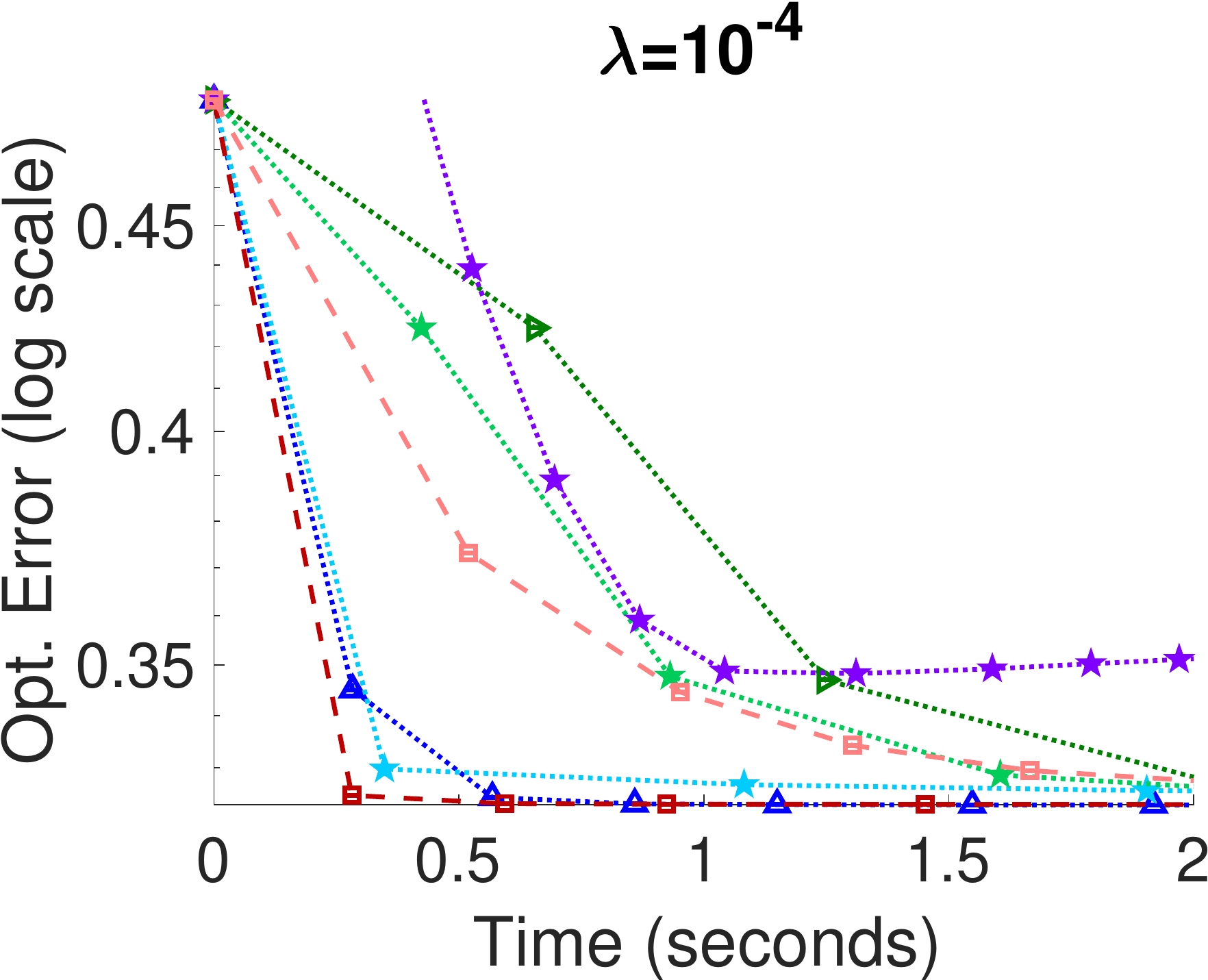} \\
		\includegraphics[width=.32\linewidth,height=4cm]{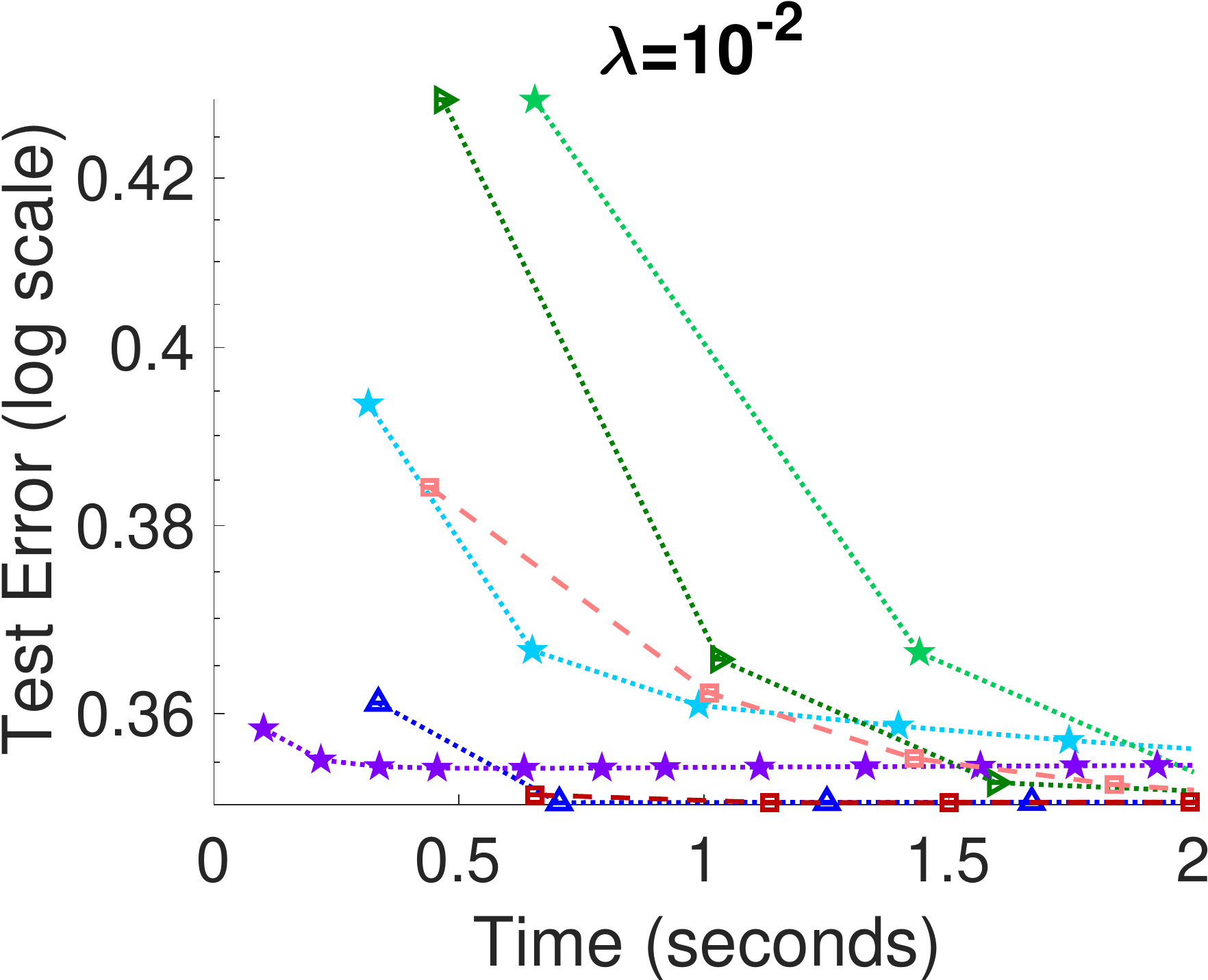}
		\includegraphics[width=.32\linewidth,height=4cm]{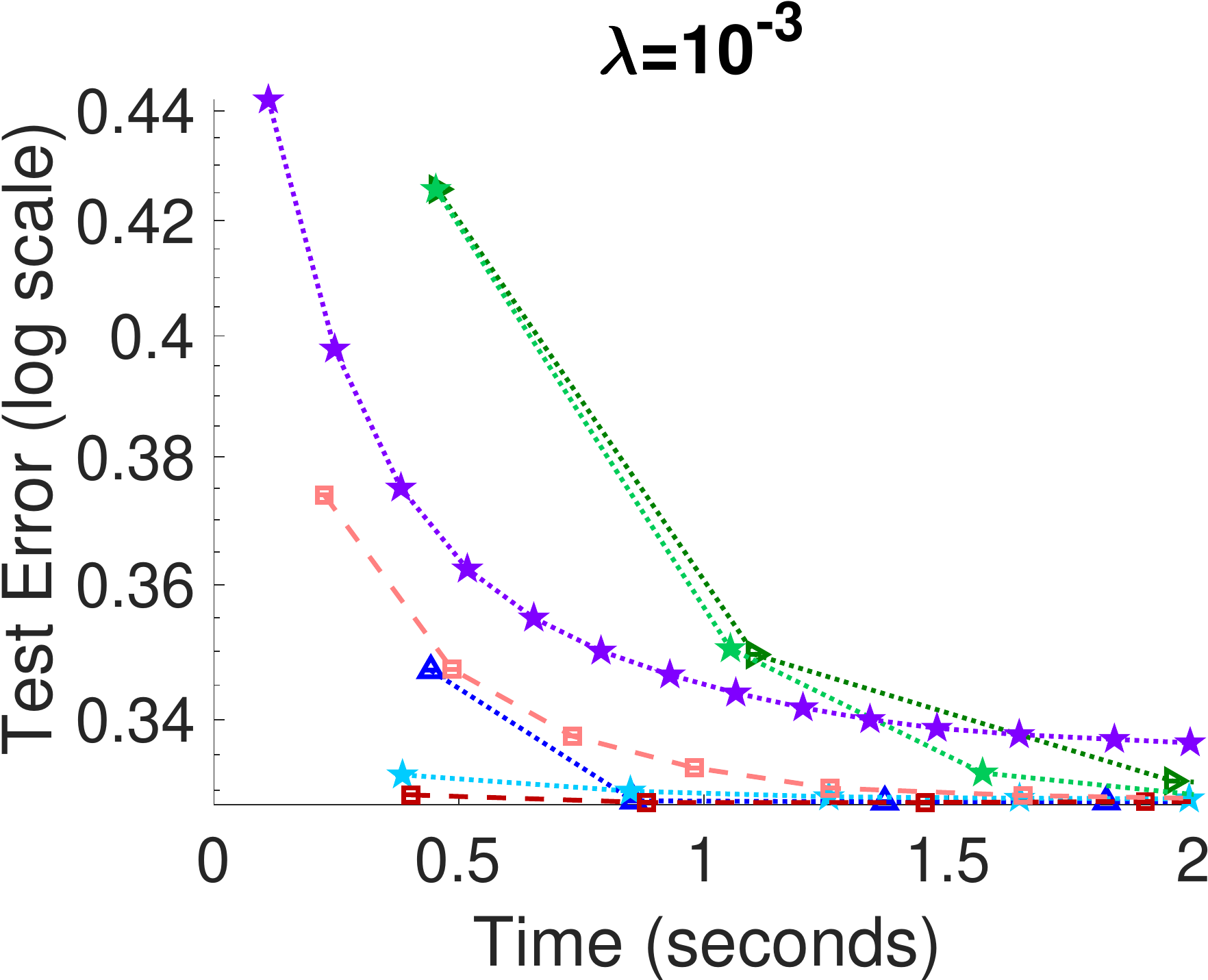}
		\includegraphics[width=.32\linewidth,height=4cm]{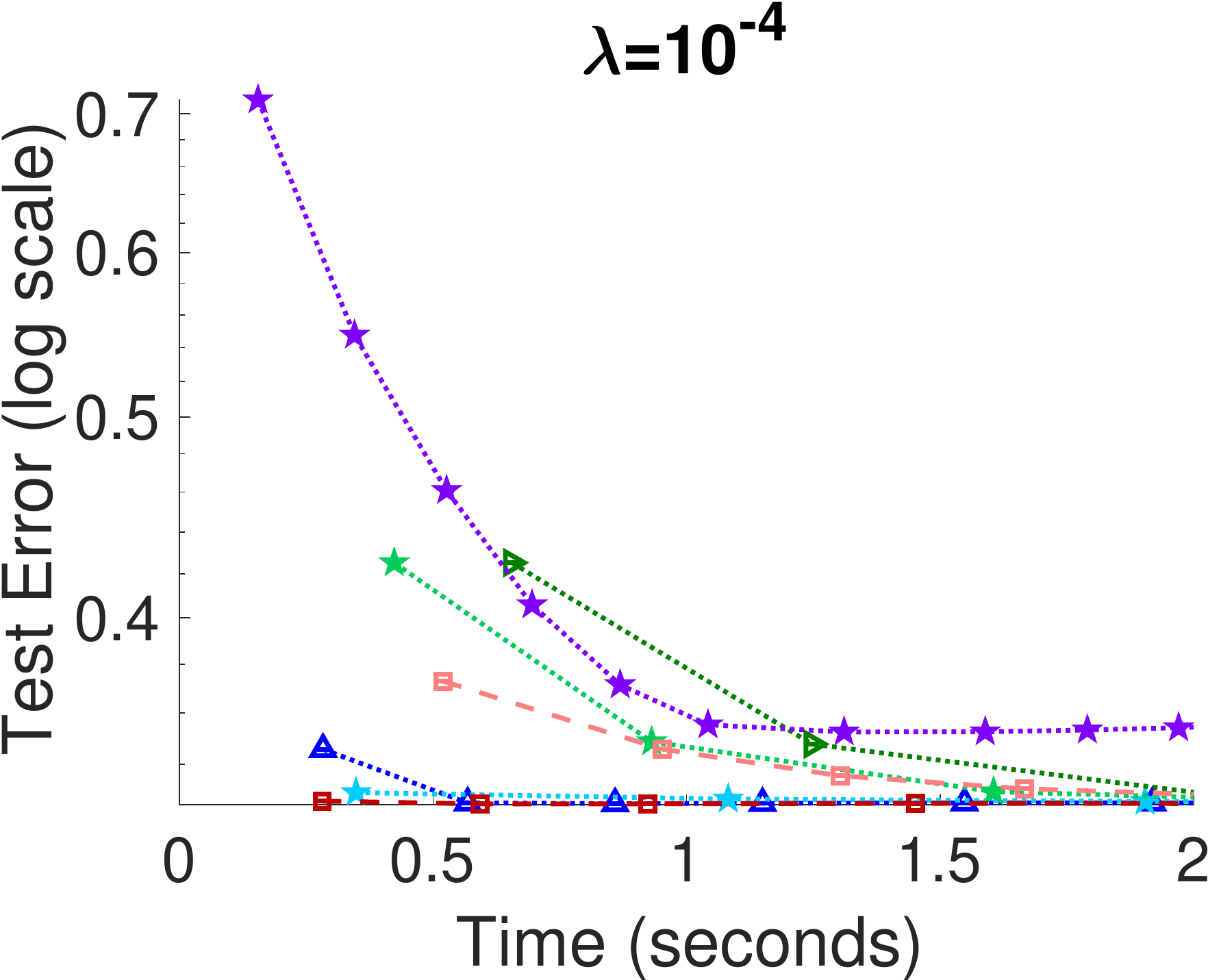} \\
		\includegraphics[width=.32\linewidth,height=4cm]{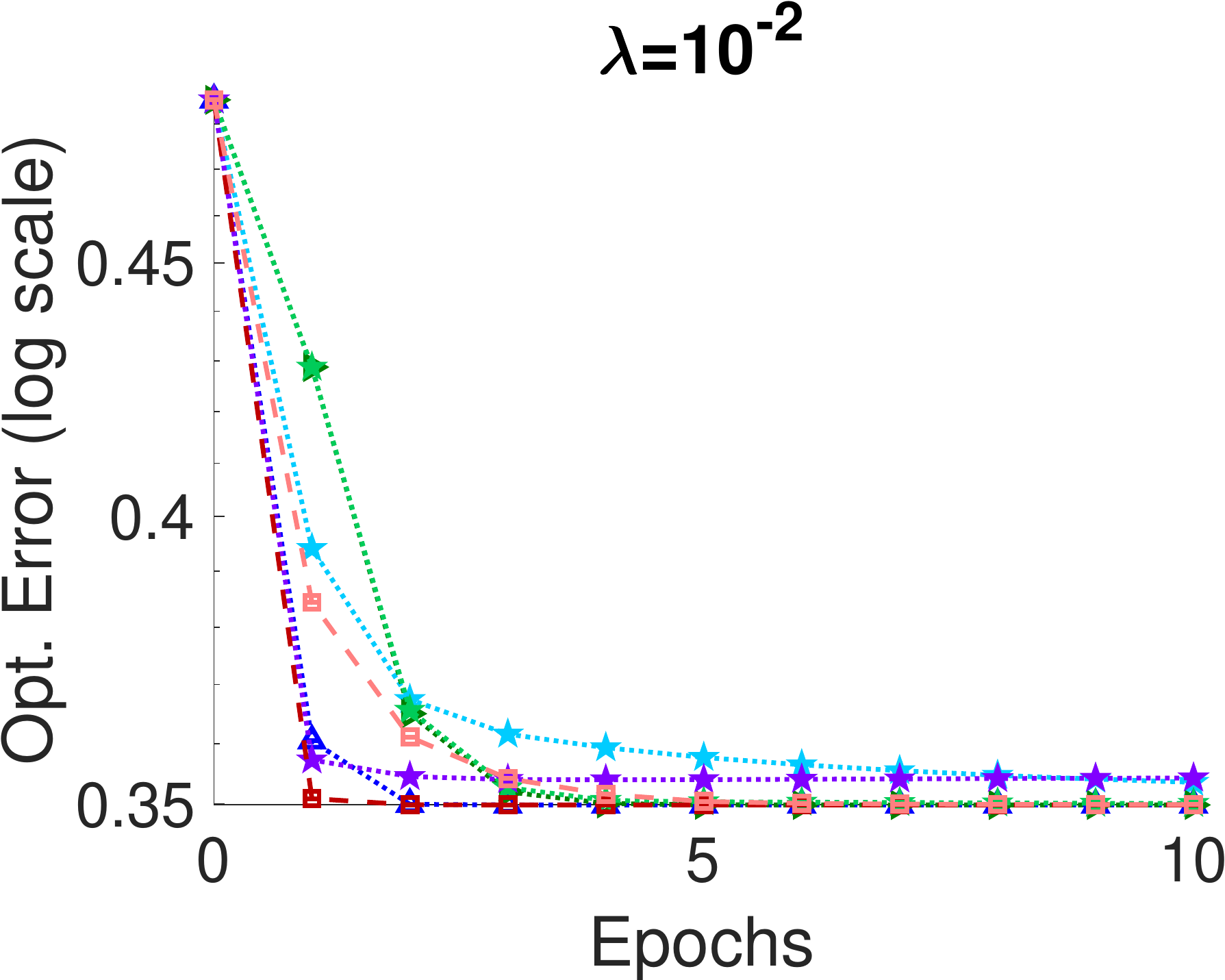}
		\includegraphics[width=.32\linewidth,height=4cm]{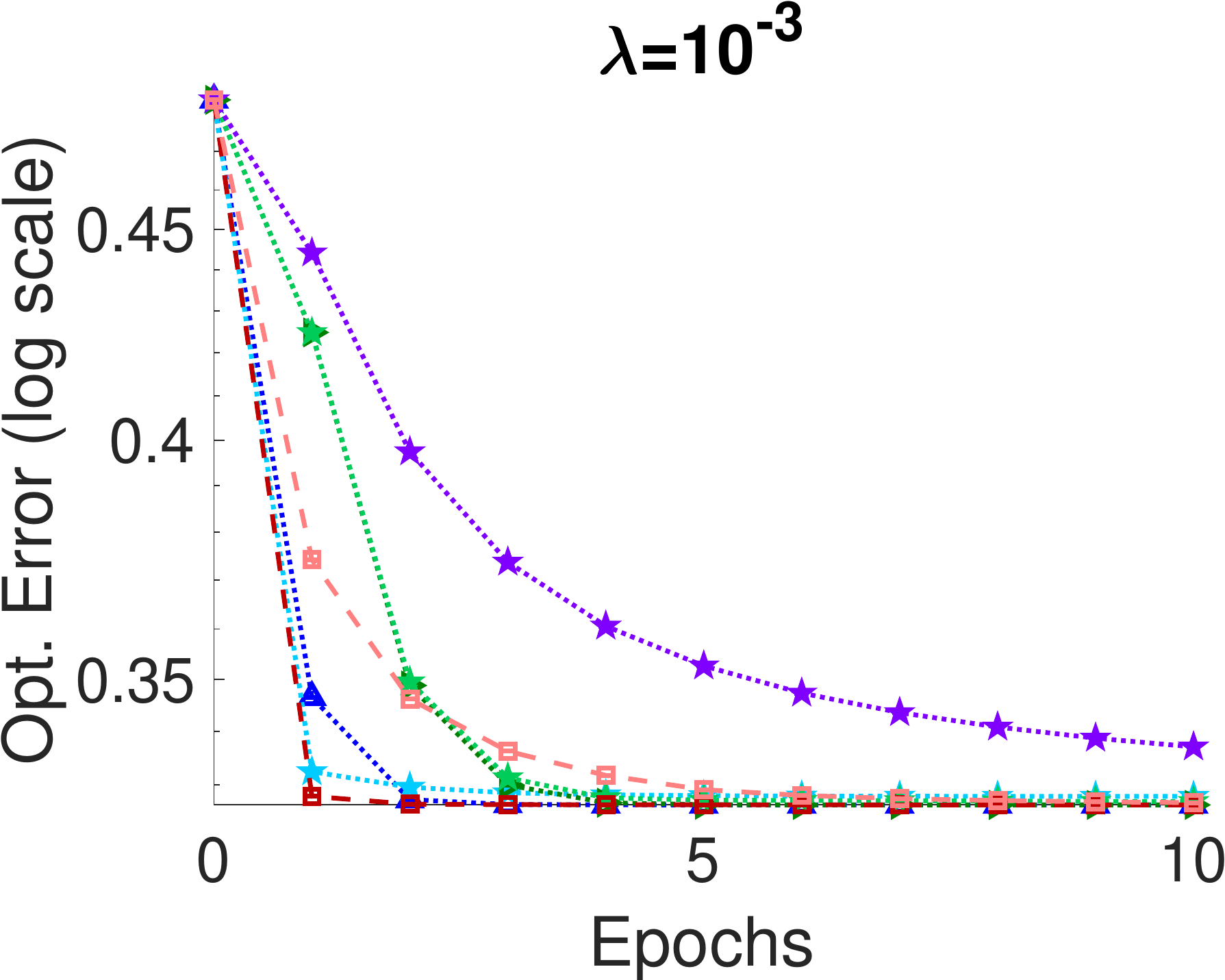}
		\includegraphics[width=.32\linewidth,height=4cm]{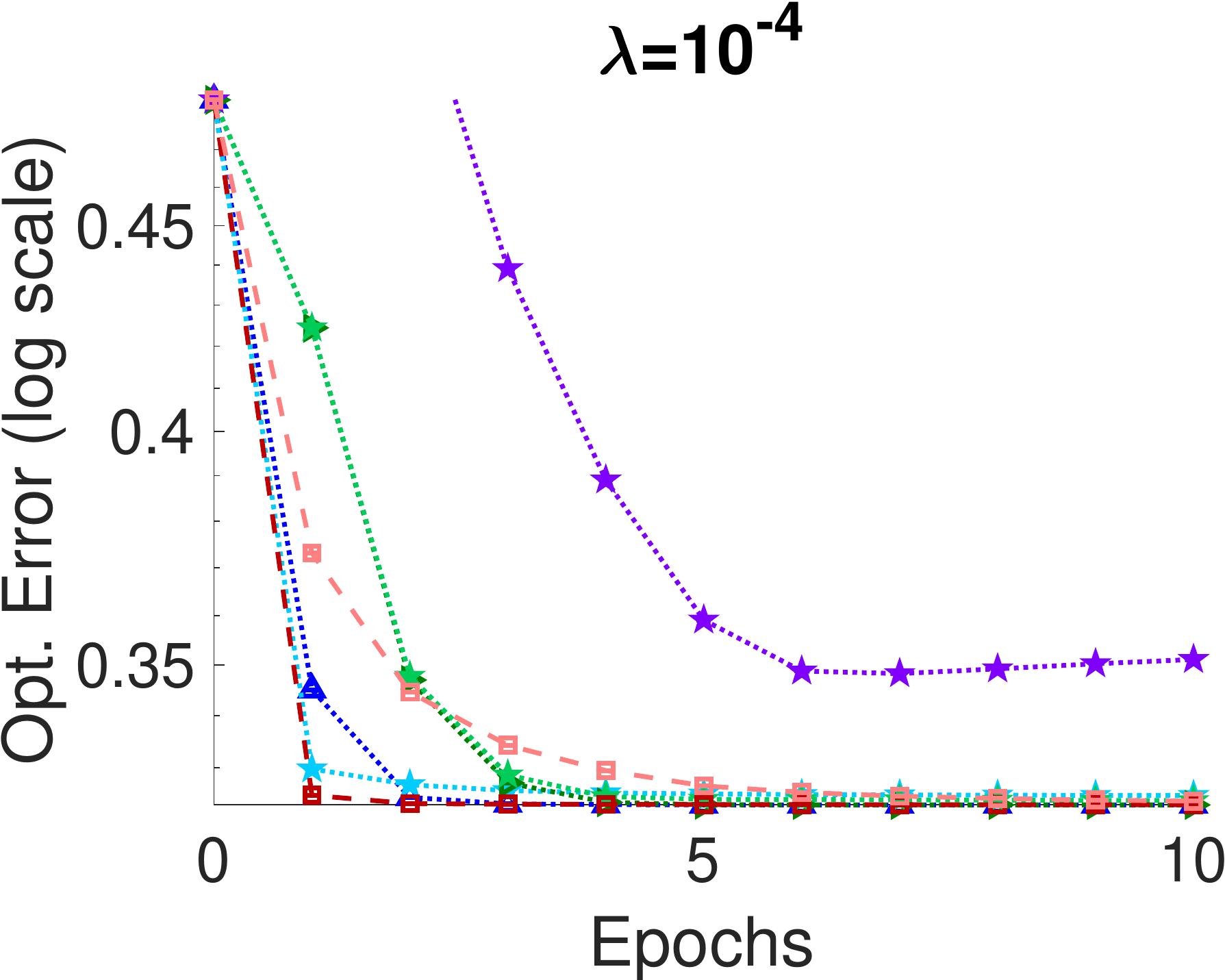} \\
		\includegraphics[width=.32\linewidth,height=4cm]{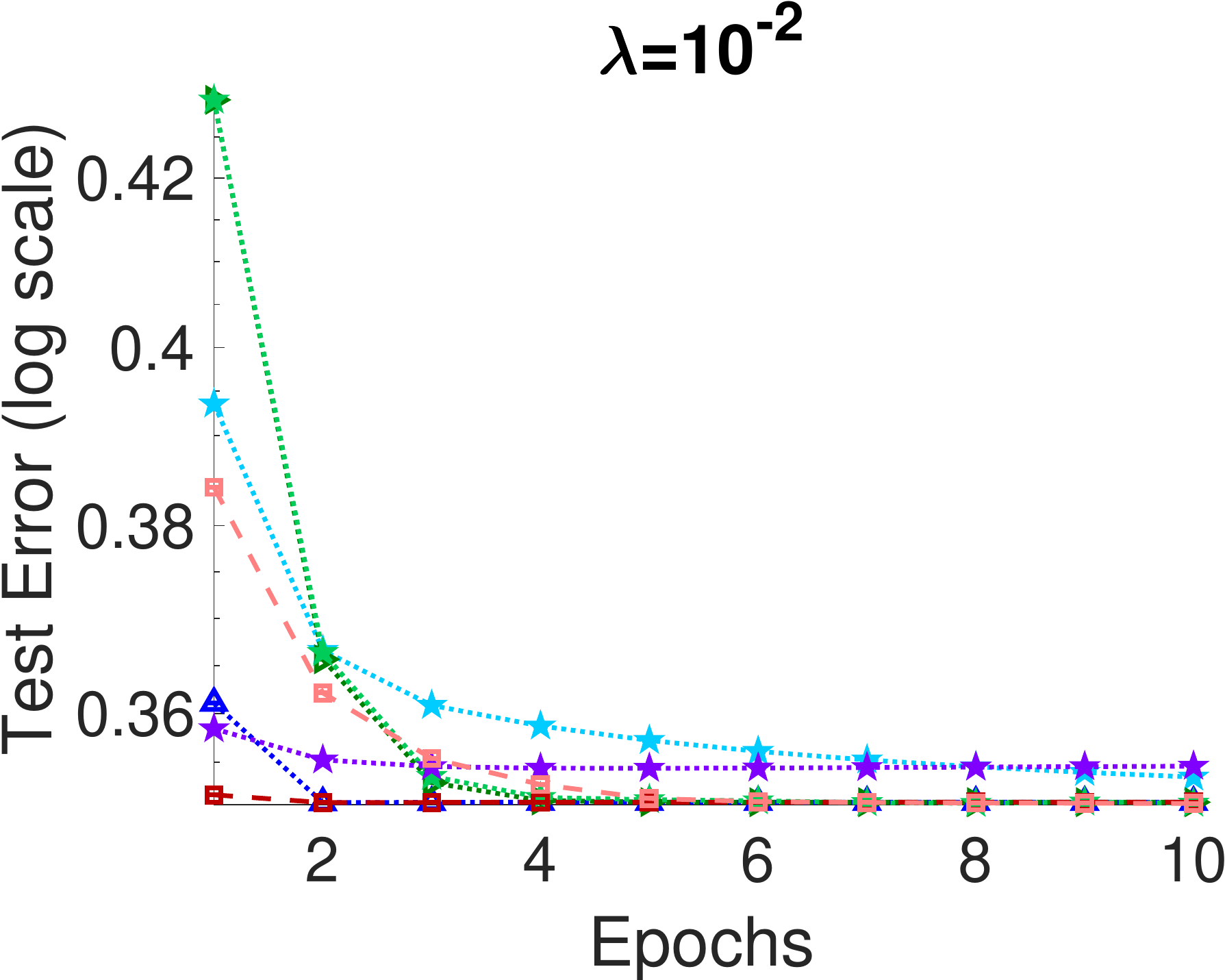}
		\includegraphics[width=.32\linewidth,height=4cm]{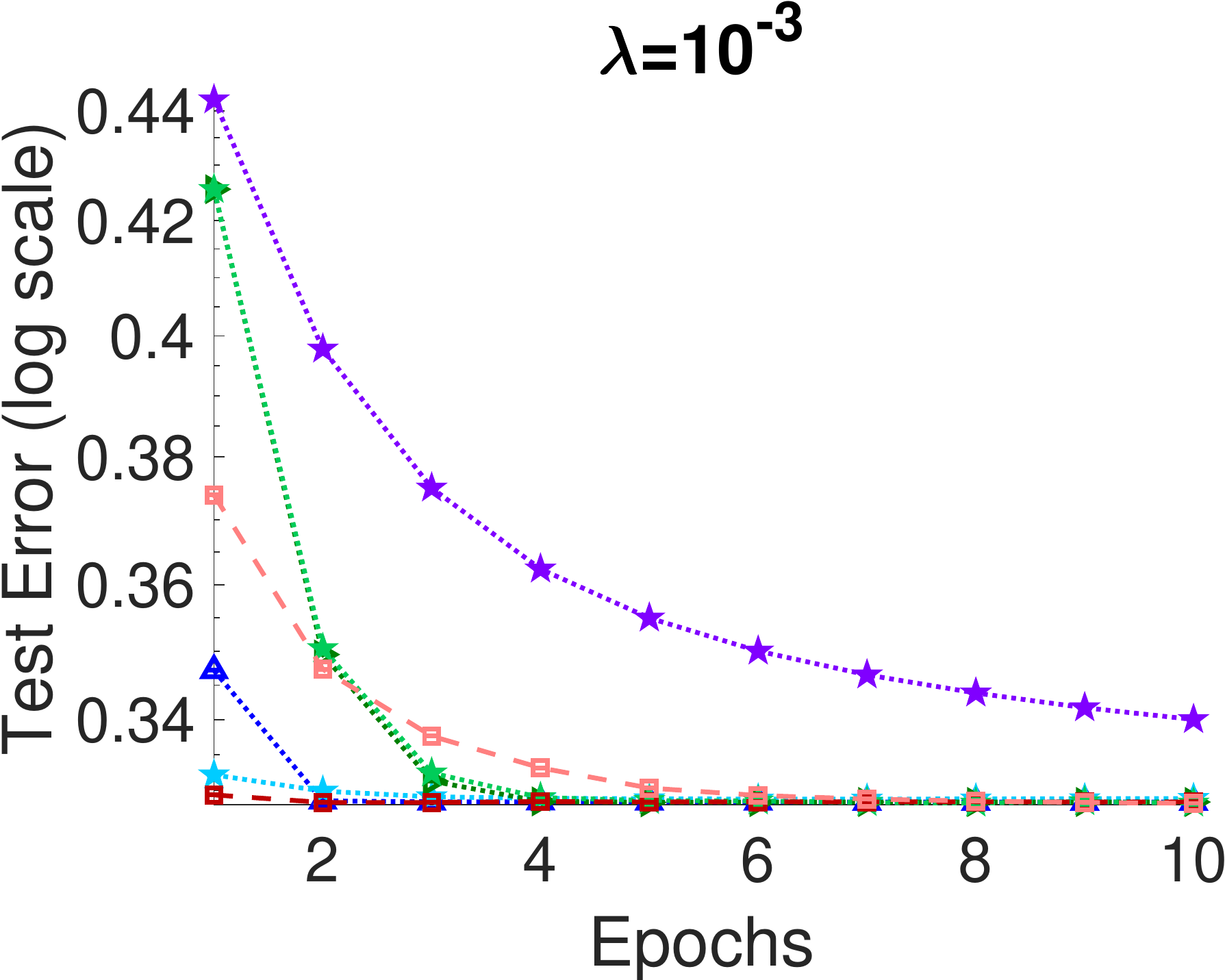}
		\includegraphics[width=.32\linewidth,height=4cm]{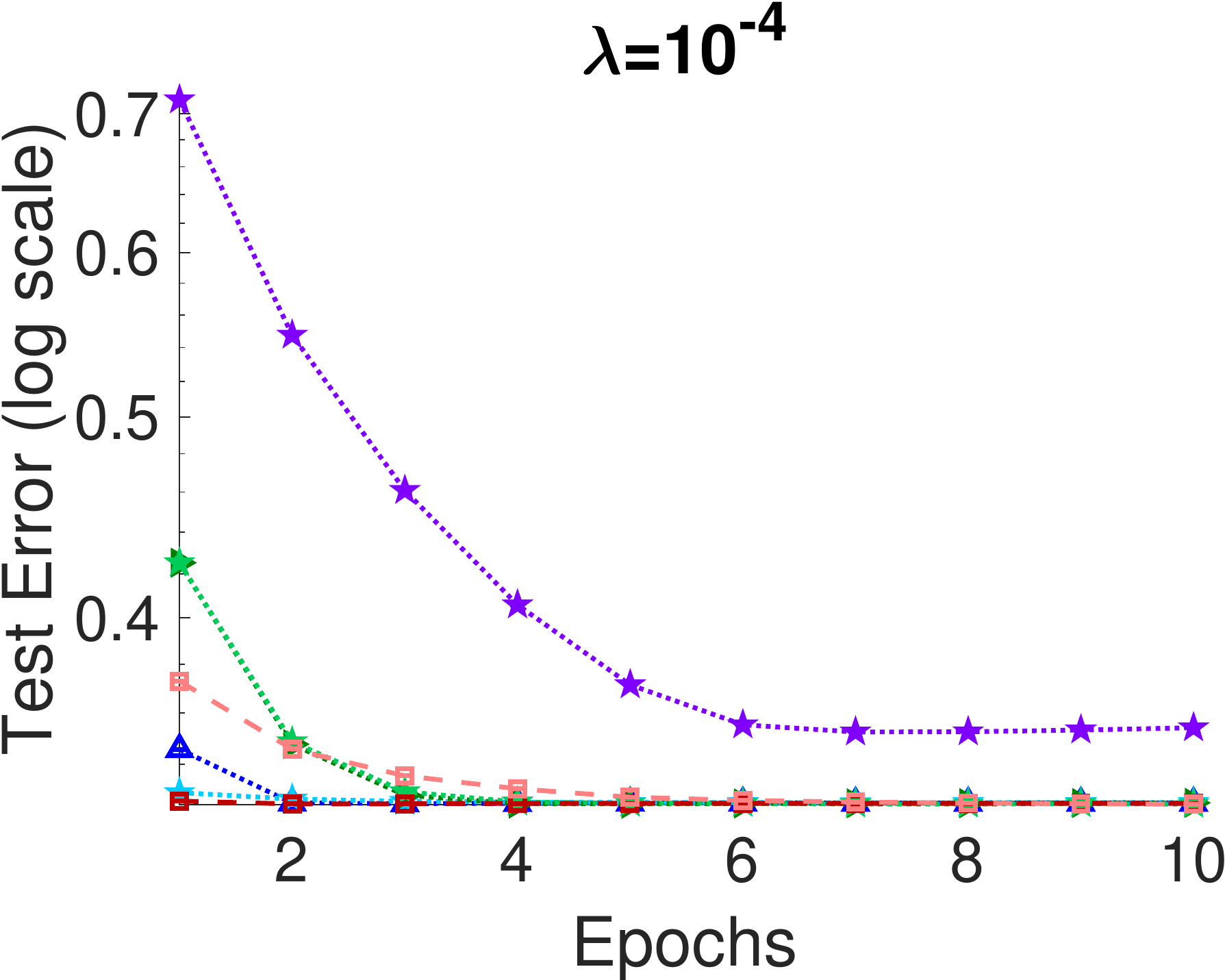}
		\caption{Comparison of proposed methods with the existing optimization methods for optimization performances achieved within 2 seconds on \textit{adult} dataset with various regularizer $\lambda$ for \textit{logistic regression}.}
	\end{figure}

	\clearpage
	
	\section{Comparison with S4QN}
	\begin{figure}[!h]
		\centering
		\begin{tabular}{cc}
			\includegraphics[width=.4\linewidth,keepaspectratio]{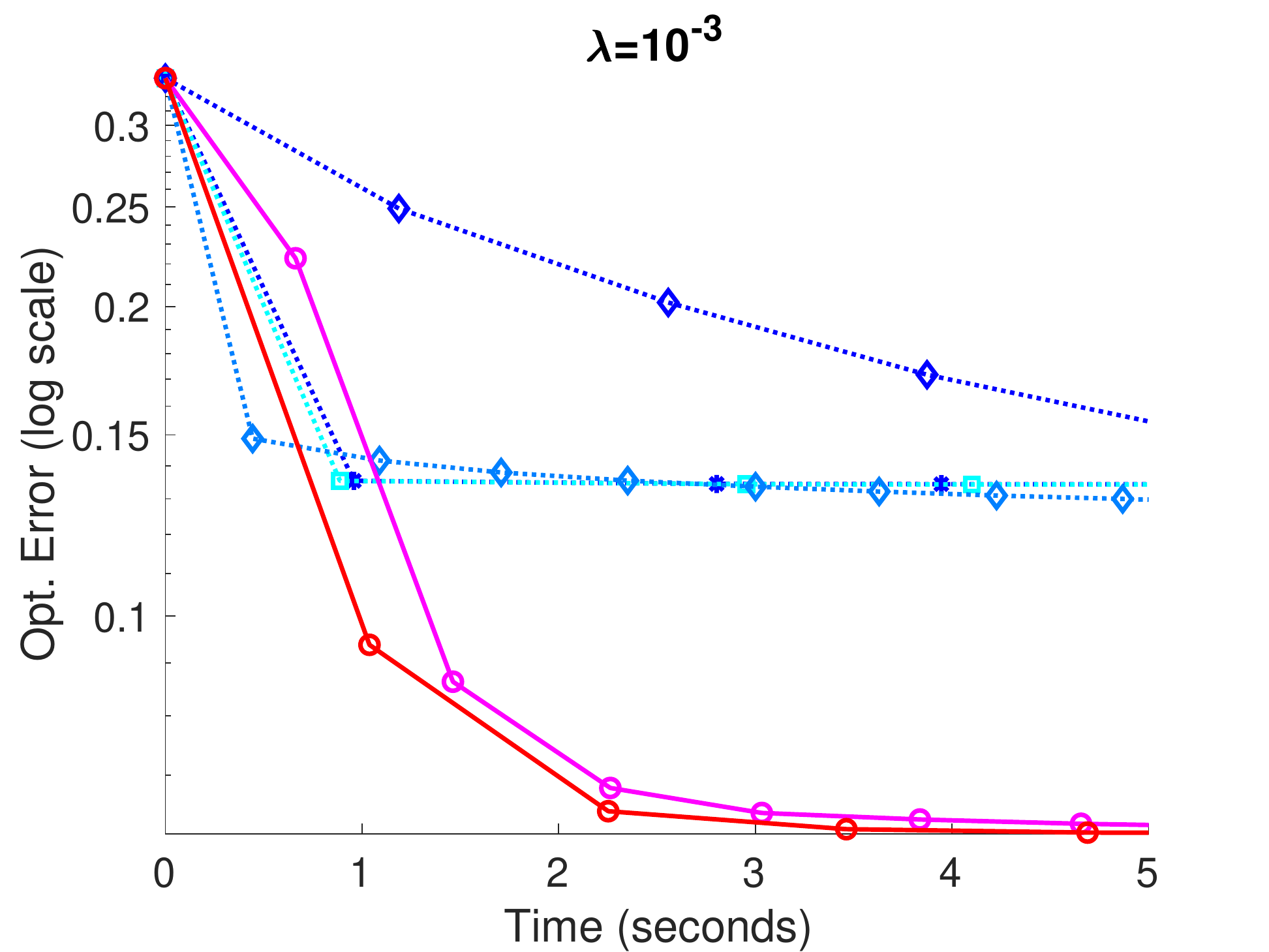}& 
			\includegraphics[width=.4\linewidth,keepaspectratio]{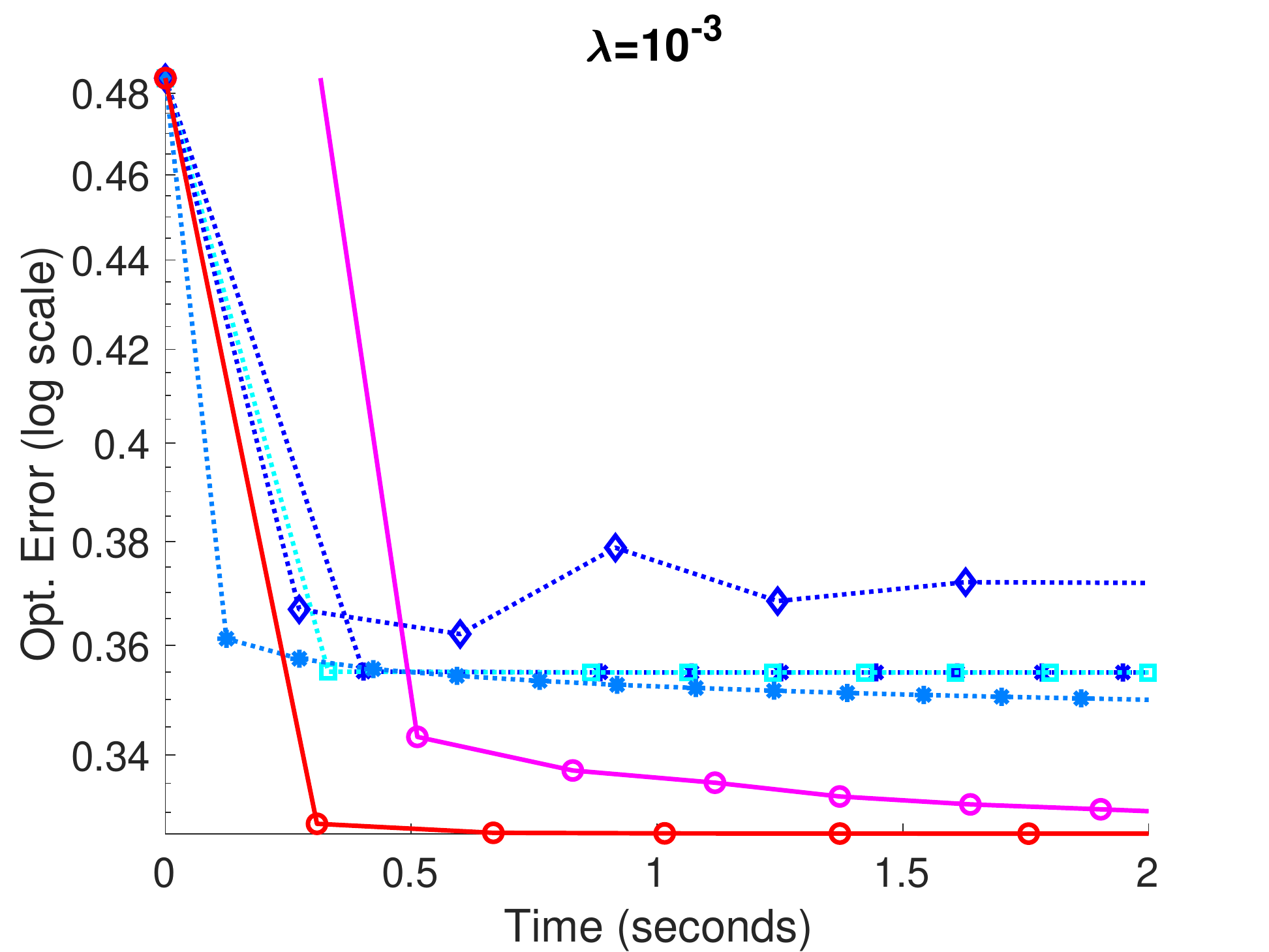}\\ 
			\includegraphics[width=.4\linewidth,keepaspectratio]{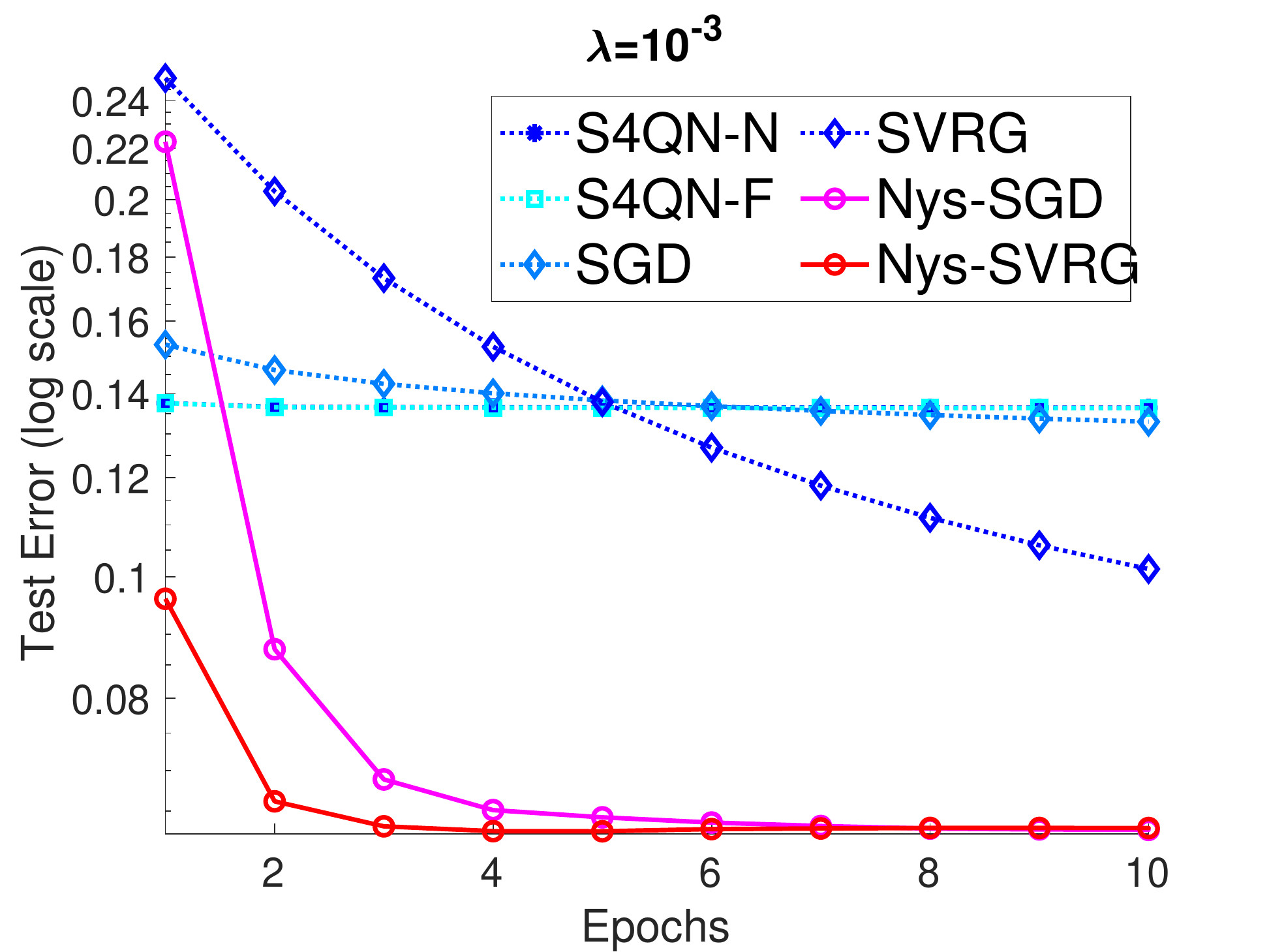}& 
			\includegraphics[width=.4\linewidth,keepaspectratio]{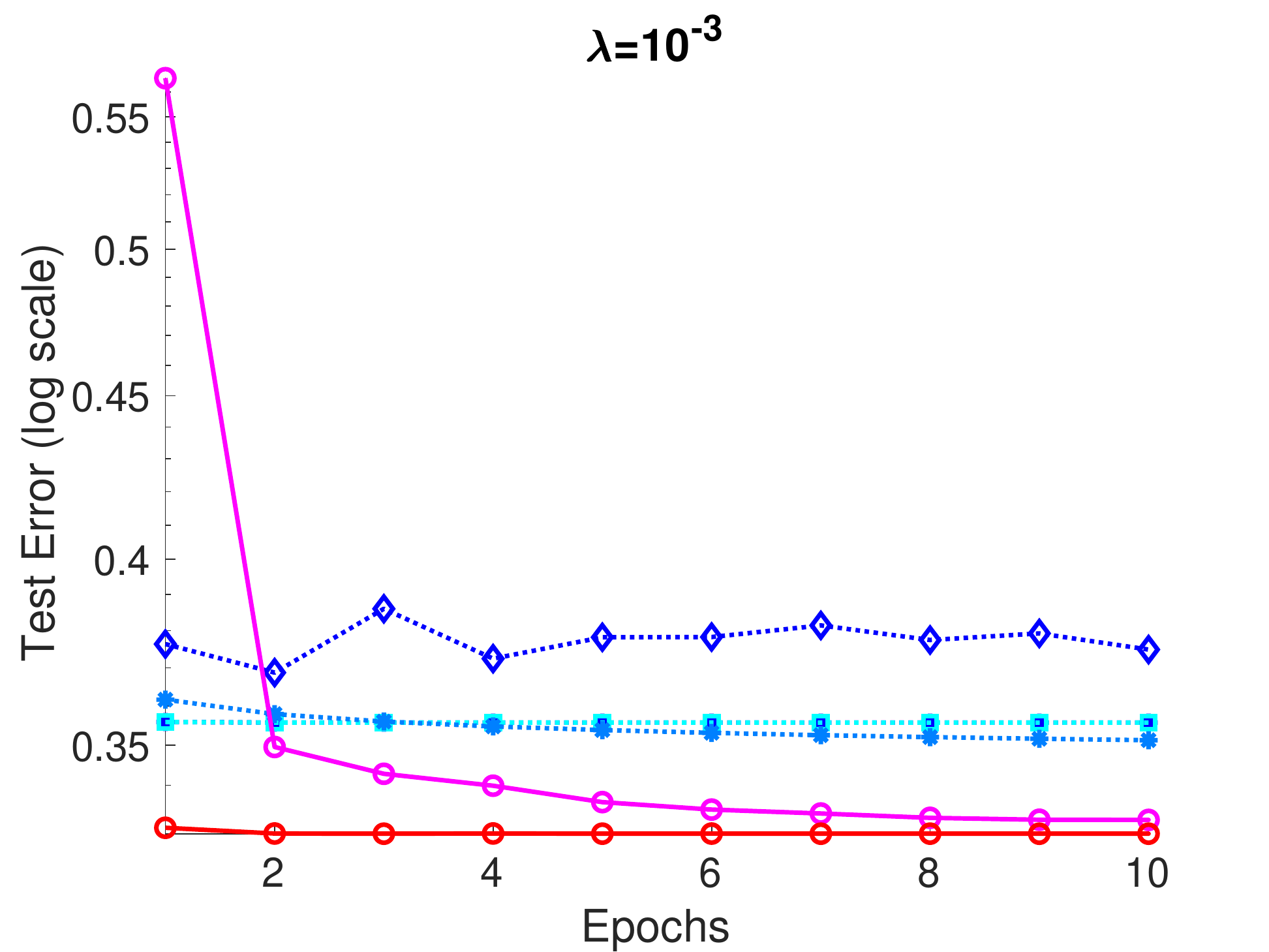}\\
			(A) \textit{w8a} & (B) \textit{adult}\\
		\end{tabular}
		\caption{Comparison of proposed methods with S4QN, SGD, and SVRG on \textit{w8a} and \textit{adult} datasets. These results are computed the results shown in Figure~\ref{fig:s4qn} on an Intel(R) Xeon(R) CPU E5-2690v4 @ 2.60GHz with 14 cores.}
		\label{fig:s4qn}
	\end{figure}
	To show the quality of the approximate curvature information used in the proposed methods, we present a comparison of the proposed methods with first-order methods (SGD and SVRG) used as the base for our algorithms in Figure~\ref{fig:s4qn} on \textit{w8a} and \textit{adult} datasets. The difference between SGD and Nys-SGD on both datasets clearly shows the benefits of using the approximate curvature information and the quality of the approximation.
	methods, Figure~\ref{fig:s4qn} also presents a performance comparison with the recently proposed state-of-the-art S4QN method~\cite{yang2020structured}. For S4QN, we used two variants with sub-sampled Newton (S4QN-N) and Fisher information matrix (S4QN-F) as the base matrix and set $r_1=0.1$ and $r_2=10$. 
	The time taken in the initial epoch by S4QN is significantly higher, as it computes Hessian on a sub-sample at each update. However, the use of a growing gradient batch size reduces the frequency of updates per epoch to speed up subsequent epochs, but decreases the quality of curvature information. Hence, no significant improvement may be observed in the optimization error. However, the proposed Nys-SGD and Nys-SVRG outperformed both variants of the S4QN, indicating that the proposed methods use a better approximation of the curvature information.
	\clearpage
	
	\section{Results on the $\ell_2$-svm loss function}
	\begin{figure}[!h]
		\centering
		\includegraphics[width=1\linewidth]{ICML_convex_figures/LEGEND_CONVEX.png}\\
		\includegraphics[width=.48\linewidth,height=5.5cm]{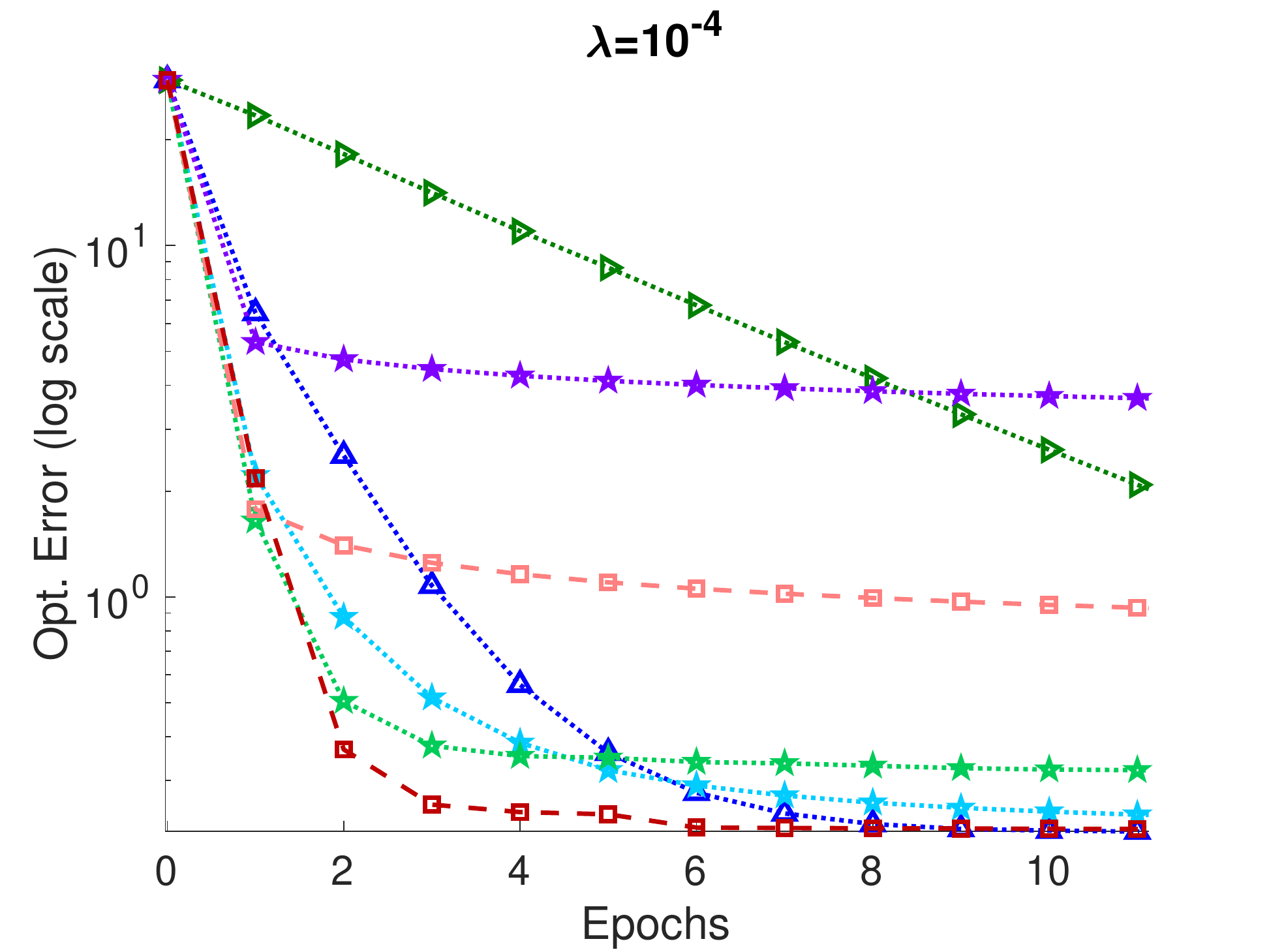}
		\includegraphics[width=.48\linewidth,height=5.5cm]{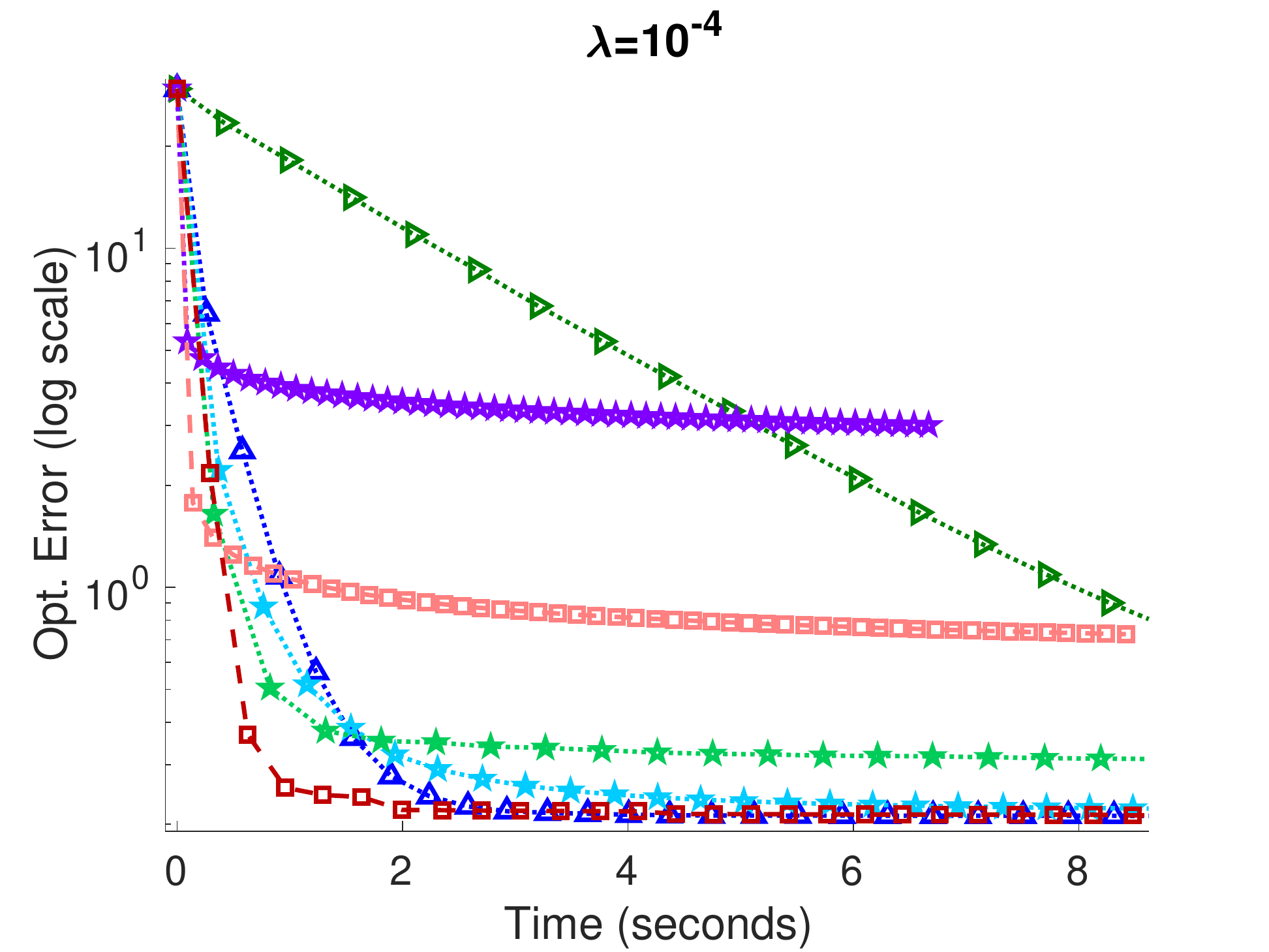}\\
		\caption{Comparision of results for $\ell_2$-svm loss function on \textit{adult} dataset.}
		\label{fig:my_labelSVMadult}
	\end{figure}
	\begin{figure}[!h]
		\centering
		\includegraphics[width=1\linewidth]{ICML_convex_figures/LEGEND_CONVEX.png}\\
		\includegraphics[width=.48\linewidth,height=5.5cm]{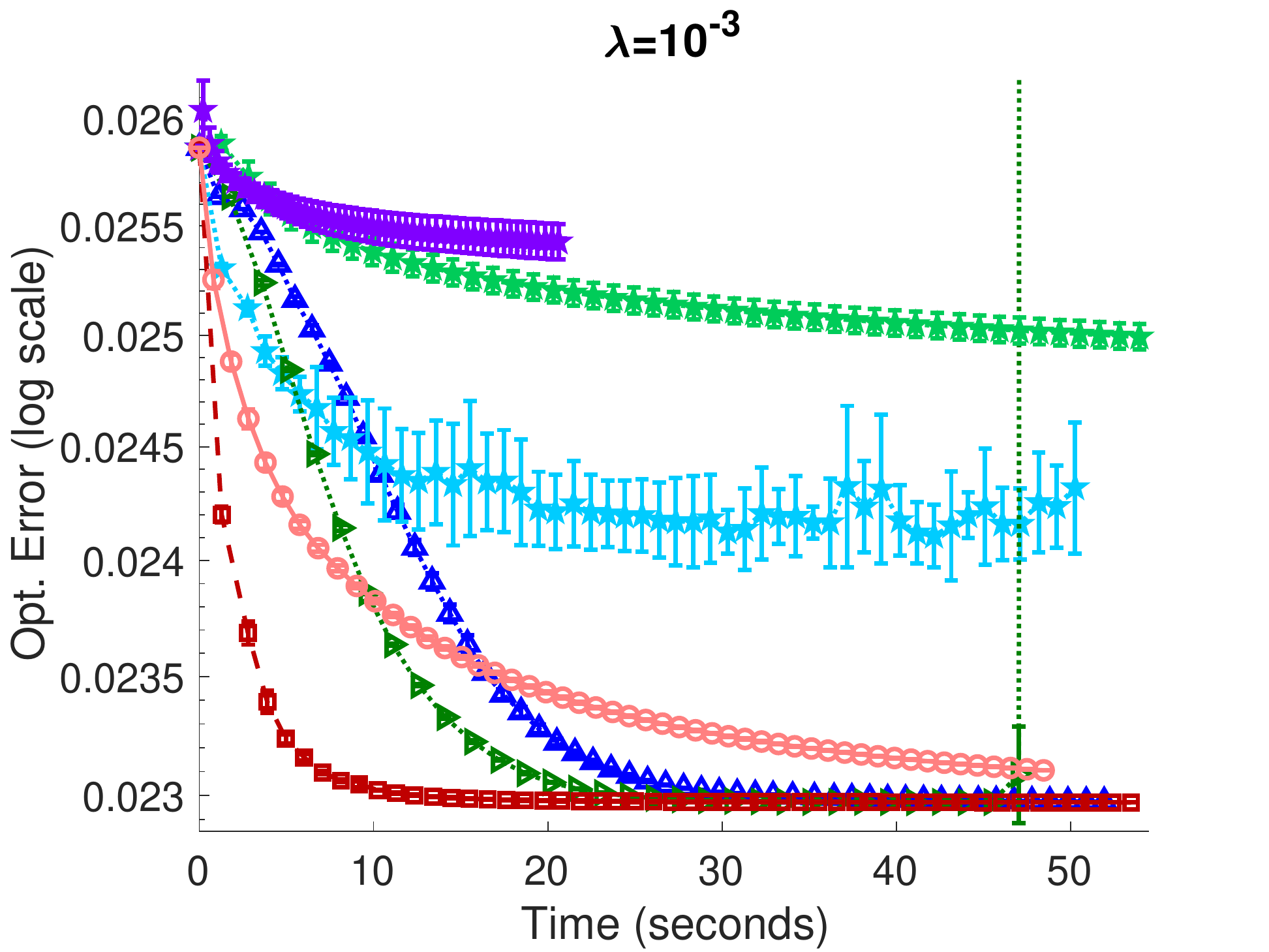}
		\includegraphics[width=.48\linewidth,height=5.5cm]{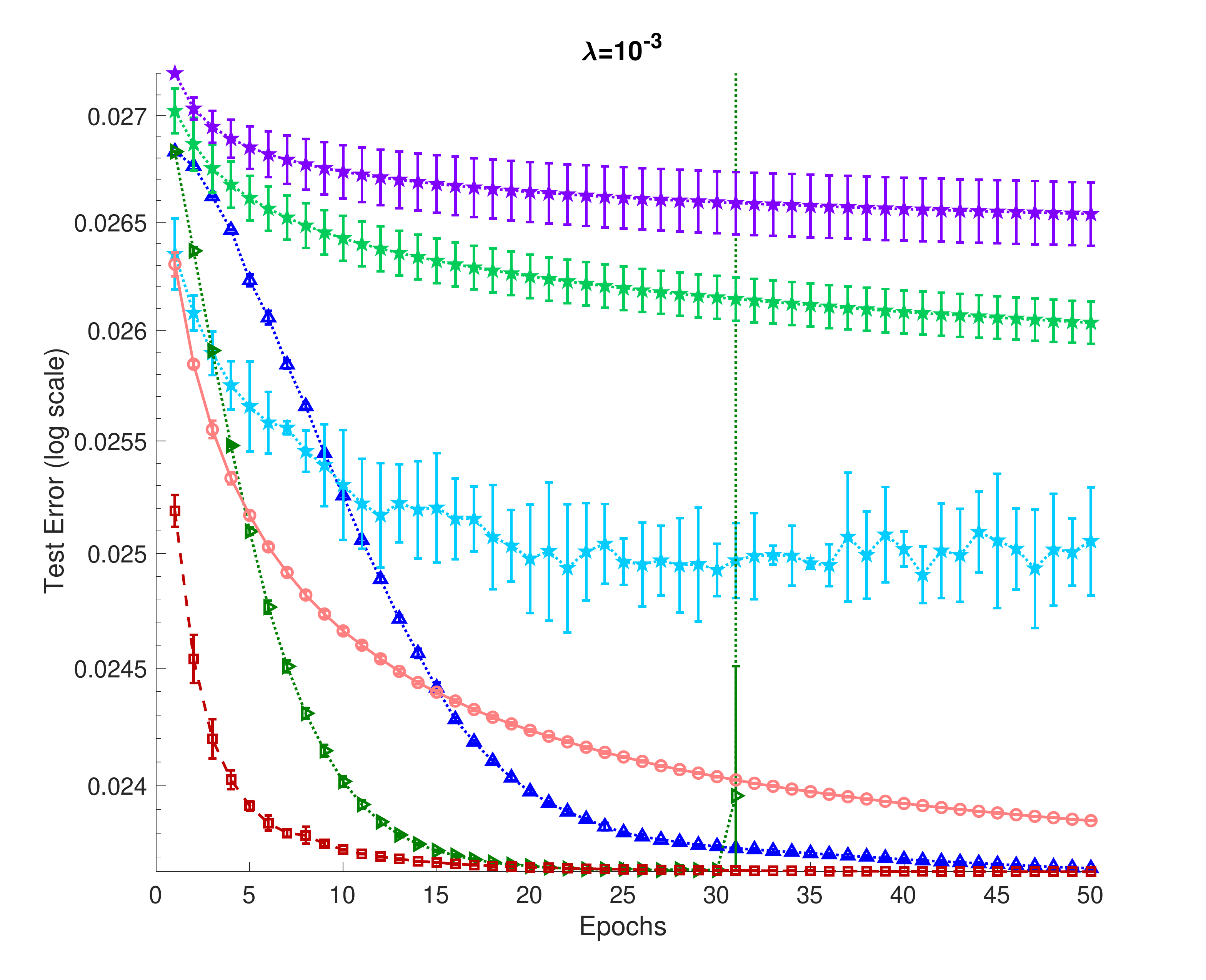}\\
		\caption{Comparision of results for $\ell_2$-svm loss function on \textit{w8a} dataset}
		\label{fig:my_labelSVMw8a}
	\end{figure}

	\clearpage
	
	\section{Quality of Hessian approximation}  \label{num:closenesstonewton}
\begin{figure}[h]
		\centering
		\begin{subfigure}[h]{0.45\textwidth}
			\centering
			\includegraphics[width=0.99\textwidth]{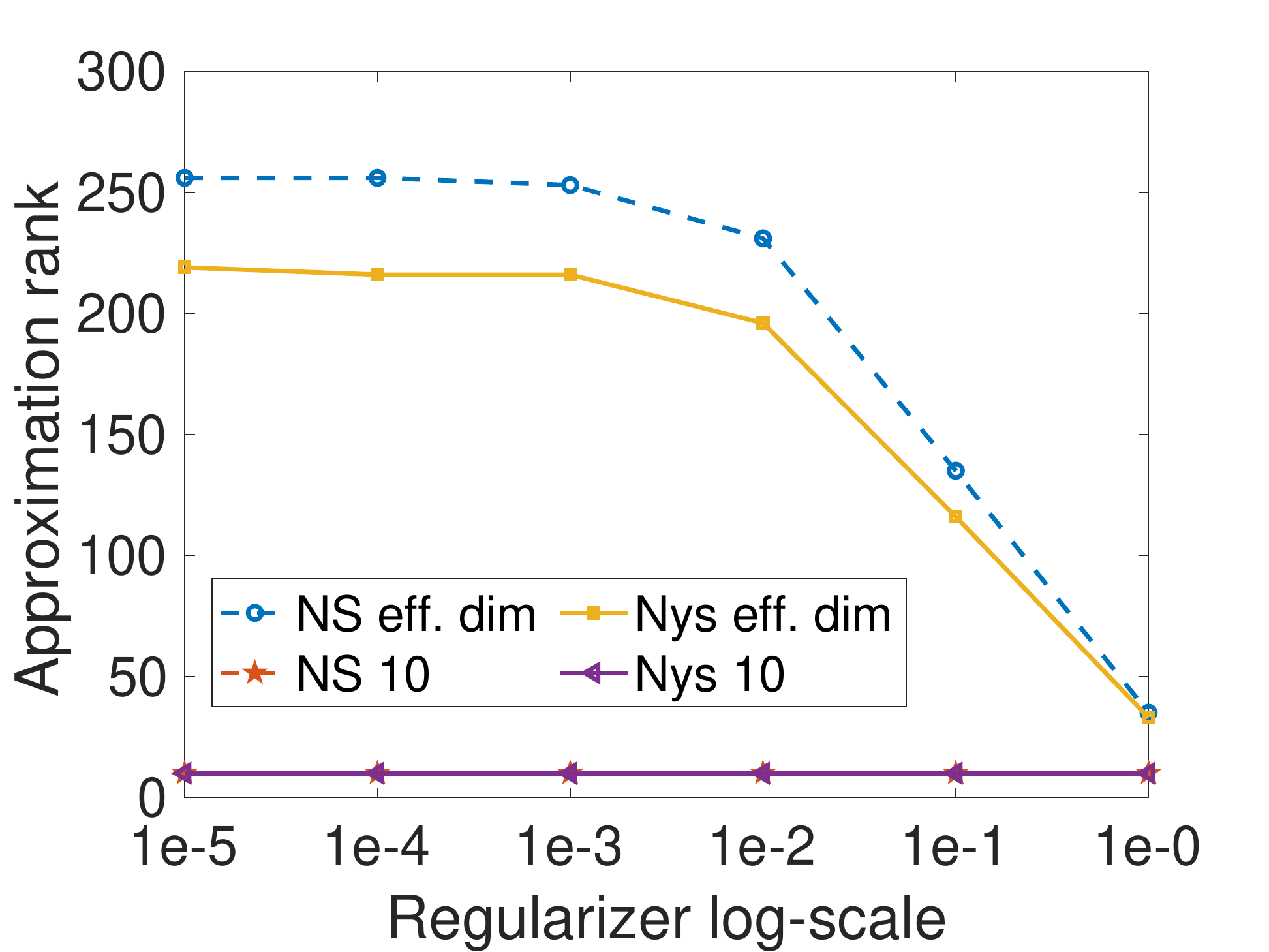}
			\caption{Rank of approximated Hessian with different $\lambda$.}
			\label{fig:rank_with_diff_lambda}
		\end{subfigure}
		\begin{subfigure}[h]{0.45\textwidth}
			\centering
			\includegraphics[width=0.99\textwidth]{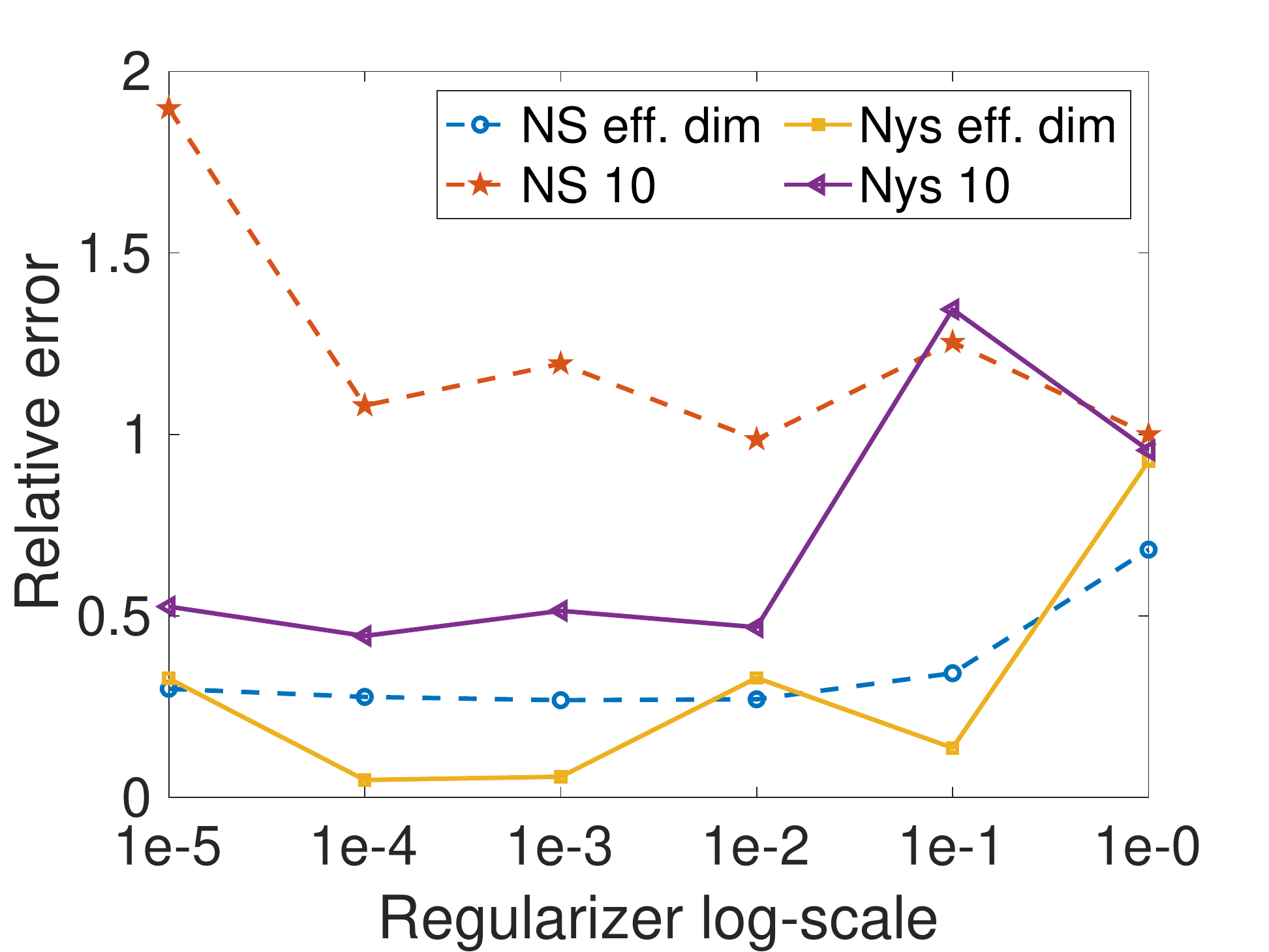}
			\caption{Relative error in Hessian approximation.}
			\label{fig:norm_diffrecne}
		\end{subfigure}
		\caption{Synthetic experimental results.}
		\vspace{-.15in}
	\end{figure}

\begin{figure*}[!h]
		\centering
		\includegraphics[width=.8\textwidth]{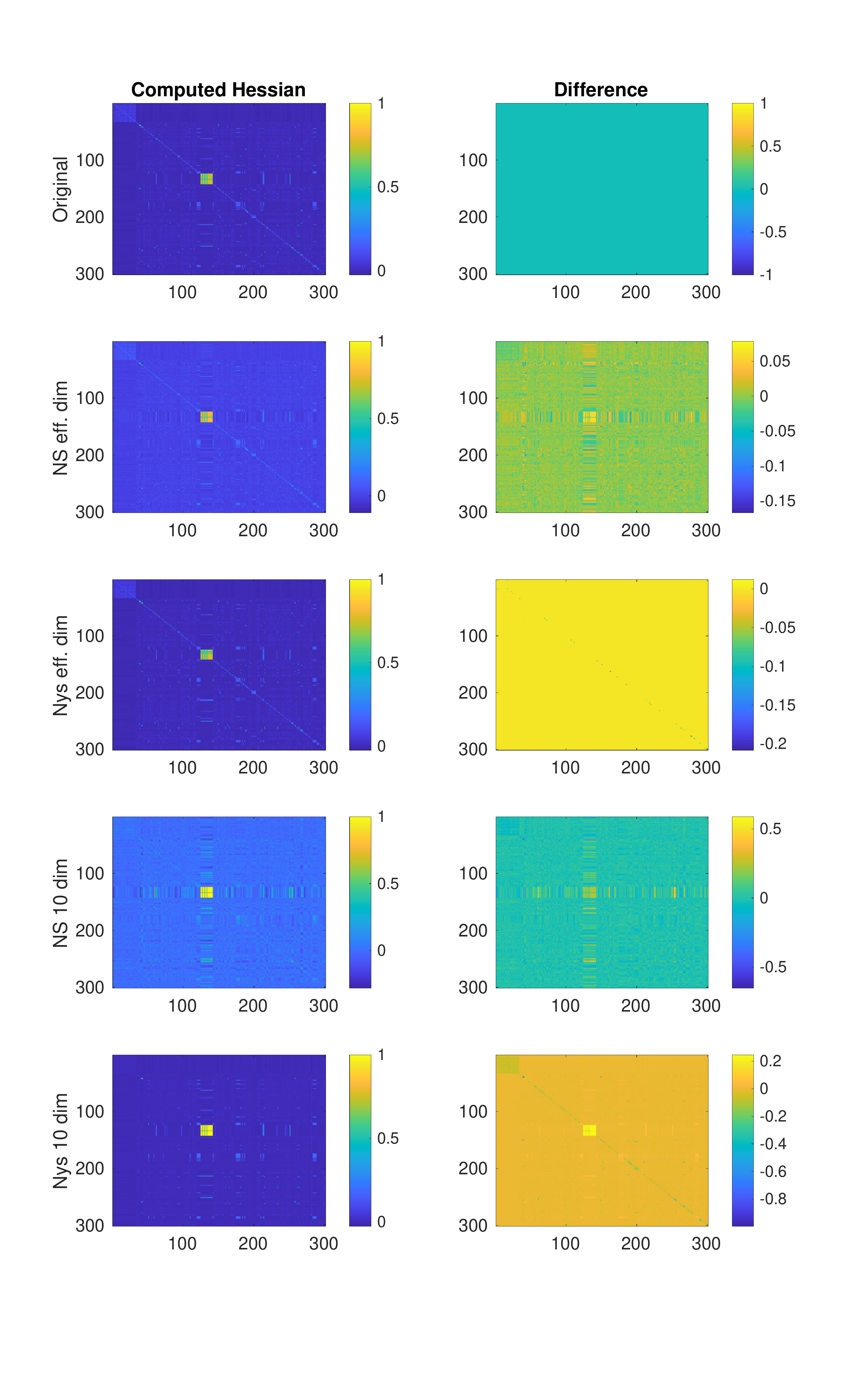}
		\caption{Closeness to the true Hessian(NS: Newton sketch, Nys: Nystr\"om)}
		\label{fig:closenessNewton}
	\end{figure*}

	In order to show the difference between the adaptive Newton sketch~\cite{Lacotte2021icml} and the Nystr\"om approximation, we compute their Hessian approximations for \textit{logistic regression} model on \textit{w8a} dataset with same $\boldw$. The adaptive Newton sketch uses the concept of effective dimension.
	\begin{equation}\label{eq:effective_dimension}
		d_{\lambda} = \text{trace}(\nabla^2 f(\boldw) (\nabla^2 f(\boldw)+ \lambda\boldI)^{-1}).
	\end{equation}
	First, we computed the effective dimensions for various $\lambda = [10^{-5},10^{-4}10^{-3}, 10^{-2},10^{-1}, 10^{0}]$. Next, we compute the adaptive Newton sketch 
	\begin{equation}
		\boldH_S =  (\nabla^2 f(\boldw) ^{1/2})^\top {\boldS^\top} \boldS (\nabla^2 f(\boldw) ^{1/2}),
	\end{equation}
	where $\boldS \in \mathbb{R}^{m\times n}$ and $\nabla^2 f(\boldw)^{1/2} \in \mathbb{R}^{n\times d}$. Since, their method is deterministic, to have a fair comparison, we compare the approximation only. 
	We computed the adaptive Newton sketch $\boldH_S$ with $m = \text{eff}(d_{\lambda})$ and  $\boldH_{S1}$ with $m = 10\ \text{eff}(d_{\lambda})$ using the randomized orthogonal matrix. Then, we computed the Nystr\"om $\boldN = \boldZ \boldZ^\top$(without using $\lambda$ and $\rho$) as in \eqref{def:Nystrom} with \textbf{$m = 10$ and with $m = d_{\lambda}$}. Figure~\eqref{fig:rank_with_diff_lambda} shows the behaviour with respect to the \textit{rank} of approximation, and Figure~\eqref{fig:norm_diffrecne} shows the relative error \emph{i.e.}, $\frac{\|\boldH-\text{App}(\boldH)\|}{\|\boldH\|}$, where $\boldH = \nabla^2 f(\boldw)$ is computed without using $\lambda$. Figure~\eqref{fig:closenessNewton} shows the approximations and its difference with original Hessian.
	
	It may be easily observed that the effective dimension can fluctuate with a larger $\lambda$ and becomes insipid to approximate the Hessian of ill-condition or high-dimensional problems. Moreover, Nystr\"om approximation showed that small column matrix with $m = 10$ is sufficient to approximate the Hessian with lowest error.

	\clearpage
	
\end{document}

%% file: mathdef.tex
\newcommand{\mathbbR}{\mathbb{R}}

\newcommand{\boldzero}{{\boldsymbol{0}}}

\newcommand{\boldA}{{\boldsymbol{A}}}
\newcommand{\boldB}{{\boldsymbol{B}}}
\newcommand{\boldC}{{\boldsymbol{C}}}
\newcommand{\boldD}{{\boldsymbol{D}}}

\newcommand{\boldG}{{\boldsymbol{G}}}
\newcommand{\boldH}{{\boldsymbol{H}}}
\newcommand{\boldI}{{\boldsymbol{I}}}
\newcommand{\boldJ}{{\boldsymbol{J}}}
\newcommand{\boldK}{{\boldsymbol{K}}}

\newcommand{\boldM}{{\boldsymbol{M}}}
\newcommand{\boldN}{{\boldsymbol{N}}}

\newcommand{\boldQ}{{\boldsymbol{Q}}}

\newcommand{\boldS}{{\boldsymbol{S}}}

\newcommand{\boldU}{{\boldsymbol{U}}}
\newcommand{\boldV}{{\boldsymbol{V}}}
\newcommand{\boldW}{{\boldsymbol{W}}}
\newcommand{\boldX}{{\boldsymbol{X}}}

\newcommand{\boldZ}{{\boldsymbol{Z}}}

\newcommand{\boldg}{{\boldsymbol{g}}}

\newcommand{\boldk}{{\boldsymbol{k}}}

\newcommand{\boldv}{{\boldsymbol{v}}}
\newcommand{\boldw}{{\boldsymbol{w}}}
\newcommand{\boldx}{{\boldsymbol{x}}}
\newcommand{\boldy}{{\boldsymbol{y}}}
\newcommand{\boldz}{{\boldsymbol{z}}}

\newcommand{\boldSigma}{{\boldsymbol{\Sigma}}}

\newcommand{\boldOmega}{{\boldsymbol{\Omega}}}

\newcommand{\calB}{{\mathcal{B}}}

\newcommand{\calK}{{\mathcal{K}}}
\newcommand{\calL}{{\mathcal{L}}}

\newcommand{\calW}{{\mathcal{W}}}






%% file: main_arxiv.bbl
\begin{thebibliography}{10}

\bibitem{AgarwalBH17jmlr}
Naman Agarwal, Brian Bullins, and Elad Hazan.
\newblock Second-order stochastic optimization for machine learning in linear
  time.
\newblock {\em Journal of Machine Learning Research, {JMLR}}, 18:116:1--116:40,
  2017.

\bibitem{amari1998natural}
Shun-Ichi Amari.
\newblock Natural gradient works efficiently in learning.
\newblock {\em Neural computation}, 10(2):251--276, 1998.

\bibitem{Bhatia2013matrix}
Rajendra Bhatia.
\newblock {\em Matrix Analysis}, volume 169.
\newblock Springer Science \& Business Media, 2013.

\bibitem{broyden1965class}
Charles~G Broyden.
\newblock A class of methods for solving nonlinear simultaneous equations.
\newblock {\em Mathematics of Computation}, 19(92):577--593, 1965.

\bibitem{broyden1969new}
Charles~G Broyden.
\newblock A new double-rank minimisation algorithm. preliminary report.
\newblock {\em Notices of the American Mathematical Society}, 16(4):670, 1969.

\bibitem{Byrd2016stochastic}
Richard~H Byrd, Samantha~L Hansen, Jorge Nocedal, and Yoram Singer.
\newblock A stochastic quasi-newton method for large-scale optimization.
\newblock {\em SIAM Journal on Optimization}, 26(2):1008--1031, 2016.

\bibitem{byrd1996analysis}
Richard~H Byrd, Humaid~Fayez Khalfan, and Robert~B Schnabel.
\newblock Analysis of a symmetric rank-one trust region method.
\newblock {\em SIAM Journal on Optimization}, 6(4):1025--1039, 1996.

\bibitem{chang2011libsvm}
Chih-Chung Chang and Chih-Jen Lin.
\newblock {LIBSVM}: A library for support vector machines.
\newblock {\em ACM Transactions on Intelligent Systems and Technology},
  2:27:1--27:27, 2011.

\bibitem{osti_4252678}
W~C Davidon.
\newblock Variable metric method for minimization.
\newblock 5 1959.

\bibitem{daxberger2021laplace}
Erik Daxberger, Agustinus Kristiadi, Alexander Immer, Runa Eschenhagen,
  Matthias Bauer, and Philipp Hennig.
\newblock Laplace redux-effortless bayesian deep learning.
\newblock {\em Advances in Neural Information Processing Systems, NeurIPS},
  2021.

\bibitem{Defazio2014saga}
Aaron Defazio, Francis Bach, and Simon Lacoste-Julien.
\newblock {SAGA:} {A} fast incremental gradient method with support for
  non-strongly convex composite objectives.
\newblock In {\em Advances in Neural Information Processing Systems, {NIPS}},
  pages 1646--1654, 2014.

\bibitem{drineas2005nystrom}
Petros Drineas and Michael~W. Mahoney.
\newblock On the nystr{\"{o}}m method for approximating a gram matrix for
  improved kernel-based learning.
\newblock {\em Journal of Machine Learning Research, {JMLR}}, 6:2153--2175,
  2005.

\bibitem{Duchi2011adaptive}
John Duchi, Elad Hazan, and Yoram Singer.
\newblock Adaptive subgradient methods for online learning and stochastic
  optimization.
\newblock {\em Journal of Machine Learning Research, {JMLR}}, 12(7), 2011.

\bibitem{fletcher1970new}
Roger Fletcher.
\newblock A new approach to variable metric algorithms.
\newblock {\em The Computer Journal}, 13(3):317--322, 1970.

\bibitem{fletcher1963rapidly}
Roger Fletcher and Michael~JD Powell.
\newblock A rapidly convergent descent method for minimization.
\newblock {\em The Computer Journal}, 6(2):163--168, 1963.

\bibitem{frangella2021randomized}
Zachary Frangella, Joel~A Tropp, and Madeleine Udell.
\newblock Randomized nystr$\backslash$" om preconditioning.
\newblock {\em arXiv preprint arXiv:2110.02820}, 2021.

\bibitem{goldfarb1970family}
Donald Goldfarb.
\newblock A family of variable-metric methods derived by variational means.
\newblock {\em Mathematics of Computation}, 24(109):23--26, 1970.

\bibitem{GrosseM16ICML}
Roger~B. Grosse and James Martens.
\newblock A kronecker-factored approximate fisher matrix for convolution
  layers.
\newblock In {\em Proceedings of the International Conference on Machine
  Learning, {ICML}}, pages 573--582, 2016.

\bibitem{he2016deep}
Kaiming He, Xiangyu Zhang, Shaoqing Ren, and Jian Sun.
\newblock Deep residual learning for image recognition.
\newblock In {\em Proceedings of the IEEE conference on computer vision and
  pattern recognition}, pages 770--778, 2016.

\bibitem{Hsieh2014icml}
Cho{-}Jui Hsieh, Si~Si, and Inderjit~S. Dhillon.
\newblock A divide-and-conquer solver for kernel support vector machines.
\newblock In {\em Proceedings of the International Conference on Machine
  Learning, {ICML}}, pages 566--574, 2014.

\bibitem{Johnson2013accelerating}
Rie Johnson and Tong Zhang.
\newblock Accelerating stochastic gradient descent using predictive variance
  reduction.
\newblock In {\em Advances in Neural Information Processing Systems, {NIPS}},
  pages 315--323, 2013.

\bibitem{karakida2020understanding}
Ryo Karakida and Kazuki Osawa.
\newblock Understanding approximate fisher information for fast convergence of
  natural gradient descent in wide neural networks.
\newblock {\em NeurIPS}, 2020.

\bibitem{Kasai2017jmlr}
Hiroyuki Kasai.
\newblock Sgdlibrary: {A} {MATLAB} library for stochastic optimization
  algorithms.
\newblock {\em Journal of Machine Learning Research, {JMLR}}, 18:215:1--215:5,
  2017.

\bibitem{Kingma2015adam}
Diederik~P. Kingma and Jimmy Ba.
\newblock Adam: {A} method for stochastic optimization.
\newblock In Yoshua Bengio and Yann LeCun, editors, {\em Proceedings of the
  International Conference on Learning Representations, {ICLR}}, 2015.

\bibitem{Kolte2015accelerating}
Ritesh Kolte, Murat Erdogdu, and Ayfer Ozgur.
\newblock Accelerating svrg via second-order information.
\newblock In {\em NIPS Workshop on Optimization for Machine Learning}, 2015.

\bibitem{Lacotte2021icml}
Jonathan Lacotte, Yifei Wang, and Mert Pilanci.
\newblock Adaptive newton sketch: Linear-time optimization with quadratic
  convergence and effective hessian dimensionality.
\newblock In Marina Meila and Tong Zhang, editors, {\em Proceedings of the 38th
  International Conference on Machine Learning, {ICML}}, volume 139, pages
  5926--5936, 2021.

\bibitem{Lan2019tnnls}
Liang Lan, Zhuang Wang, Shandian Zhe, Wei Cheng, Jun Wang, and Kai Zhang.
\newblock Scaling up kernel {SVM} on limited resources: {A} low-rank
  linearization approach.
\newblock {\em {IEEE} Transactions on Neural Networks and Learning Systems,
  {TNNLS}}, 30(2):369--378, 2019.

\bibitem{lin2021structured}
Wu~Lin, Frank Nielsen, Mohammad~Emtiyaz Khan, and Mark Schmidt.
\newblock Structured second-order methods via natural gradient descent.
\newblock {\em arXiv preprint arXiv:2107.10884}, 2021.

\bibitem{lin2021tractable}
Wu~Lin, Frank Nielsen, Mohammad~Emtiyaz Khan, and Mark Schmidt.
\newblock Tractable structured natural gradient descent using local
  parameterizations.
\newblock {\em ICML}, 2021.

\bibitem{Liu1989limited}
Dong~C Liu and Jorge Nocedal.
\newblock On the limited memory bfgs method for large scale optimization.
\newblock {\em Mathematical Programming}, 45(1):503--528, 1989.

\bibitem{mishkin2018slang}
Aaron Mishkin, Frederik Kunstner, Didrik Nielsen, Mark Schmidt, and
  Mohammad~Emtiyaz Khan.
\newblock Slang: fast structured covariance approximations for bayesian deep
  learning with natural gradient.
\newblock In {\em NeurIPS}, pages 6248--6258, 2018.

\bibitem{Moritz2016linearly}
Philipp Moritz, Robert Nishihara, and Michael Jordan.
\newblock A linearly-convergent stochastic l-bfgs algorithm.
\newblock In {\em Proceedings of the International Conference on Artificial
  Intelligence and Statistics, {AISTATS}}, pages 249--258, 2016.

\bibitem{Nguyen2017sarah}
Lam~M Nguyen, Jie Liu, Katya Scheinberg, and Martin Tak{\'a}{\v{c}}.
\newblock Sarah: A novel method for machine learning problems using stochastic
  recursive gradient.
\newblock In {\em Proceedings of the International Conference on Machine
  Learning, {ICML}}, pages 2613--2621, 2017.

\bibitem{Nocedal2006numerical}
Jorge Nocedal and Stephen Wright.
\newblock {\em Numerical optimization}.
\newblock Springer Science \& Business Media, 2006.

\bibitem{pilanci2017newton}
Mert Pilanci and Martin~J Wainwright.
\newblock Newton sketch: A near linear-time optimization algorithm with
  linear-quadratic convergence.
\newblock {\em SIAM Journal on Optimization}, 27(1):205--245, 2017.

\bibitem{Robbins1951sgd}
Herbert Robbins and Sutton Monro.
\newblock A stochastic approximation method.
\newblock {\em The Annals of Mathematical Statistics}, 22(3):400--407, 1951.

\bibitem{Schraudolph2007stochastic}
Nicol~N Schraudolph, Jin Yu, and Simon G{\"u}nter.
\newblock A stochastic quasi-newton method for online convex optimization.
\newblock In {\em Proceedings of the International Conference on Artificial
  Intelligence and Statistics {AISTATS}}, pages 436--443, 2007.

\bibitem{scieur2021generalization}
Damien Scieur, Lewis Liu, Thomas Pumir, and Nicolas Boumal.
\newblock Generalization of quasi-newton methods: application to robust
  symmetric multisecant updates.
\newblock In {\em Proceedings of the International Conference on Artificial
  Intelligence and Statistics, {AISTATS}}, pages 550--558, 2021.

\bibitem{shanno1970conditioning}
David~F Shanno.
\newblock Conditioning of quasi-newton methods for function minimization.
\newblock {\em Mathematics of computation}, 24(111):647--656, 1970.

\bibitem{talwalkar2010matrix}
Ameet Talwalkar.
\newblock {\em Matrix approximation for large-scale learning}.
\newblock PhD thesis, New York University, 2010.

\bibitem{tan2019efficientnet}
Mingxing Tan and Quoc Le.
\newblock Efficientnet: Rethinking model scaling for convolutional neural
  networks.
\newblock In {\em International Conference on Machine Learning (ICML)}, pages
  6105--6114, 2019.

\bibitem{tran2020bayesian}
M-N Tran, Nghia Nguyen, David Nott, and Robert Kohn.
\newblock Bayesian deep net glm and glmm.
\newblock {\em Journal of Computational and Graphical Statistics},
  29(1):97--113, 2020.

\bibitem{yang2020structured}
Minghan Yang, Dong Xu, Hongyu Chen, Zaiwen Wen, and Mengyun Chen.
\newblock Enhance curvature information by structured stochastic quasi-newton
  methods.
\newblock In {\em {IEEE} Conference on Computer Vision and Pattern Recognition,
  {CVPR}}, pages 10654--10663, 2021.

\bibitem{yao2021adahessian}
Zhewei Yao, Amir Gholami, Sheng Shen, Mustafa Mustafa, Kurt Keutzer, and
  Michael Mahoney.
\newblock Adahessian: An adaptive second order optimizer for machine learning.
\newblock In {\em Proceedings of the AAAI Conference on Artificial
  Intelligence}, pages 10665--10673, 2021.

\bibitem{Zhang2012aistats}
Kai Zhang, Liang Lan, Zhuang Wang, and Fabian Moerchen.
\newblock Scaling up kernel {SVM} on limited resources: {A} low-rank
  linearization approach.
\newblock In Neil~D. Lawrence and Mark~A. Girolami, editors, {\em Proceedings
  of the International Conference on Artificial Intelligence and Statistics,
  {AISTATS}}, pages 1425--1434, 2012.

\end{thebibliography}
